\documentclass[english,ruled]{article}
\usepackage[T1]{fontenc}
\usepackage[latin9]{inputenc}
\usepackage{verbatim}
\usepackage{algorithm2e}
\usepackage{amsmath}
\usepackage{amsthm}
\usepackage{amssymb}
\usepackage{graphicx}
\usepackage{xcolor}
\usepackage{multicol}
\usepackage{multirow}
\usepackage{geometry}
\usepackage{enumitem}
\usepackage{setspace}
\makeatletter
\usepackage[toc,page,header]{appendix}
\usepackage{minitoc}
\usepackage{ifthen}
\newboolean{doublecolumn}
\setboolean{doublecolumn}{false}
\newboolean{arxiv}
\setboolean{arxiv}{true}
\usepackage{pgfplots}
\pgfplotsset{compat=1.15}

\ifthenelse{\boolean{doublecolumn}}{
\renewcommand\vspace[1]{#1}
}

\renewcommand \thepart{}
\renewcommand \partname{}

\theoremstyle{plain}
\newtheorem{lem}{\protect\lemmaname}[section]
\theoremstyle{remark}
\newtheorem{rem}{\protect\remarkname}
\theoremstyle{plain}
\newtheorem{thm}{\protect\theoremname}[section]
\theoremstyle{plain}

\providecommand{\corollaryname}{Corollary}
\theoremstyle{plain}

\theoremstyle{plain}
\newtheorem{exple}{\protect\examplename}[section]

\usepackage{babel}
\providecommand{\propositionname}{Proposition}
\usepackage[colorlinks]{hyperref}       
\usepackage{url}            
\usepackage{booktabs}       
\usepackage{amsfonts}       
\usepackage{nicefrac}       
\usepackage{authblk}
\usepackage{optidef}

\global\long\def\norm#1{\lVert#1\rVert}%

\global\long\def\vertiii#1{\left\vert \kern-0.25ex  \left\vert \kern-0.25ex  \left\vert #1\right\vert \kern-0.25ex  \right\vert \kern-0.25ex  \right\vert }%

\global\long\def\diam{\mathrm{diam}}%

\global\long\def\argmin{\operatornamewithlimits{arg\,min}}%

\global\long\def\sign{\operatornamewithlimits{sign}}%

\global\long\def\dom{\mathrm{dom}}%

\global\long\def\and{\mathrm{and}}%

\global\long\def\vep{\varepsilon}%

\global\long\def\dom{\operatornamewithlimits{dom}}%

\global\long\def\Ebb{\mathbb{E}}%

\global\long\def\Rbb{\mathbb{R}}%

\global\long\def\Gcal{\mathcal{G}}%

\global\long\def\Ocal{\mathcal{O}}%

\global\long\def\Xcal{\mathcal{X}}%

\newcommand\condspace{{\ifthenelse{\boolean{doublecolumn}}{}{\\}}}
\newcommand\condsep{{\ifthenelse{\boolean{doublecolumn}}{\vspace{-\topsep}}{}}}
\newcommand\myfrac[2]{\frac{#1}{#2}}

\newcommand\tmop[1]{\mathrm{#1}}
\newcommand\tmtextbf[1]{\textbf{#1}}

\makeatother
\usepackage{babel}

\usepackage[numbers]{natbib}
\usepackage[resetlabels]{multibib}

\newcites{app}{References in appendix}

\providecommand{\lemmaname}{Lemma}
\providecommand{\remarkname}{Remark}
\providecommand{\theoremname}{Theorem}
\providecommand{\examplename}{Example}

\begin{document}

\global\long\def\sgm{\texttt{SGD}}%
\global\long\def\sgd{\texttt{SGD}}%
\global\long\def\spl{\texttt{SPL}}%
\global\long\def\rpc{\texttt{RPC}}%

\global\long\def\assign{\coloneqq}%
\global\long\def\tmstrong#1{{\bf #1}}%
\global\long\def\tmem#1{#1}%
\global\long\def\tmop#1{#1}%

\global\long\def\tmverbatim#1{#1}%
\global\long\def\revised#1{{#1}}%
\global\long\def\changed#1{\textcolor{blue}{#1}}%

\newtheorem{definition}{Definition}
\newcommand{\TODO}[1]{\textbf{\textcolor{red}{TODO: #1}}}

\font\myfont=cmr12 at 16.1pt
\title{\myfont {Stochastic Weakly Convex Optimization Beyond Lipschitz Continuity}}
\setlength{\parindent}{0pt}
\author[1]{Wenzhi Gao}
\author[2]{Qi Deng}
\affil[1]{\texttt{gwz@stanford.edu}, ICME, Stanford University}
\affil[2]{\texttt{contact.qdeng@gmail.com}, Shanghai University of Finance and Economics}

\maketitle

\begin{abstract}
This paper considers stochastic weakly convex optimization without the standard Lipschitz continuity assumption. Based on new adaptive regularization (stepsize) strategies, we show that a wide class of stochastic algorithms, including the stochastic subgradient method,  preserve the $\mathcal{O} ( 1 / \sqrt{K})$ convergence rate with constant failure rate.  Our analyses rest on rather weak assumptions: the Lipschitz parameter can be either bounded by a general growth function of $\norm{x}$ or locally estimated through independent random samples. 
Numerical experiments demonstrate the efficiency and robustness of our proposed stepsize policies. 
\end{abstract}

\section{Introduction \label{sec:intro}}
This paper studies the following stochastic optimization problem
\begin{eqnarray}
  \min_{x \in \mathbb{R}^n} & \psi (x) \assign f (x) + \omega (x) \assign
  \mathbb{E}_{\xi \sim \Xi} [f (x, \xi)] + \omega (x). & 
\end{eqnarray}
Here, $f(x, \xi)$ is a continuous function in $x$, with $\xi$ being a random sample drawn from a particular distribution $\Xi$. 
The function $\omega(x)$ is lower-semicontinuous, and its proximal mapping is easy to evaluate.
We assume both $f(x,\xi)$ and $\omega(x)$ are weakly convexfunctions. A function $g$ is defined as $\lambda$-weakly convex if $g + \frac{\lambda}{2} \|\cdot\|^2$ is convex, for some $\lambda\ge 0$.
 When $\lambda$ is unspecified,  $g$ is called weakly convex. Weak convexity has found many important applications, including phase retrieval, robust PCA,  reinforcement learning, and many others~\citep{duchi2019solving, charisopoulos2021low, wang2023policy}. And recent years witnessed a surge in interest regarding weakly convex optimization, leading to a substantial body of work on efficient algorithms with finite time complexity guarantees~\citep{davis2019stochastic, davis2018subgradient, deng2021minibatch, davis2019stochastic2, mai2020convergence}.
In particular, under the global Lipschitz continuity assumption, \citet{davis2019stochastic2} develop a model-based approach and analyzes the convergence of several stochastic algorithms under a unified framework. \\

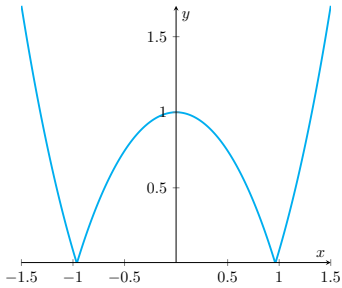
\begin{figure}[h]
\begin{center}
\scalebox{0.6}{
	\begin{tikzpicture}
    \begin{axis}[
            axis on top,
            legend pos=outer north east,
            axis lines = center,
            xlabel = $x$,
            ylabel = $y$,
        ]
        \addplot[very thick,cyan,samples=500,domain=-1.5:1.5] {abs(e^x + e^(-x) - 3)};
    \end{axis}
\end{tikzpicture}
}
\end{center}
\caption{$f(x) = |e^x  + e^{-x} - 3|$ exhibits exponential growth as $\|x\| \rightarrow + \infty$ \label{fig:fx}}
\end{figure}

While this global Lipschitz assumption is valid for many problems, such as piece-wise linear functions, it can be overly restrictive. To illustrate, consider the weakly convex function  $\psi(x) = |e^x + e^{-x} - 3|$ whose subgradient $\psi'(x)$ explodes exponentially as $\|x\|$ grows (\textbf{Figure \ref{fig:fx}}). 
Hence, in algorithm design, treating the Lipschitz constant as any fixed constant can lead to highly unstable iterations and, potentially, to the algorithm divergence.\\

To address this issue,  a straightforward strategy is to impose an explicit convex compact set constraint, such as $\{x:\|x\| \le B\}$ to address this issue. However, it introduces extra parameter tuning and may lead to a significantly overestimated Lipschitz constant. The latter phenomenon is evident in the toy example, where the Lipschitz constant grows exponentially with the domain set's diameter.
To deal with non-globally Lipschitz settings, one research direction is to shift from standard Euclidean geometry to Bregman divergence. \citet{lu2019relative} shows that when convex non-Lipschitz functions exhibit ``relative'' Lipschitz continuity under a carefully chosen divergence kernel, mirror descent still obtains the desired sublinear rate of convergence to optimality defined in the sense of Bregman divergence. 
However,  there are trade-offs to consider. Compared to {\sgd}, a mirror descent update is more expensive, often involving a nontrivial root-finding procedure. 
Additionally, choosing the right kernel is a nuanced and critical task, heavily reliant on an in-depth understanding of the subgradient's growth dynamics~(\citep{davis2018stochastic,zhang2018convergence}).\\

Alternatively, recent works aim to develop new algorithms/analyses under relaxed Lipschitz assumptions. For example, \citet{asi2019stochastic} show that using stochastic proximal point update, algorithmic dependency on the global Lipschitz constant can be relaxed to $\Ebb[\|f'(x^\star, \xi)\|^2]$, which only relies on the optimal solution $x^\star$. \citet{mai2021stability} show that, for stochastic convex optimization with quadratic growth, subgradient methods incorporating a clipping stepsize still ensures convergence, even if Lipschitz constant exhibits arbitrary growth.
In weakly convex optimization, \citet{li2023revisiting} shows that when the Moreau envelope of objective has a bounded level-set,  local Lipschitz continuity alone is sufficient to ensure convergence of the subgradient method. Nevertheless, extending their analysis to stochastic optimization remains challenging.
\citet{grimmer2019convergence} have established the convergence of \sgd{} in convex optimization,  without Lipschitz continuity. In \cite{zhu2023unified}, the authors propose a relaxed subgradient bound for weakly convex optimization:
\begin{equation}\label{eqn:subgradbound}
	\Ebb[\|f'(x, \xi)\|^2] \leq c_0 + c_1 \|x\|,\ c_0,c_1\ge 0,
\end{equation} 
which naturally induces a bound on the local Lipschitzness as a function of $\|x\|$.
Whether \sgd{} still converges in case of arbitrary non-Lipschitzness, especially those not conforming to the bounded assumption in~\eqref{eqn:subgradbound}, remains an open area of investigation.
 The major difficulty in analyzing stochastic optimization without the standard Lipschitz assumption stems from stability issues. In this paper, a stochastic algorithm is considered stable if it produces iterations in a bounded set with probability greater than $0$. Unlike the deterministic case, establishing stability in the face of randomness is not straightforward, especially when dealing with non-convex functions. This challenge motivates the exploration of an appropriate definition of non-Lipschitzness and the development of efficient algorithms for stochastic weakly convex optimization in this non-standard setting.

\subsection{Contributions}
We show that carefully chosen adaptive stepsizes can effectively deal with arbitrary non-Lipchitzness in stochastic model-based weakly convex optimization. Our contributions are as follows:
\begin{enumerate}[leftmargin=*, label={\arabic*)}]
\item When the Lipschitz constant is not uniformly bounded above but instead depends on a general growth function $\mathcal{G}(\cdot)$, we design a novel adaptive stepsize strategy such that stochastic weakly convex optimization achieves the $\mathcal{O}(1/\sqrt{K})$ convergence rate with a constant failure probability. Our analysis does not assume any specific form of $\mathcal{G}$, such as those implied by \eqref{eqn:subgradbound}. To our knowledge, this is the first result of stochastic weakly convex optimization for arbitrary non-Lipschitz objectives. Our analysis applies to a broad class of model-based algorithms~\cite{davis2019stochastic, deng2021minibatch}, including \sgd{} as a special case.
\item Even if the growth function $\mathcal{G}$ is unknown, we show that it is still possible to achieve the same convergence guarantee. To this end, we introduce a new adaptive stepsize based on the concept of ``reference Lipschitz continuity'', which allows us to estimate the Lipschitz parameter of a stochastic model function using local samples. Our algorithm is highly flexible, and it can be applied to most of the weakly convex problems of interest. Moreover, our analyses can be extended to solving convex stochastic optimization without Lipschitz continuity. A more detailed discussion is left to \textbf{Section \ref{app:cvx}}.
\end{enumerate}

\paragraph{Other related works}
Adaptive stepsize strategy and gradient clipping are two important tools adopted in our algorithm framework. On the one hand, stepsize selection has been an important topic in stochastic optimization, and it has been justified that adaptive stepsize benefits stochastic first-order methods both in theory and in practice \cite{duchi2011adaptive, kingma2014adam, li2019convergence, hinton2012neural, pmlr-v202-defazio23a, ivgi2023dog, malitsky2023adaptive}.  On the other hand, gradient clipping {\cite{zhang2019gradient}} will be employed as a technique in the paper. In theory, gradient clipping was initially identified as a tool to solve problems with generalized Lipschitz smoothness condition~\citep{li2023convex, xie2023trust, zhang2019gradient}. Recent works \citep{gorbunov2020stochastic, koloskova2023revisiting} show that gradient clipping can effectively deal with problems with heavy-tail noise. It is also observed that gradient clipping improves robustness and stability of {\sgd} {\cite{mai2020convergence}} in stochastic convex optimization.  In our analysis, a generalized version of gradient clipping is developed to  alleviate the instability arising from stochastic noise.
\section{Preliminaries \label{sec:preliminary}}

\paragraph{Notations} Throughout the paper $\| \cdot \|$ and $\langle \cdot, \cdot \rangle$
denote the Euclidean inner product and norm. Subdifferential of $f$ is given
by $\partial f (x) \assign \{ v : f (y) \geq f (x) + \langle v, y - x \rangle
+ o (\| x - y \|), y \rightarrow x \}$ and $f' (x) \in \partial f (x)$ is called
a subgradient. A growth function $\mathcal{G}(\cdot)$ is a continuous non-decreasing function mapping from $\Rbb_+$ to $\Rbb_+$.

\paragraph{Envelope smoothing}

For weakly convex optimization, our analysis adopts the Moreau envelope as the potential function. Let $f$ be a $\lambda$-weakly convex function. Given $\rho > \lambda$, the Moreau envelope and the associated proximal mapping of $f$
are given by
\[f_{1/\rho}(x)\assign\min_{y}\big\{ f(y)+\tfrac{{\rho}}{2}\|x-y\|^{2}\big\} \quad \text{and} \quad\text{prox}_{f/\rho}(x)\assign\argmin_{y}\big\{ f(y)+\tfrac{{\rho}}{2}\|x-y\|^{2}\big\}.\]

The Moreau envelope can be interpreted as a smooth approximation
of the original function. $f_{1/\rho}(x)$ is differentiable and its gradient is $\nabla f_{1/\rho}(x)=\rho(x-\text{prox}_{f/\rho}(x))$. If $\norm{\nabla f_{1/\rho}(x)}\le\vep$, then $x$
is in the proximity  of a near stationary point $\text{prox}_{f/\rho}(x)$ \cite{davis2019stochastic}. An important observation is that the existence of $f_{1/\rho}(x)$ only relies on weak convexity instead of Lipschitz continuity of $f$. We denote $\hat{x} := \text{prox}_{f/\rho}(x)$ in the rest of the paper.

\paragraph{Model-based optimization} 
Our main algorithm will be presented in a ``model-based'' fashion, which encompasses several first-order methods, including the most widely used (proximal) subgradient method. Model-based optimization \cite{davis2019stochastic} contains two components: a stochastic model function and a stepsize (regularization) parameter. In each iteration, we can construct a local approximation of $f(x)$ based on random sample $\xi^k$ and the current iterate $x^k$. The stochastic function, denoted by $f_{x^k}(\cdot, \xi^k) + \omega(x)$, is called a model function. Then we take parameter $\gamma_k$ and minimize this local approximation under quadratic regularization $\frac{\gamma_k}{2} \|x - x^k\|^2$ to obtain the next iterate $x^{k+1}$. Typical models include 

\begin{itemize}[leftmargin=20pt]
	\item (Sub)gradient. $f_x (y, \xi) = f (x, \xi) + \langle f' (x,  \xi), y - x \rangle, \Ebb[f'(x, \xi)] = f'(x) \in \partial f(x)$
	\item {Prox-linear}~\citep{drusvyatskiy2019efficiency}. $f(x, \xi) = h(c(x, \xi))$ and $f_x (y, \xi) = h (c (x, \xi) + \langle \nabla c (x, \xi),
  y - x \rangle)$
  \item {Truncated}~\citep{asi2019stochastic}.
 $f_x (y, \xi) = \max \{ f (x, \xi) + \langle \nabla f (x,
  \xi), y - x \rangle, \ell \}$, where $\ell$ is a known lower-bound of model 
\end{itemize}

\textbf{Algorithm \ref{alg:1}} summarizes model-based optimization.\\

\begin{algorithm}[h]\label{alg:1}
{\textbf{Input $x^{1}$}}\\
\For{k =\rm{ 1, 2,...}}{
Sample data $\xi^k$ and choose the proximal regularizer weight $\gamma_k>0$
\begin{align}
x^{k + 1} = \argmin_x \, \Big\{ f_{x^k}(x, \xi^k)+ \omega (x) + \myfrac{\gamma_k}{2} \|x - x^k\|^2 \Big\}. \label{eq:update-model-based}
\end{align} 
}
\caption{Stochastic model-based optimization\label{alg:model-based}}
\end{algorithm}

We see that \textbf{1)}  model function $f_{x}(\cdot, \xi)$; \textbf{2)} regularization parameter $\gamma_k$ are two core components for our algorithm design. Throughout this paper, we show how properly chosen $\gamma_k$ improves convergence beyond Lipschitz continuity. We start by making assumptions.

\paragraph{Assumptions.} We make the following assumptions across the paper.
\begin{enumerate}[leftmargin=35pt,itemsep=2pt,label=\textbf{A\arabic*:},ref=\rm{\textbf{A\arabic*}}]
  \item It is possible to generate i.i.d. samples $\{ \xi^k \}$. \label{A1}
  
  \item $\omega (x)$ is $\kappa$-weakly convex and
  $L_{\omega}$-Lipschitz continuous for all $x \in \dom \omega$.
  \label{A2}
  \item $\mathbb{E}_{\xi} [f_x (x, \xi)] = f (x)$ for all $x \in \dom
  \omega$ and $\mathbb{E}_{\xi} [f_x (y, \xi) - f (y)] \leq \frac{\tau}{2} \| x - y
  \|^2$ for all $x, y \in \dom \omega$. Moreover, $f_x (y, \xi)$ is convex for all $x, y \in \dom \omega$ and $\xi
  \sim \Xi$.
  \label{A5}
  \item $\psi_{1 / \rho} (x) = \min_z \{ \psi (z) + \frac{\rho}{2} \| z -
  x \|^2 \}$ is lower bounded by $- \Lambda \leq 0$.
  \label{A6}
  \item The $v$-level-set of $\psi_{1 / \rho} (x)$: $\mathcal{L}_v = \{ x :
  \psi_{1 / \rho} (x) \leq v \}$ has a bounded diameter $\text{diam}
  (\mathcal{L}_v) \leq \mathsf{B}_v < \infty$.\label{A7}
\end{enumerate}

\begin{rem}
The first four assumptions are general for weakly convex optimization. \ref{A7} is satisfied when $\psi$ is coercive \cite{li2023revisiting}, and it will be invoked to deal with non-Lipschitzness. Although \ref{A7} excludes some special cases such as interpolation, often these cases can be handled by exploiting different analyses~\citep{li2023revisiting}. In \ref{A5}, while we adopt a convex model function to simplify the analysis, it can be extended to weakly convex model functions.
\end{rem}

\begin{rem}
	Given \ref{A1} to \ref{A5}, it follows  that $f(x)$ is $\tau$-weakly convex and $
	\psi(x)$ is $(\tau + \kappa)$ weakly convex. 
\end{rem}

\paragraph{Structure of the paper} \textbf{Section
\ref{sec:known-lip}} discusses the convergence of weakly convex optimization with the standard Lipschitz condition, which serves as  a benchmark to provide sufficient intuition for the algorithm design in more challenging scenarios. \textbf{Section
\ref{sec:known-nonlip}} and \textbf{\ref{sec:unknown-nolip}} discuss two
cases where standard Lipschitz condition fails to hold. {\tmstrong{Section \ref{sec:exp}}}
conducts numerical experiments to verify our results.
\section{Weakly convex optimization under standard Lipschitzness}\label{sec:known-lip}

In this section, we consider the standard case as a benchmark and assume that
\ifthenelse{\boolean{doublecolumn}}{
\begin{enumerate}[leftmargin=25pt,itemsep=2pt,label=\textbf{B\arabic*:},ref=\rm{\textbf{B\arabic*}}]
}{
\begin{enumerate}[leftmargin=35pt,itemsep=2pt,label=\textbf{B\arabic*:},ref=\rm{\textbf{B\arabic*}}]
}
  	\item $f_x (x, \xi) - f_x (y, \xi) \leq L_f (\xi) \| x - y \|$ for all $x, y$ and $\xi \sim \Xi$. $\mathbb{E} [L_f (\xi)^2] \leq L_f^2$. \label{B1}
\end{enumerate}
This assumption is common in nonsmooth optimization, and the following descent property is known.
\begin{lem}
  \label{lem:1}Suppose that \ref{A1} to \ref{A5} as well as
  \ref{B1} holds, then given $\rho > \kappa  + \tau, \gamma_k > \rho$, 
  $\mathbb{E}_k [\psi_{1 / \rho} (x^{k + 1})] \leq \psi_{1 / \rho} (x^k) -\tfrac{\rho (\rho - \tau - \kappa)}{2
  (\gamma_k - \kappa)} \| \hat{x}^k - x^k \|^2 + \tfrac{2 \rho
  L_f^2}{(\gamma_k - \rho) (\gamma_k - \kappa)},$
  where $\mathbb{E}_k [\cdot] \assign \mathbb{E} [\cdot | \xi^1, \ldots,
  \xi^{k}]$ denotes the conditional expectation taken over $\xi^1, \ldots,
  \xi^{k}$.
\end{lem}
Taking into account that $\gamma_k$ is generally taken to be much larger than the other constants in the algorithm,  \textbf{Lemma \ref{lem:1}} reveals the following relation
\begin{equation}
	\mathbb{E}_k [\psi_{1 / \rho} (x^{k + 1})] \leq \psi_{1 / \rho} (x^k) - \mathcal{O}(\gamma_k^{-1}) \| \nabla \psi_{1/\rho} (x^k) \|^2 +  \mathcal{O}(L_f^2\gamma_k^{-2}), \label{eqn:recursion}
\end{equation}
where $\mathcal{O}(L_f^2\gamma_k^{-2})$ characterizes the error from both stochastic noise and nonsmoothness. Taking $\gamma_k \equiv \mathcal{O}(\sqrt{K})$ and telescoping the relation of \textbf{Lemma \ref{lem:1}}, we get the following convergence result.
\begin{thm}
  \label{thm:1}Under the same conditions as \textbf{Lemma \ref{lem:1}}, if
  we take $\gamma_k \equiv \rho + \kappa + \alpha \sqrt{K}$, then we have
  $\min_{1 \leq k \leq K} \mathbb{E} [\| \nabla \psi_{1 / \rho} (x^{k}) \|^2] \leq \tfrac{2
     \rho}{\rho - \tau - \kappa} \big[ \tfrac{ \rho D}{K} +
     \tfrac{1}{\sqrt{K}} \big( \alpha D + \tfrac{2\rho L^2_f}{\alpha} \big)
     \big],$
where $D = \psi (x^1)  - \inf_x \psi (x)$.
\end{thm}

\textbf{Theorem \ref{thm:1}} is standard in the literature \cite{davis2019stochastic, deng2021minibatch}. One important intuition we want to establish is that, the choice of $\gamma_k \equiv \mathcal{O}(\sqrt{K})$ is a consequence of the following trade-off: suppose we telescope over \eqref{eqn:recursion} directly, then 
\[ \frac{1}{K}\sum_{k=1}^K \mathcal{O}(\gamma_k^{-1}) \Ebb [\|\nabla\psi_{1/\rho}(x^k) \|^2] \leq \mathcal{O}\Big(\frac{1}{K}\Big) + \frac{1}{K}\sum_{k=1}^K \mathcal{O}(L_f^2\gamma_k^{-2}). \]
First, we need large, in other words, conservative $\gamma_k$, such that the error of potential reduction $\mathcal{O}(\sum_k L_f^2\gamma_k^{-2})$ is properly bounded. But the cost of being conservative is the reduced amount of reduction $\mathcal{O}(\gamma_k^{-1}) \|\nabla\psi_{1/\rho}(x^k) \|^2$. Now that due to \ref{B1}, $L_f$ is exactly bounded by constant. Finally, the optimal trade-off finally gives the choice of $\gamma_k \equiv \mathcal{O}(\sqrt{K})$, and $\mathcal{O}({1/\sqrt{K}})$ rate of convergence. As we will show in the following sections, in face of non-Lipschitzness, the error $\mathcal{O}(L_f^2\gamma_k^{-2})$ cannot be bounded by choosing some constant $\gamma_k$. What we do is find suitably large $\gamma_k$, adaptively, to reduce this error. Using \ref{A7} and a probabilistic analysis, we manage to do this without compromising the convergence rate.
\section{Weakly convex optimization under generalized Lipschitzness}\label{sec:known-nonlip}

Before achieving our goal of solving non-Lipschitz weakly convex optimization problems, we start from a less challenging case of non-Lipschitzness characterized as follows.
\begin{enumerate}[leftmargin=35pt,itemsep=2pt,label=\textbf{C\arabic*:},ref=\rm{\textbf{C\arabic*}}]
  \item $f_x (y, \xi) - f_x (z, \xi) \leq L_f (\xi) \mathcal{G} (\| x \|) \| y
  - z \|$ for all $x, y, z$; $\xi \sim \Xi$ and $\mathbb{E} [L_f (\xi)^2] \leq
  L_f^2$; Recall that growth function $\mathcal{G}$ is monotonically increasing. \label{C1}
\end{enumerate}

This assumption implies that our model function is globally Lipschitz, but the Lipschitz constant has a known dependency on the norm of the expansion point $x$. Our analysis applies if we can estimate an upper bound of $\mathcal{G}$. But for the brevity of analysis, we take this upper bound to be $\mathcal{G}$ itself. Many real-life applications have this structure, especially if the source
of non-Lipschitzness is a high-order polynomial.
\begin{exple}[Phase retrieval]
  Consider  $f (x, \xi) = | \langle a, x
  \rangle^2 - b |, a \in \mathbb{R}^n, b \in \mathbb{R}_+$. The subgradient
  model $f_x (y, \xi) = \langle f' (x, \xi), y - x \rangle = \langle 2 \cdot
  \sign (\langle a, x \rangle - b) \langle a, x \rangle a, y - x
  \rangle$ satisfies
  \[ f_x (y, \xi) - f_x (z, \xi) \leq 2 \| a \|^2 \| x \| \cdot \| y - z \|
     = L_f (\xi)\, \mathcal{G} (\| x \|) \,\| y- z \|, \]
  where $L_f (\xi) = 2 \| a \|^2, \mathcal{G} (\| x \|) = \| x \|$.
\end{exple}
\begin{exple}[Subgradient method] \sgd{} corresponds to the model function $f_x(y,\xi)=f(x,\xi)+ \langle f'(x,\xi), y-x \rangle$. Then, \ref{C1} is satisfied when $\|f'(x,\xi)\|\le L_f(\xi) \Gcal(\|x\|)$. It follows that $\Ebb[\|f'(x,\xi)\|^2]\le L_f^2\Gcal^2(\|x\|)$. One can readily see that relation \eqref{eqn:subgradbound} corresponds to the special case of $\Gcal^2(\cdot)$ being a linear function.  
\end{exple}

One direct consequence of \ref{C1} is that the Lipschitz constant of $f_x(\cdot, \xi)$
can go to $\infty$ when $\| x \| \rightarrow \infty$. Moreover, we cannot directly rely on Lipschitzness of $f(x)$. Taking subgradient update as an example, this implies $f' (x, \xi)$ has a large norm, leading to
a higher chance of divergence. From the perspective of convergence analysis, the error term $\mathcal{O}(L_f^2\gamma_k^{-2})$ in \eqref{eqn:recursion} becomes
hard to bound, and \textbf{Lemma \ref{lem:2}} quantifies this hardness.
\begin{lem}   \label{lem:2}
Suppose \ref{A1} to \ref{A5} as well as \ref{C1} holds, then given $\rho > \kappa + \tau, \gamma_k > \rho$,
  \[ \mathbb{E}_k [\psi_{1 / \rho} (x^{k + 1})] \leq \psi_{1 / \rho} (x^k) -
     \frac{\rho (\rho - \tau - \kappa)}{2(\gamma_k - \kappa)} \| \hat{x}^k - x^k \|^2 +
     \frac{\rho (\mathcal{G} (\| x^k \|) L_f + L_{\omega})^2}{2 \gamma_k(\gamma_k - \kappa)},
  \]
  where $\gamma_k$ is independent of $\xi^k$.
\end{lem}

According to \textbf{Lemma \ref{lem:2}}, the error of potential reduction involves the norm of 
the current iteration, which makes the previous constant stepsize analysis invalid, since we cannot assume $\mathcal{G}(\|x^k\|)$ is bounded. As a natural fix, we can take  \[\gamma_k = \mathcal{O} \big( \mathcal{G} (\| x^k \|) \sqrt{K} \big)\] to reduce the error. However, according to the trade-off we previously mentioned, unless $\mathcal{G} (\| x^k \|)$ is bounded by some constant independent of $K$, the reduction in the potential function can be arbitrarily small, and we still cannot obtain a valid convergence rate. To resolve this issue, we essentially need to show the boundedness of $\{\|x^k\|\}$, and our solution is to associate boundedness of $\{\|x^k\|\}$ with another bounded quantity during the algorithm: $\Ebb[\psi_{1 / \rho}(x^k)]$. Intuitively, since $\Ebb[\psi_{1 / \rho}(x^k)]$ is reduced during the algorithm, it remains bounded on expectation, and  from \ref{A7} we know that boundedness of $\{\psi_{1 / \rho}(x^k)\}$ implies boundedness of $\|x^k\|$, giving bounded $\mathcal{G} (\| x^k \|)$ and the $\mathcal{O}(1/\sqrt{K})$ rate we want. The following asymptotic result confirms our intuition.

\begin{thm}
  \label{thm:2}Under the same conditions as \textbf{Lemma \ref{lem:2}} and \ref{A6}, \ref{A7}, if $\gamma_k = \rho + \kappa + (\mathcal{G} (\| x^k \|) + 1) k^{\zeta}, \zeta
  \in ( \frac{1}{2}, 1 )$, then as $k\rightarrow \infty$, $\{ \| x^k \| \}$ is bounded with
  probability $1$; Moreover, the sequence $\{\inf_{j \leq k} \| \nabla \psi_{1 / \rho} (x^j) \|\}$ converges to $0$ almost surely.
\end{thm}

While it is relatively easy to show convergence asymptotically, to obtain a convergence rate, we need a more careful analysis of the algorithm behavior. One major difficulty is that boundedness of $\Ebb[\psi_{1 / \rho}(x^k)]$ is unable directly to give us information of $\|x^k\|$, since this relation holds only on expectation. To deal with this issue, we resort to probabilistic tools. And we establish a new probabilistic argument in the following subsection.

\subsection{Stability of the iterations}
In this subsection, we aim to analyze the stability of the iterates of a stochastic algorithm on a non-Lipschitz function. The intuition is very simple: if a stochastic algorithm reduces some potential function with a bounded level-set, then the iterates will stay in a bounded region with high probability. We provide the basic sketch of proof and leave a more rigorous argument to the appendix. \\

Our analysis relies on two simple facts that we gain from the adaptive stepsize.
\begin{lem}[Informal] \label{lem:2.22}
Under the same conditions as \textbf{Lemma \ref{lem:2.2}}, if $\gamma_k = \mathcal{O} \big( (\mathcal{G} (\| x^k \|) + 1) \sqrt{K} \big)$, then we have, for all $k = 2,\ldots, K$, that
\begin{align}
	\mathbb{E} [\| x^{k + 1} - x^k \|] \leq{} & \mathcal{O}(\tfrac{1}{\sqrt{K}}) \label{eqn:intuition-1}\\
	\mathbb{E} [\psi_{1 / \rho} (x^k)] \leq{} &\mathcal{O}(1) \label{eqn:intuition-2}
\end{align}
\end{lem}

Relation \eqref{eqn:intuition-1} says that with our adaptive stepsize strategy, we are very careful exploring the feasible region and at each iteration we only take a tiny step of $\mathcal{O}(1/\sqrt{K})$. The second relation \eqref{eqn:intuition-2} comes directly from \textbf{Lemma~\ref{lem:2}}.  Indeed, with $\gamma_k = \mathcal{O} \big( (\mathcal{G} (\| x^k \|) + 1) \sqrt{K} \big)$ we have 
$\tfrac{\rho (\mathcal{G} (\| x^k \|) L_f + L_{\omega})^2}{2 \gamma_k(\gamma_k - \kappa)} = \mathcal{O}(1/K)$.
And if we telescope \textbf{Lemma \ref{lem:2}} and take expectation over all the randomness for $k = 2,\ldots K$, we have 
\[\textstyle\mathbb{E} [\psi_{1 / \rho} (x^k)] \leq \sum_{j=1}^k \mathcal{O}(1/K)=\mathcal{O}(1), ~~~k = 2,\ldots, K.\]

Each of the two relations alone may not tell us useful information, since they both hold on expectation. However, when they are combined, a more useful result is available. Our argument is as follows:\\

Consider the event ``$\|x^k\| \text{ is large}$'' and we wish to upper-bound its probability. We have the following facts:

\begin{enumerate}
  
  \item If $\| x^k \| $ is large,  $\| x^{k + 1} \| \geq \|x^k\|
  -\mathcal{O} ( 1 / \sqrt{K} )$ by triangle inequality and $\| x^{k + 1} \|$ is also large.

   \item If $\| x^{k+1} \|$ is large, then $\psi_{1 / \rho} (x^{k + 1})$ is large by \ref{A7}.
   \item $\mathbb{E} [\psi_{1 / \rho} (x^{k+1})]$ is bounded by some constant.
\end{enumerate}

In other words, conditioned on the event ``$\|x^k\| \text{ is large}$'', to ensure that $\mathbb{E} [\psi_{1 / \rho} (x^{k + 1})]$ is still bounded, either \textbf{case 1)} the event happens with low probability, or \textbf{case 2)}  $x^{k+1}$ has to immediately jump back to a bounded region of smaller radius. However, since our adaptive stepsize restricts the ``jump'' between two consecutive iterations, \textbf{case 2)} cannot happen. Therefore, it is unlikely that ``$\|x^k\| \text{ is large}$''.\\

This argument brings us the following tail-bound characterizing the behavior of $\|x^k\|$ as a random variable.
\begin{lem}
  \label{lem:2.2}Under the same conditions as \textbf{Lemma \ref{lem:2}} as well as \ref{A6}, \ref{A7},
  if we take $\gamma_k = \rho +  \kappa + \tau + \alpha (\mathcal{G}(\| x^k \|) + 1) \sqrt{K}$, then the tail bound
  \[ \mathbb{P} \Big\{ \| x^k \| \geq \mathsf{B}_{a \Delta} + \frac{4(L_f +
     L_{\omega})}{\alpha \sqrt{K}} \Big\} \leq \frac{2\Delta}{a \Delta +
     \Lambda}, \]
  holds for all $2 \leq k \leq K$, where $\Delta = \psi_{1 / \rho} (x^1) + \Lambda + \frac{\rho(L_f +
  L_{\omega})^2}{\alpha^2} > 0$ and recall that $\diam
  (\mathcal{L}_{a\Delta}) \leq \mathsf{B}_{a\Delta}$.
\end{lem}

\textbf{Lemma \ref{lem:2.2}} provides a useful characterization of the tail probability on the norm of the iterations. Now that the bound holds for all $x^k$, we can immediately condition on the event that $\Theta(K)$ iterations lie in the bounded set, to retrieve an $\mathcal{O}(1/\sqrt{K})$ convergence rate.

\begin{thm}
  \label{thm:2.2}Under the same conditions as \textbf{Lemma
  \ref{lem:2.2}}, given $\delta \in (0, 1/4)$, at least with probability $1-p, p \in (2\delta, 1)$,
  $(1 - 2p^{-1} \delta) K$ iterations will lie in the ball with radius $\mathsf{R} (\delta) =
  \mathsf{B}_{\delta^{- 1} \Delta} + \tfrac{4(L_f + L_{\omega})}{\alpha \sqrt{K}}$ and
\[ \min_{1 \leq k \leq K} \mathbb{E} [\| \nabla \psi_{1/\rho}(x^k) \|^2] \leq
   \frac{p}{p - 2 \delta} \cdot \frac{2 \rho}{\rho - \tau - \kappa}
  \left[ D + \frac{\rho }{2 \alpha^2} (L_f + L_{\omega})^2 \right] \left(
  \frac{\rho + \tau + \kappa}{K} + \frac{\alpha (\mathsf{G}_{\delta} +
  1)}{\sqrt{K}} \right), \]
where $\mathsf{G}_\delta \assign \max_{z \leq \mathsf{R}(\delta)} \mathcal{G} (z).$
\end{thm}

\textbf{Theorem \ref{thm:2.2}} shows that, with constant probability, we can retrieve $\Ocal(1/\sqrt{K})$ convergence rate after $K$ iterations. This probability argument can be further improved, for example, by running the algorithm independently multiple times \cite{davis2019proximally}. The analysis in this section serves as a step-stone for the next section, where we deal with non-Lipschitz optimization even without knowing $\mathcal{G}(\cdot)$.

\section{Weakly convex optimization under unknown Lipschitzness}\label{sec:unknown-nolip}

While the analysis in \textbf{Section~\ref{sec:known-nonlip}} already
extends the solvability of weakly convex optimization to non-Lipschitz
functions, it  relies on the knowledge of an explicit growth function $\mathcal{G}
(\| x \|)$, which bounds the local Lipschitzness. However, in many real scenarios, either access to $\mathcal{G}
(\| x \|)$ is not viable, or the bound lacks a predefined functional form. In these scenarios, we assume that the growth function is unknown a priori.
\begin{enumerate}[leftmargin=35pt,itemsep=2pt,label=\textbf{D\arabic*:},ref=\rm{\textbf{D\arabic*}}]
  \item For all fixed $x \in \dom \omega$, $f_x (z, \xi) - f_x (y, \xi) \leq \textup{\textsf{Lip}} (x, \xi) \| z - y \|$ for all $y, z$; $\xi \sim \Xi$. \label{D1}
\end{enumerate}

Although \ref{D1} and \ref{C1} look similar, they are very different in nature. The most direct consequence is that $\textup{\textsf{Lip}} (x, \xi)$ is sample-dependent, which means any stepsize strategy based on it will introduce bias in the stochastic algorithm. To resolve this issue, our new stepsize policy relies on constructing an
estimator of $\textup{\textsf{Lip}} (x, \xi)$. We introduce the
property of ``allowing a good estimation of Lipschitz constant'', which we call the
{\textit{reference Lipschitz continuity}}.

\subsection{Reference Lipschitz continuity}
\begin{definition}[Reference Lipschitz continuity]
A stochastic model function $f_x (y, \xi)$ satisfies reference Lipschitz continuity if
  \begin{itemize}
    \item Given an anchor $x$, $f_x(\cdot,\xi)$ is globally Lipschitz with a Lipschitz constant $\textup{\textsf{Lip}} (x, \xi)$
    
    \item Given $\xi,\xi' \sim \Xi$,
    $\mathbb{E}_{\xi, \xi'} [| \textup{\textsf{Lip}} (x, \xi) - \textup{\textsf{Lip}} (x, \xi') |^2] \leq \sigma^2 < \infty$..
  \end{itemize}
\end{definition}

The first property is determined by the stochastic model function, and is common since most model functions are compositions of Lipschitz functions and linear expansions. The second property states that the expected difference between models' Lipschitzness is bounded by noise parameter $\sigma$.
One direct outcome of reference Lipschitz continuity is that we can cheaply estimate
$\textup{\textsf{Lip}} (x, \xi)$ based on $\textup{\textsf{Lip}} (x, \xi')$, where $\xi$
and $\xi'$ are two independent samples drawn from $\Xi$.

\begin{exple}[(Sub)gradient] \label{example-2}
\rm{ $f_x (y, \xi) = f (x, \xi) + \langle \nabla f (x,
  \xi), y - x \rangle$\\
  The model is globally Lipschitz with $\textup{\textsf{Lip}} (x, \xi) = \| \nabla
  f (x, \xi) \|$. If $\mathbb{E} [\| \nabla f (x) - \nabla f
  (x, \xi) \|^2] \leq \sigma^2$, then \[\mathbb{E} [| \textup{\textsf{Lip}} (x, \xi)
  - \textup{\textsf{Lip}} (x, \xi') |] \leq \mathbb{E} [\| \nabla f (x, \xi) -
  \nabla f (x, \xi') \|] \leq 2 \sigma.\]
  Even if $f$ is nonsmooth, the property may still hold. One
  example is the composition problem $f (x, \xi) = h (c (x, \xi))$, where $h$
  is $L_h$-Lipschitz continuous and $c$ is smooth. Then $\partial f (x, \xi) =
   \nabla c (x, \xi) \partial h (c (x, \xi))$. If we estimate the Lipschitz
  constant with $L_h \| \nabla c (x, \xi) \|$, then $\mathbb{E} [| \textup{\textsf{Lip}}
  (f_x (\cdot, \xi)) - \textup{\textsf{Lip}} (x, \xi') |] \leq 2 L_h \sigma$.}
\end{exple}

\begin{exple}[Proximal linear] \label{example-3}$f_x (y, \xi) = h (c (x, \xi) + \langle \nabla c (x, \xi),
  y - x \rangle)$\\
\rm{When $h$ is $L_h$ Lipschitz continuous, the model is globally Lipschitz with
  $\textup{\textsf{Lip}} (x, \xi) = L_h \| \nabla c (x, \xi) \|$. If $\mathbb{E} [\| \nabla c (x) - \nabla c
  (x, \xi) \|^2] \leq \sigma^2$, then $\mathbb{E} [|
  \textup{\textsf{Lip}} (x, \xi) - \textup{\textsf{Lip}} (x, \xi') |] \leq L_h
  \mathbb{E} [\| \nabla c (x, \xi) - \nabla c (x, \xi') \|] \leq 2 L_h \sigma
  .$}
\end{exple}

\begin{exple}[Truncated model] \label{example-4}  \rm{  $f_x (y, \xi) = \max \{ f (x, \xi) + \langle \nabla f (x,
  \xi), y - x \rangle, \ell \}.$ \\The model is $\| \nabla f (x, \xi) \|$-Lipschitz and the reasoning of reference Lipschitz continuity is the same as in \textbf{Example \ref{example-2}}.} Note that the truncated model encompasses stochastic Polyak stepsize as a special case \cite{schaipp2023a}.
\end{exple}

In this section, we would assume that our model satisfies the reference Lipschitz continuity.

\begin{enumerate}[leftmargin=35pt,itemsep=2pt,start=2,label=\textbf{D\arabic*:},ref=\rm{\textbf{D\arabic*}}]
  \item The stochastic model $f_x (\cdot, \xi)$ satisfies the reference Lipschitz continuity with noise parameter $\sigma$. \label{D2}
\end{enumerate}

\subsection{Algorithm design and analysis}

As we did in \textbf{Section \ref{sec:known-nonlip}}, before getting down to the algorithm design, we first 
need to see what happens to our potential reduction in this new setting. Firstly, \textbf{Lemma \ref{lem:3.4}} characterizes the descent property of our potential function under the assumption that $\gamma_k$ is independent of $\xi^k$.
\begin{lem}
  \label{lem:3.4}Suppose \ref{A1} to \ref{A5} as well as \ref{D1}, \ref{D2} hold, then given $\rho > \kappa + \tau$, 
\[ \mathbb{E}_k [\psi_{1 / \rho} (x^{k + 1}) ] \leq \psi_{1 / \rho} (x^k) -
     \myfrac{\rho (\rho - \tau - \kappa)}{2 (\gamma_k - \kappa)} \| \hat{x}^k - x^k \|^2
     + \mathbb{E}_k \Big[ \myfrac{\rho }{2 \gamma_k(\gamma_k - \kappa)} (\textup{\textsf{Lip}} (x^k, \xi^k) +
     L_{\omega})^2 \Big], \]
     where $\gamma_k$ is chosen to be independent of $\xi^k$ and is considered deterministic here.
\end{lem}

Compared to \textbf{Lemma \ref{lem:2}}, bounding the error term $\myfrac{\rho (\textup{\textsf{Lip}} (f_{x^k} (\cdot, \xi)) +
     L_{\omega})^2}{2 \gamma_k (\gamma_k - \kappa)}$ becomes more challenging due to the lack of information about its growth. However, thanks to
      the reference Lipschitz continuity,  by sampling  $\xi'$, an independent copy of $\xi^k$, we can utilize $\textup{\textsf{Lip}}(x^k, \xi')$ as a surrogate for $\textup{\textsf{Lip}}(x^k, \xi^k)$. 
      The preference of $\textup{\textsf{Lip}}(x^k, \xi')$ over $\textup{\textsf{Lip}}(x^k, \xi^k)$  is driven by the fact that $\textup{\textsf{Lip}}(x^k, \xi^k)$ is correlated with $x^{k +
1}$, which significantly complicates the analysis~\citep{andradottir1996scaled}. 
To mitigate the
impact of large noise $\sigma$ on the accuracy of our estimation, we  clip the estimator by a threshold $\alpha>0$: $\max \{ \textup{\textsf{Lip}} (x^k, \xi'), \alpha \}$. 
Consequently, we set \[ \gamma_k = \mathcal{O}( \max \{ \textup{\textsf{Lip}} (x^k, \xi'), \alpha \} \cdot \sqrt{K}).\]

\begin{rem}
When $f_x (y, \xi)
= f (x, \xi) + \langle \nabla f (x, \xi), y - x \rangle$ and $\xi = \xi'$, we
retrieve the gradient clipping technique. Our approach can be seen
as a generalization of gradient clipping to the proximal setting.
\end{rem}

The next theorem establishes the asymptotic result and confirms our intuition.

\begin{thm}
  \label{thm:3}Under the same conditions as \textbf{Lemma \ref{lem:3.4}}, \ref{A6}, \ref{A7},
if $\gamma_k = \rho + \kappa + \tau + \max \{ \textup{\textsf{Lip}} (x^k, \xi'), \alpha \} k^{\zeta}$, $\zeta \in ( \myfrac{1}{2}, 1 )$, as $k \rightarrow \infty$,
  $\{ \| x^k \| \}$ is bounded with probability $1$; $ \{\inf_{j \leq k} \| \nabla \psi_{1
  / \rho} (x^j) \| \}$ converges to $0$ almost surely.
\end{thm}

To obtain a non-asymptotic result, we again apply the probabilistic analysis to obtain the tail bound.

\begin{lem}
  \label{lem:3.5}Under the same conditions of~\textbf{Lemma \ref{lem:3.4}} as well as \ref{A6}, \ref{A7},
  if we take $\gamma_k = \rho + \tau + \kappa + \max \{ \textup{\textsf{Lip}} (f_{x^k} (\cdot,
  \xi')), \alpha \}  \sqrt{K}$, then the tail bound
  \[ \mathbb{P} \Big\{ \| x^k \| \geq \mathsf{B}_{a \Delta} + \myfrac{4(\alpha + \sigma +
     L_{\omega})}{\alpha \sqrt{K}} \Big\} \leq \myfrac{2\Delta}{a \Delta +
     \Lambda}, \]
  holds for all $2 \leq k \leq K$, where $\Delta = \psi_{1 / \rho} (x^1) + \Lambda + \myfrac{\rho}{\alpha^2} (\alpha + \sigma+
  L_{\omega})^2 > 0$.
\end{lem}

\begin{thm}
  \label{thm:3.2}Assuming the  conditions of \tmtextbf{Lemma
  \ref{lem:3.5}} hold, then given $\delta \in (0, 1/4)$,  with probability at least $1-p, p \in (2\delta, 1)$, $(1 - 2p^{-1} \delta) K$ iterations will lie in the ball with radius $\mathsf{R} (\delta) =
  \mathsf{B}_{\delta^{- 1} \Delta} + \myfrac{4 (\alpha + \sigma +
  L_{\omega})}{\alpha \sqrt{K}}$ and
  \[ \min_{1 \leq k \leq K} \mathbb{E} [\| \nabla
     \psi_{1 / \rho} (x^k) \|^2] \leq \myfrac{p}{p - 2 \delta} \cdot \myfrac{2 \rho}{\rho - \tau - \kappa}
  \left[ D + \myfrac{\rho}{\alpha^2} (\alpha + \sigma + L_{\omega})^2 \right]
  \left( \myfrac{\rho + \tau + \kappa}{K} + \myfrac{\alpha + \mathsf{G}_{\delta}}{\sqrt{K}} \right), \]
  where $\mathsf{G}_\delta \assign \max_x \sup_{\xi\sim\Xi} \textup{\textsf{Lip}} (x, \xi), \| x
  \| \leq \mathsf{R}(\delta)$.
\end{thm}

\section{Experiments} \label{sec:exp}

In this section, we perform numerical experiments to demonstrate the effectiveness of our proposed methods.
We consider the following robust non-linear regression problem:
\begin{equation}
  \min_x \quad \frac{1}{m} \sum_{i = 1}^m | r (x, a_i) - b_i | =:
  \frac{1}{m} \sum_{i = 1}^m f (x, \xi_i),
\end{equation}
where, given observations $\{a_i\}$ from $A\in\mathbb{R}^{m\times n}$, regression model $r(x, a)$ and  target label $b_i$, we aim to fit the model coefficient $x$ given problem data.  The following table summarizes our tested regression models and their Lipschitz properties. These three functions exhibit increasing non-Lipschitzness. 

\ifthenelse{\boolean{doublecolumn}}{
\begin{table*}[h]
}{
\begin{table}[h]
}
\caption{Nonlinear regression models. $r(a,x)=\langle a, x \rangle^2$ represents the standard robust phase retrieval problem; $r(a,x)=\langle a, x \rangle^5 + \langle a, x \rangle^3 + 1$ is a high-order polynomial of $\langle a, x \rangle$; while $e^{\langle a, x \rangle} + 10$ exhibits exponential growth.}
  \begin{tabular}{ccccc}
    \toprule
    Loss& $r (x, a)$ & $\partial f (x, \xi)$ & $\mathcal{G} (\| x \|)$ &
    $\textsf{Lip} (f_x (y, a))$\\
    \midrule
    $r_1$ & $\langle a, x \rangle^2$ & $\sign (\langle a, x \rangle^2 - b)
    \cdot 2 \langle a, x \rangle a$ & $\| x \|$ & $2 | \langle a, x \rangle
    | \cdot \| a \|$\\
    $r_2$ &$\langle a, x \rangle^5 + \langle a, x \rangle^3 + 1$ & $\sign
    (r_2 - b) \cdot (5 \langle a, x \rangle^4 + 3 \langle a, x \rangle^2) a$ &
    $5 (\| x \|^4 + \| x \|^2)$ & $5 (\| a \|^4 \| x \|^4 + \| a \|^2 \| x
    \|^2)$\\
    $r_3$ &$e^{\langle a, x \rangle} + 10$ & $\sign (e^{\langle a, x \rangle} -
    b) \cdot e^{\langle a, x \rangle} \cdot a$ & $e^{A \| x \|}$ ($A = 3.0)$ &
    $e^{\langle a, x \rangle} \cdot \| a \|$\\
    \bottomrule
  \end{tabular}

\ifthenelse{\boolean{doublecolumn}}{
\end{table*}
}{
\end{table}
}

\subsection{Experiment setup}
\textbf{Dataset}. We let $m=300,n=100$. Data generation is consistent with {\cite{deng2021minibatch}}, where, given condition number parameter $\kappa\geq 1$, we compute
$A=QD,Q\in\mathbb{R}^{m\times n}$. Here each element of $Q$ is drawn from $\mathcal{N}(0,1)$;
$D=\tmop{\text{{diag}}}(d),d\in\mathbb{R}^{n},d_{i}\in[1/\kappa,1]$ for all $i$.
Then a true signal $\hat{x}\sim\mathcal{N}(0,I)$ is generated, giving the measurements $b$ by formula $b_{i}=r(x, a_i)$.
We randomly perturb $p_{\text{fail}}$-fraction of the measurements
with $\mathcal{N}(0,25)$ noise added to them to simulate data corruption.

\begin{enumerate}[label=\textbf{\arabic*)}, ref=Ex\arabic*, leftmargin=*]
\item \textbf{Dataset}. We follow {\cite{deng2021minibatch}} and set
set $\kappa\in\{1,10\}$ and $p_{\text{fail}}\in\{0.2,0.3\}$.
\item \textbf{Initial point}. We generate $x'\sim\mathcal{N}(0,I_{ n})$
and start from $x^{1}=\frac{10x'}{\|x'\|}$ for $r_1$, $x^{1}=\frac{x'}{\|x'\|}$ for $r_2, r_3$. 
\item \textbf{Stopping criterion}. We run algorithms for 400 epochs ($K=400m$).
Algorithms stop if $f \leq 1.2f(\hat{x})$ . 
\item \textbf{Stepsize}. We let $\gamma_k = \theta \cdot \sqrt{K}$ for vanilla algorithms; $\gamma_k = \theta \cdot \mathcal{G}(\|x^k\|) \sqrt{K}$ for adaptive stepsize with known growth condition; $\gamma_k = \theta \cdot \max \{ \textsf{Lip}(x^k, \xi'), \alpha  \}\sqrt{K}$ for adaptive stepsize with unknown growth condition. $\theta\in[10^{-2},10^{1}]$ serves as a hyper-parameter.
\item \textbf{Clipping}. Clipping parameter $\alpha$ is set to $1.0$.
\item \textbf{Mirror descent}. For experiments on mirror descent, we construct kernels for $r_1, r_2$ according to \cite{davis2018stochastic}.
\end{enumerate}

\subsection{Comparing different stepsizes}

We conduct experiments to investigate the number of iterations for each stepsize to converge, under different choices of $\theta$. As the experiments suggest, when the function exhibits low-order growth, our adaptive choices tend to be conservative. However, when the function exhibits high-order growth, our adaptive stepsize tends to converge within a reasonable range of stepsizes. It is worth noticing that for problem $r_2$, {\sgd} diverges for $\theta \sim 10^8$, while our proposed approaches work robustly. Moreover, we notice that our adaptive stepsize based on reference Lipschitz property never diverges in practice, although it is sometimes conservative on problems where function growth is mild (such as $r_1$).

\begin{figure*}[h] \label{fig-1}
\centering
	\includegraphics[scale=0.20]{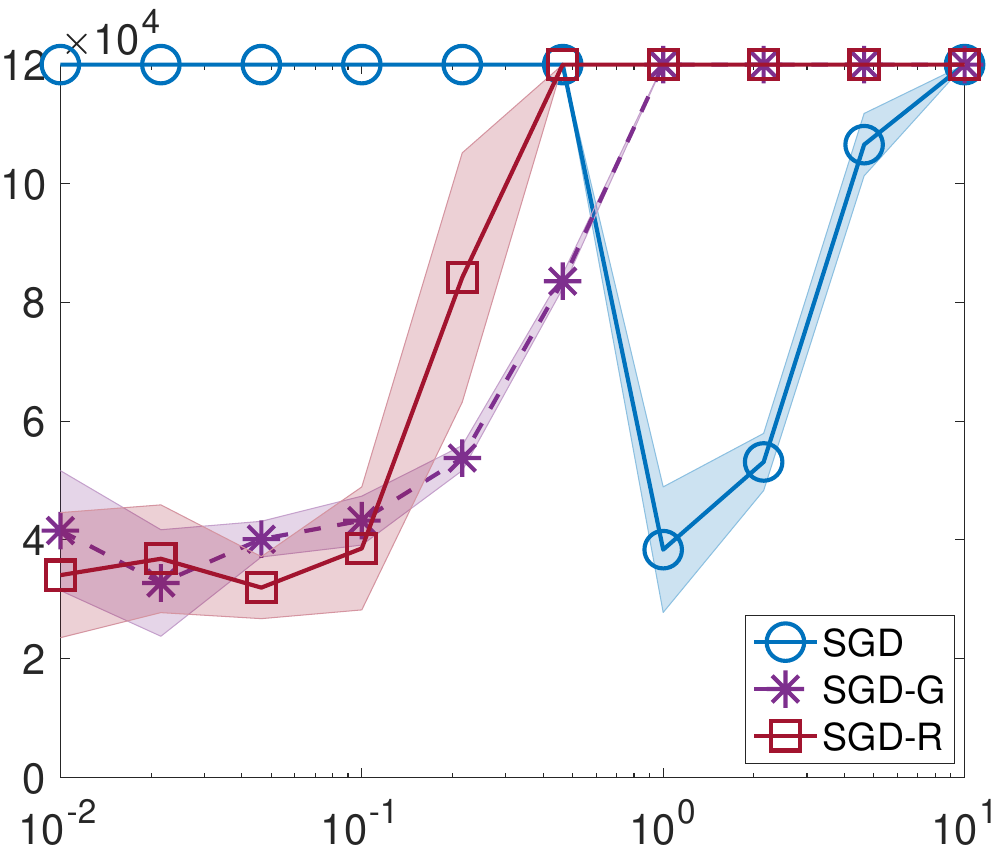}
	\includegraphics[scale=0.20]{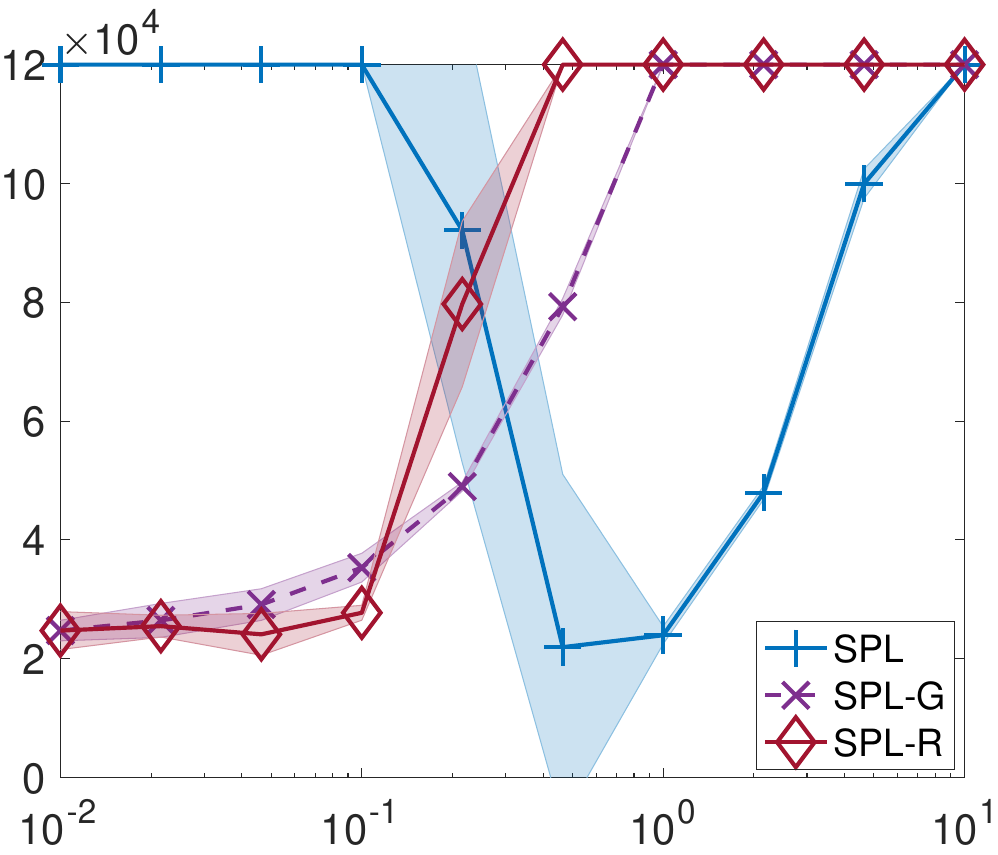}
	\includegraphics[scale=0.20]{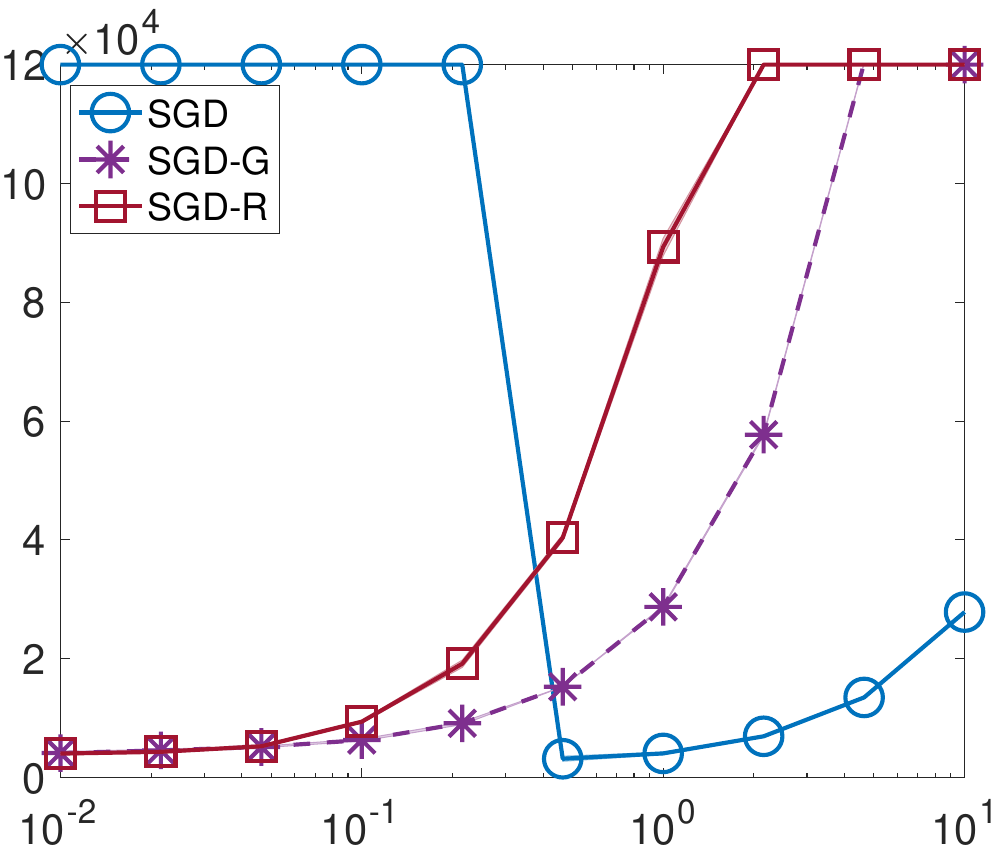}
	\includegraphics[scale=0.20]{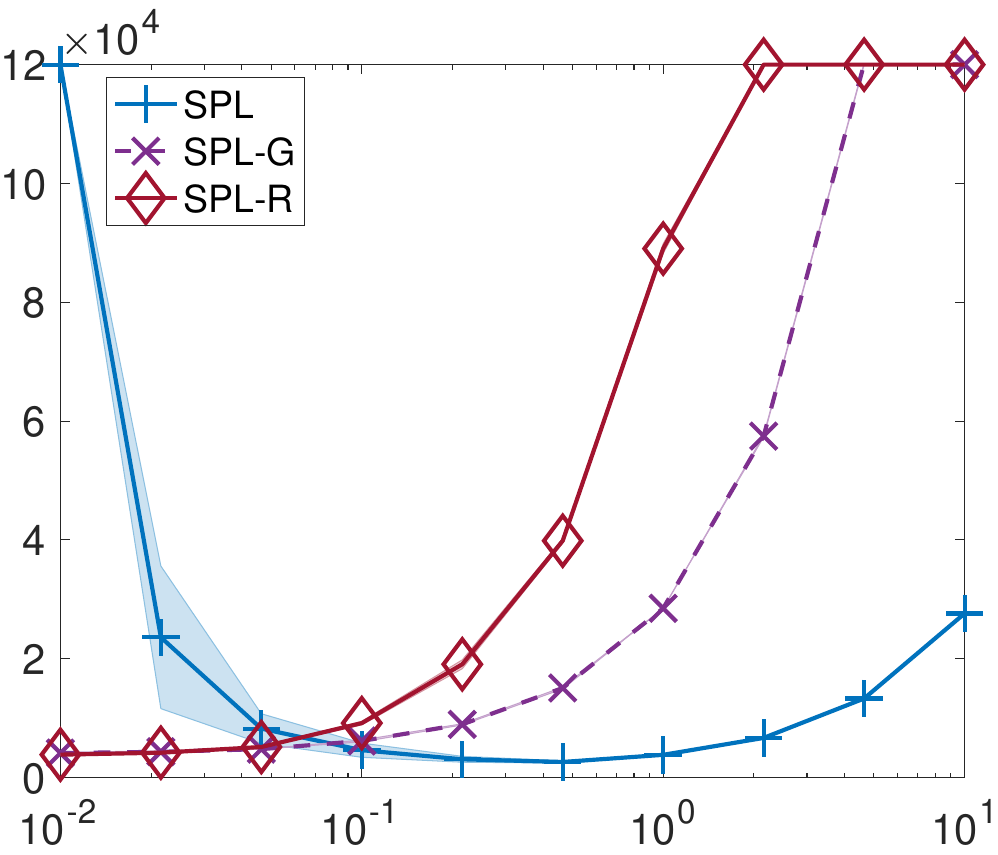}
	\caption{Problem $r_1$. Left two: $(\kappa,p_{\text{{fail}}})=(10,0.2)$; Right two: $(\kappa,p_{\text{{fail}}})=(10,0.3)$. x-axis: parameter $\theta$; y-axis: number of iterations. {\sgd} denotes vanilla SGD; {\sgd\texttt{-G}} denotes SGD adaptive to known Lipschitzness; {\sgd\texttt{-R}} denotes SGD adaptive to unknown Lipschitzness. The same applies to {\spl}.}
\end{figure*}

\begin{figure*}[h] \label{fig-2}
\centering
	\includegraphics[scale=0.20]{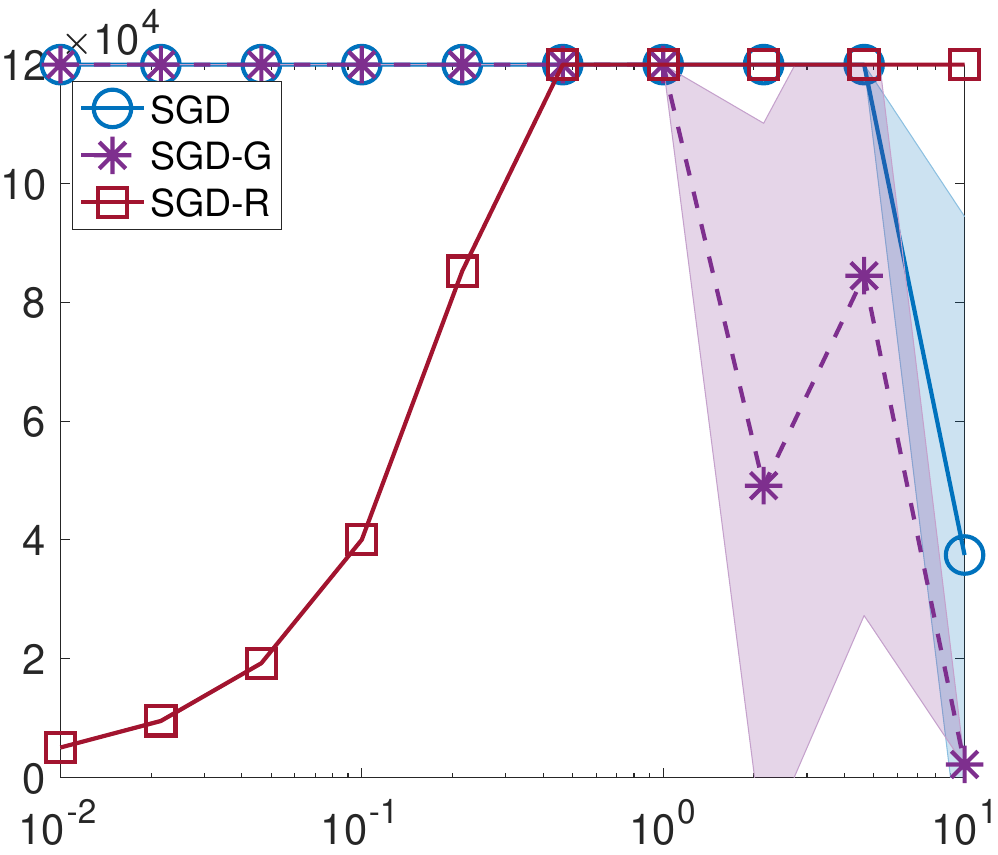}
	\includegraphics[scale=0.20]{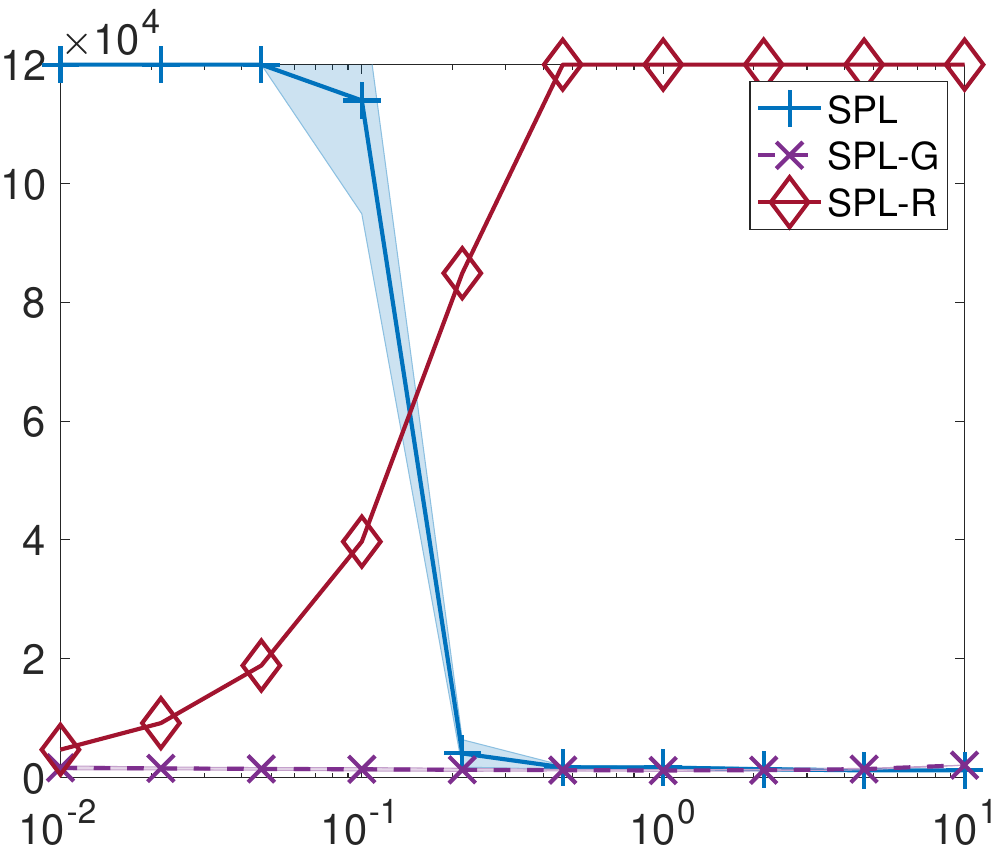}
	\includegraphics[scale=0.20]{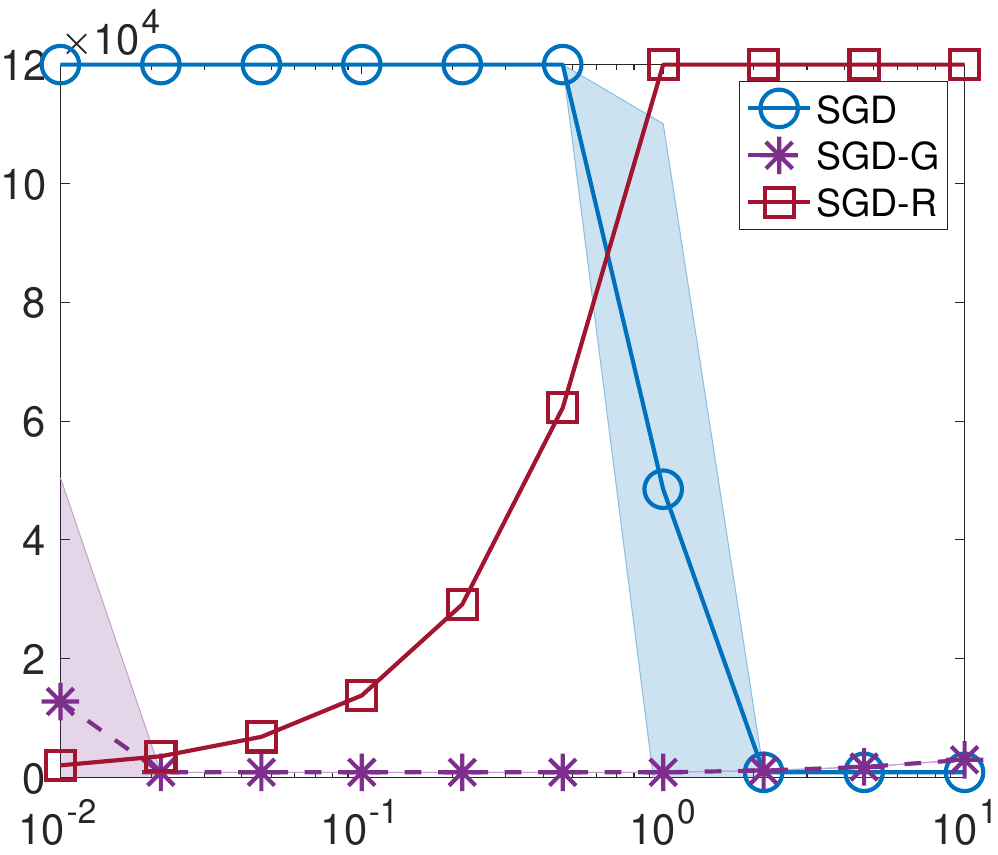}
	\includegraphics[scale=0.20]{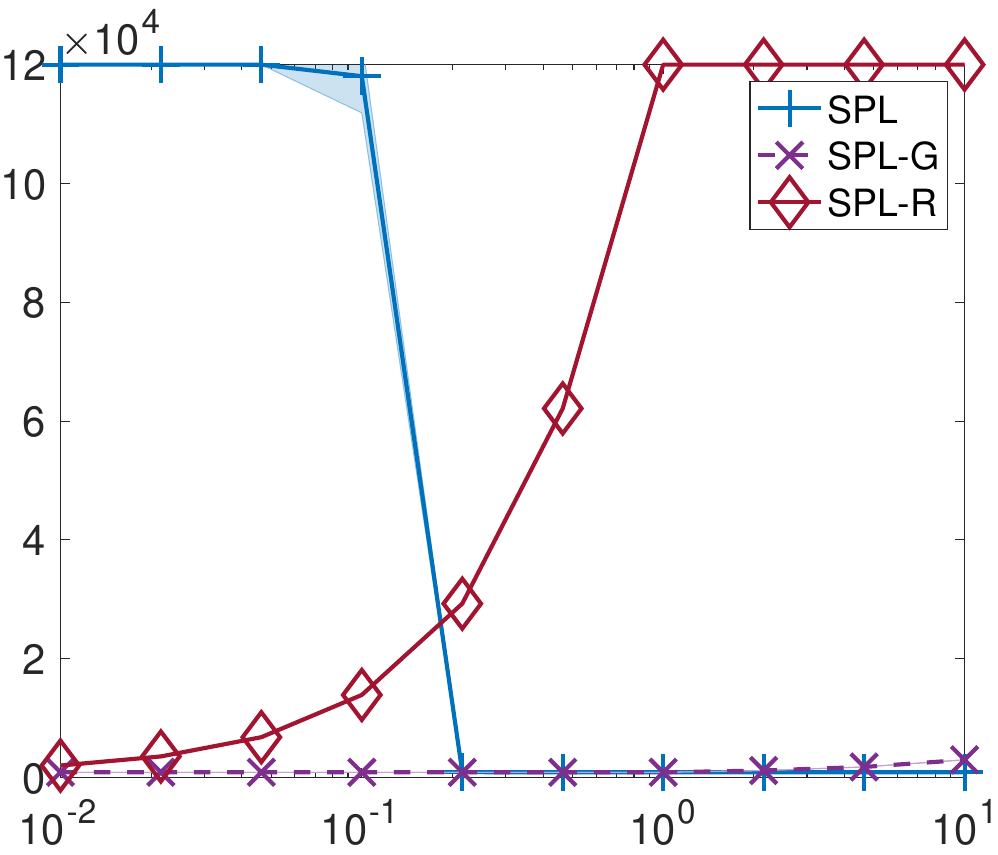}
	\caption{Problem $r_2$. Left two: $(\kappa,p_{\text{{fail}}})=(1,0.2)$; Right two: $(\kappa,p_{\text{{fail}}})=(10,0.3)$. x-axis: parameter $\theta$; y-axis: number of iterations.}
\end{figure*}

\begin{figure*}[h] \label{fig-3}
\centering
	\includegraphics[scale=0.20]{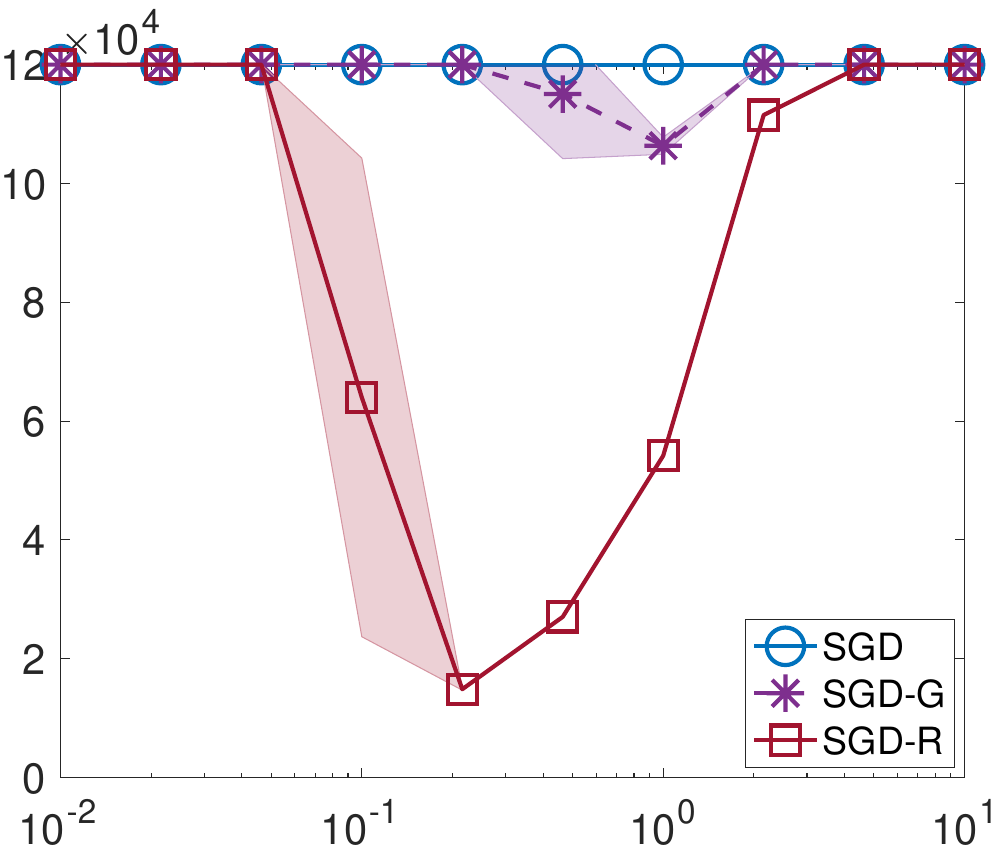}
	\includegraphics[scale=0.20]{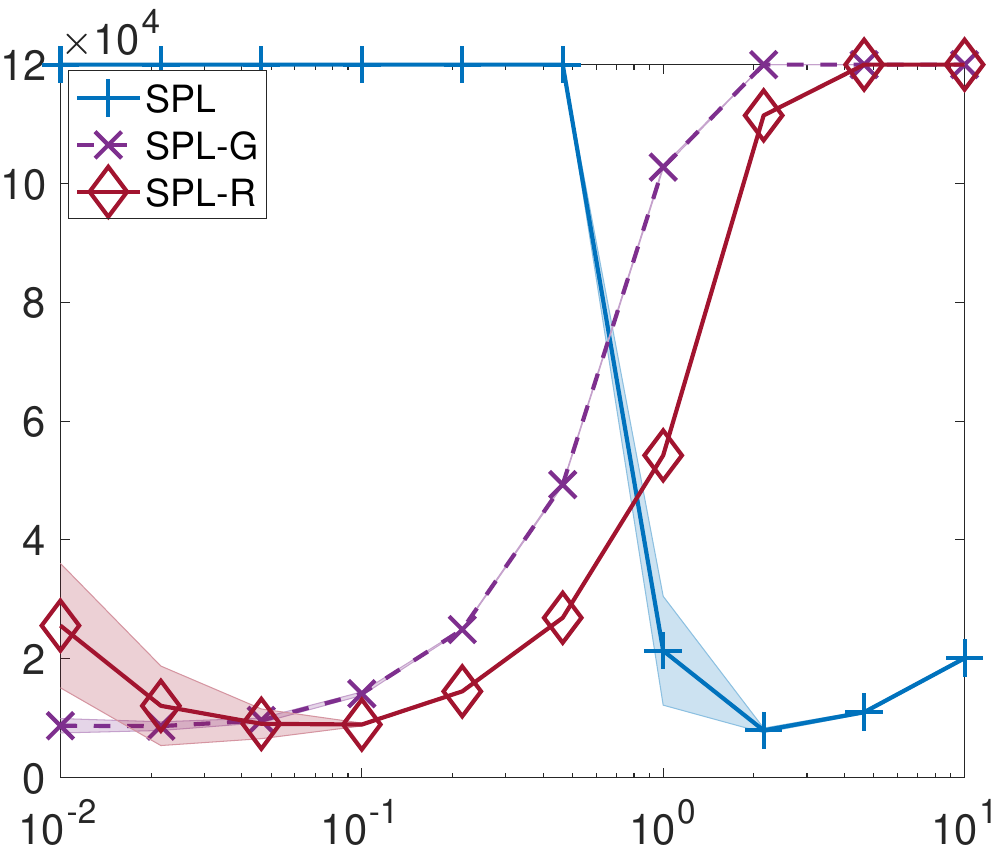}
	\includegraphics[scale=0.20]{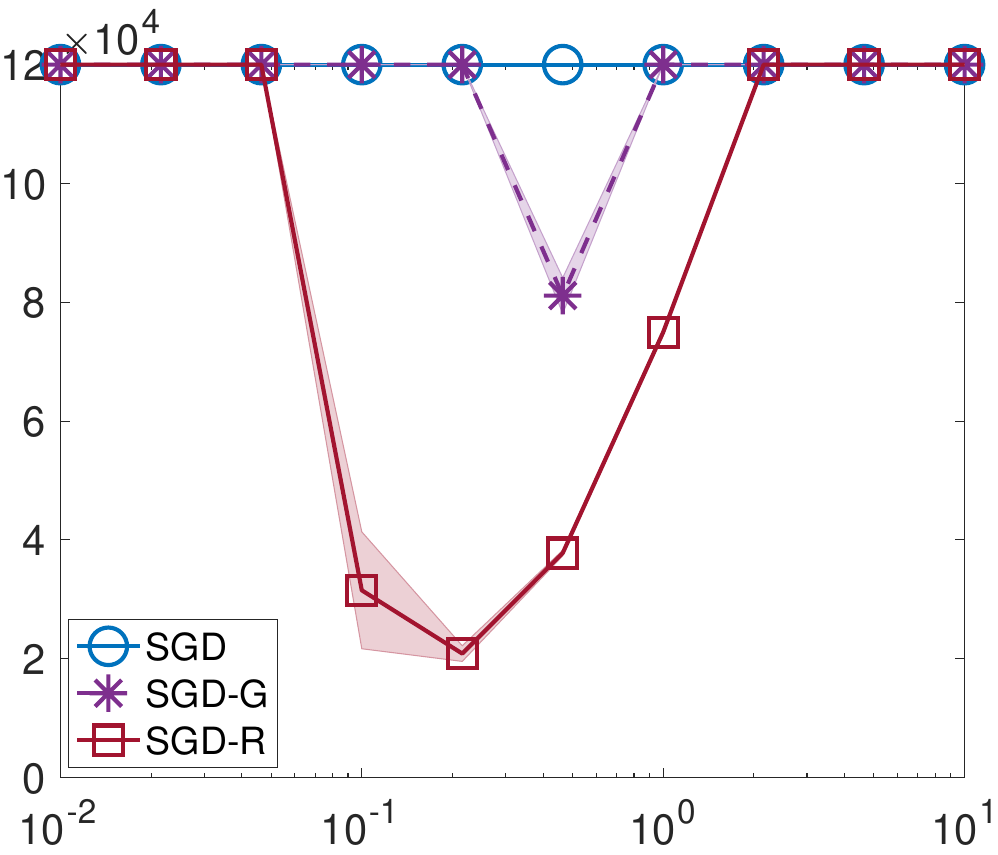}
	\includegraphics[scale=0.20]{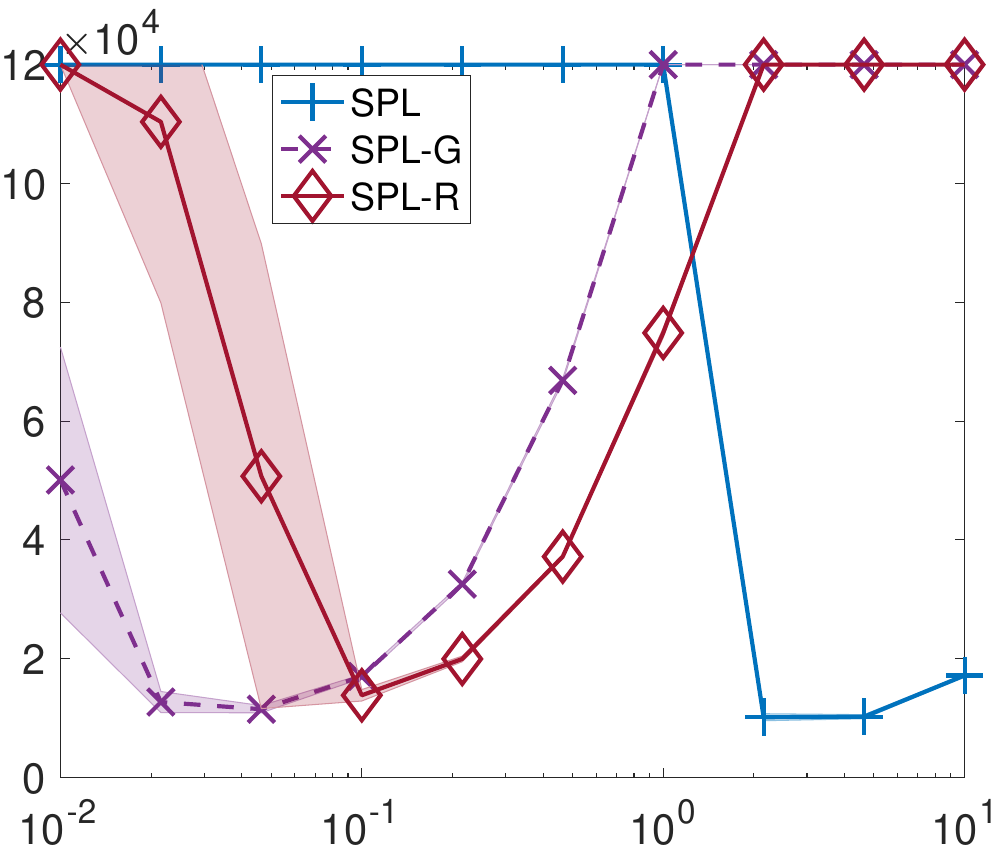}
	\caption{Problem $r_3$. Left two: $(\kappa,p_{\text{{fail}}})=(1,0.2)$; Right two: $(\kappa,p_{\text{{fail}}})=(1,0.3)$. x-axis: parameter $\theta$; y-axis: number of iterations.}
\end{figure*}

\subsection{Comparison with mirror descent}
Last, we compare our proposed method with the commonly adopted mirror descent approach  for non-Lipschitz problems. We test both mirror descent and our proposed {\sgd}-based approaches. \\

\begin{figure*}[h]
\centering
	\includegraphics[scale=0.20]{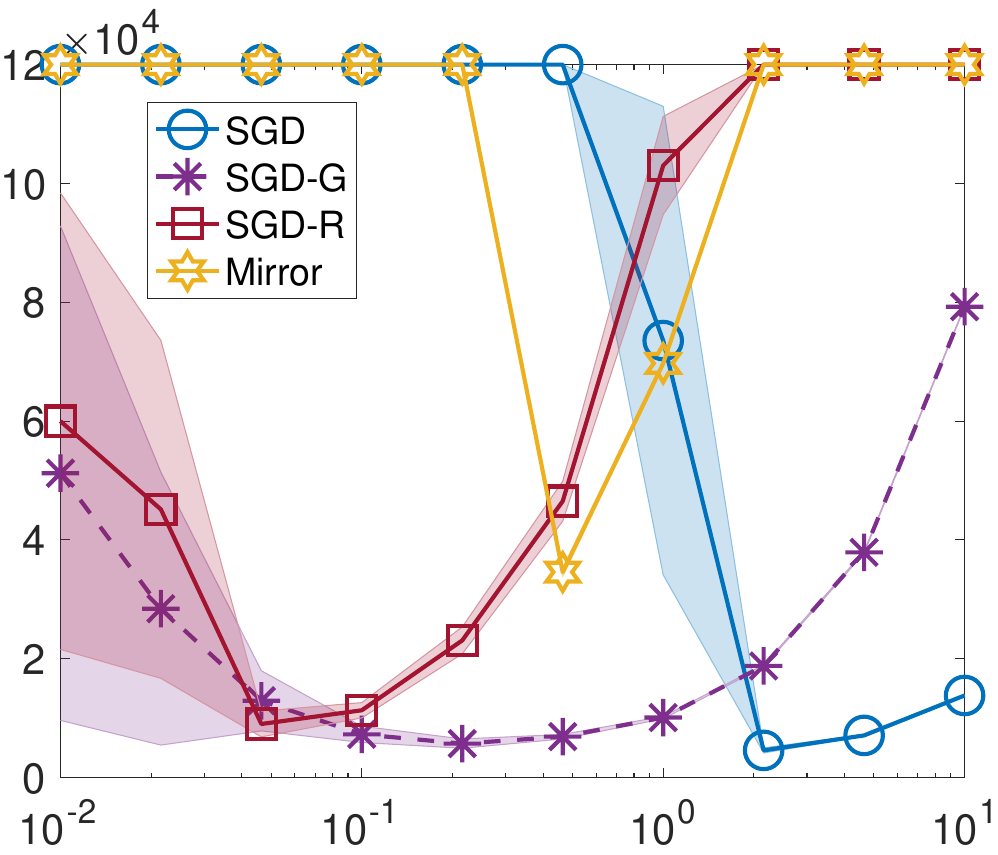}
	\includegraphics[scale=0.20]{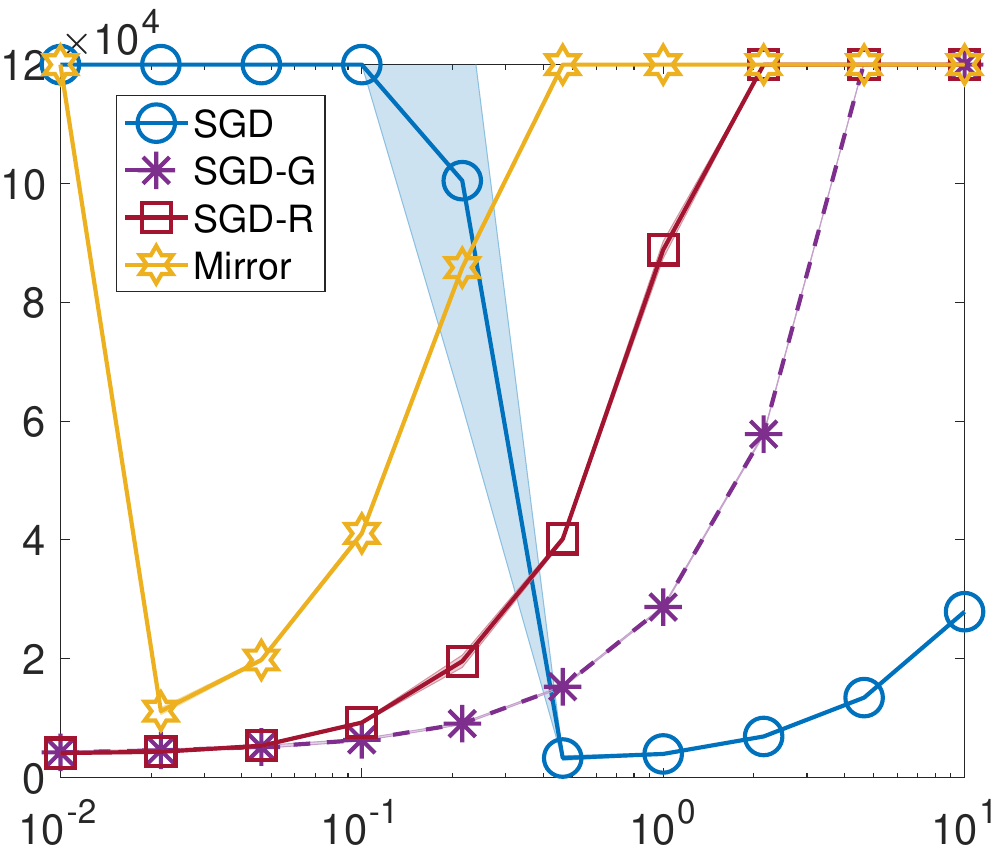}
	\includegraphics[scale=0.20]{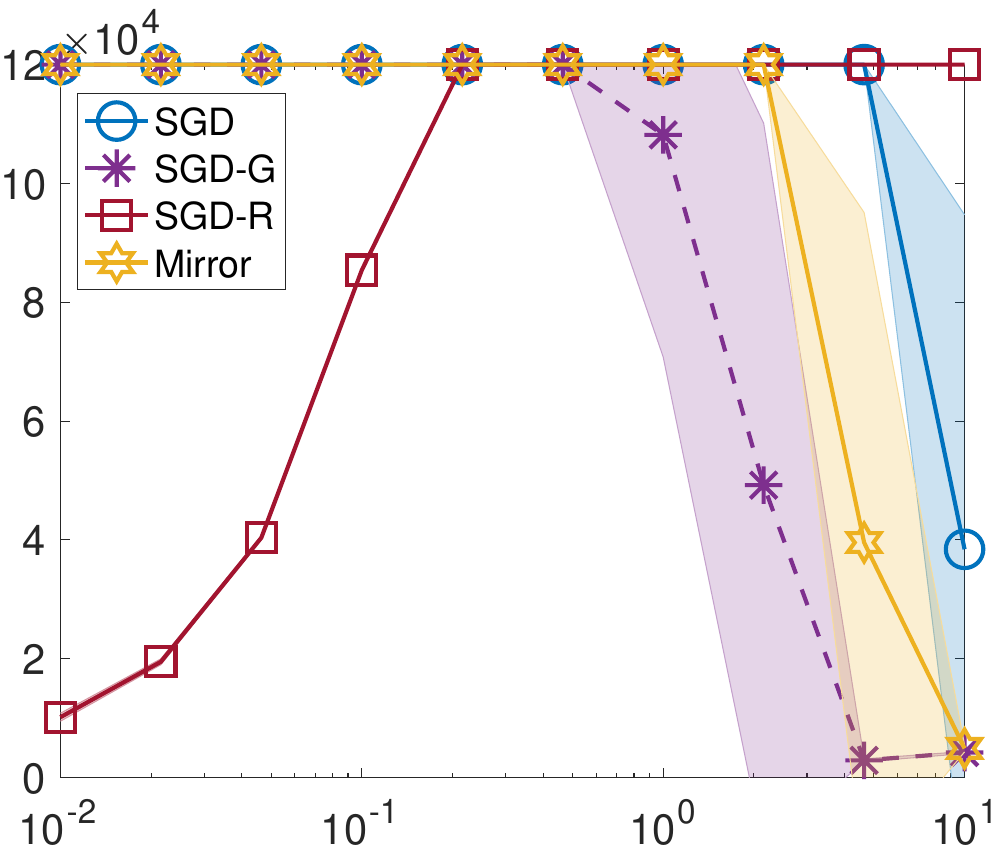}
	\includegraphics[scale=0.20]{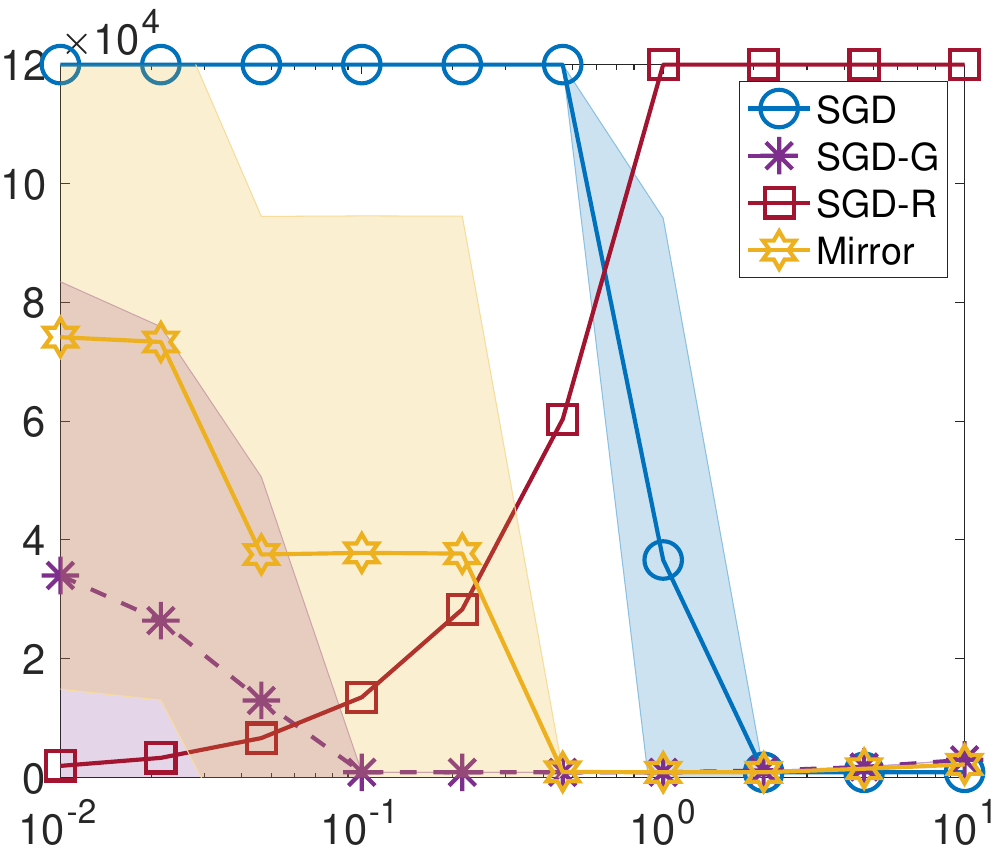}
	\caption{Left two: Problem $r_1$, $(\kappa,p_{\text{{fail}}})=(1,0.3)$; Right two: Problem $r_2$, $(\kappa,p_{\text{{fail}}})=(1,0.3)$. x-axis: parameter $\theta$; y-axis: number of iterations.\label{fig-4}}
\end{figure*}

As \textbf{Figure \ref{fig-4}} suggests, we see that mirror descent indeed often exhibits more stable performance compared to vanilla {\sgd}. However, we see that our approaches still exhibit superior convergence performance. Moreover, each mirror descent iteration involves a root-finding subproblem, which is far more complicated than our approaches.

\section{Conclusions}
We develop novel adaptive stepsize (regularization) strategies and show that  for weakly convex objectives without Lipschitz continuity, stochastic model-based methods can still converge at the desirable $\mathcal{O}(1/\sqrt{K})$ rate  with constant failure probability. To our knowledge, this is achieved under the least restrictive assumptions known to date. A promising direction for future research is the adaptation of our analyses to more sophisticated methods, such as momentum-based or adaptive gradient methods.
\renewcommand \thepart{}
\renewcommand \partname{}

\newpage 
\bibliographystyle{plainnat}
\bibliography{ref.bib}


\doparttoc
\faketableofcontents
\part{}

\newpage
\appendix
\addcontentsline{toc}{section}{Appendix}
\part{Appendix} 
\parttoc

\paragraph{Structure of the appendix}

The appendix is organized as follows. In \textbf{Section \ref{app:known-lip}}, \textbf{Section \ref{app:known-nolip}} and \textbf{Section \ref{app:unknown-nolip}}, we respectively prove the results for weakly convex optimization under different Lipschitz continuity assumptions. In \textbf{Section \ref{app:cvx}} and its subsections, we extend our results to convex optimization.
 
\newpage
\section{Proof of results in Section \ref{sec:known-lip}}
\label{app:known-lip}
\subsection{Proof of Lemma \ref{lem:1}}

By the optimality conditions of the proximal 
subproblems \eqref{eq:update-model-based}, we have, for any $\xi^k\sim\Xi$, that
\[ f_{x^k} (x^{k + 1}, \xi^k) + \omega (x^{k + 1}) + \frac{\gamma_k}{2} \|
   x^{k + 1} - x^k \|^2 \leq f_{x^k} (\hat{x}^k, \xi^k) + \omega (\hat{x}^k) +
   \frac{\gamma_k}{2} \| \hat{x}^k - x^k \|^2 - \frac{\gamma_k - \kappa}{2} \| x^{k +
   1} - \hat{x}^k \|^2 \]
\[ f (\hat{x}^k) + \omega (\hat{x}^k) + \frac{\rho}{2} \| \hat{x}^k - x^k \|^2
   \leq f (x^{k + 1}) + \omega (x^{k + 1}) + \frac{\rho}{2} \| x^{k + 1} - x^k
   \|^2  \]
Summing over the above two relations, we deduce that
\begin{align}
  & \frac{\gamma_k - \rho}{2} \| x^{k + 1} - x^k \|^2 - \frac{\gamma_k -
  \rho}{2} \| \hat{x}^k - x^k \|^2 + \frac{\gamma_k  - \kappa}{2} \|
  x^{k + 1} - \hat{x}^k \|^2 \nonumber\\
  \leq{} & f (x^{k + 1}) - f_{x^k} (x^{k + 1}, \xi^k) + f_{x^k} (\hat{x}^k,
  \xi^k) - f (\hat{x}^k) \nonumber
\end{align}
Conditioned on $\xi^1, \ldots, \xi^{k-1}$ and taking expectation with respect to $\xi^k$, we have
\begin{align}
  & \frac{\gamma_k - \rho}{2} \mathbb{E}_k [\| x^{k + 1} - x^k \|^2] -
  \frac{\gamma_k - \rho}{2} \| \hat{x}^k - x^k \|^2 + \frac{\gamma_k -
  \kappa}{2} \mathbb{E}_k [\| x^{k + 1} - \hat{x}^k \|^2] \nonumber\\
  \leq{} & \mathbb{E}_k [f (x^{k + 1}) - f_{x^k} (x^k, \xi^k)] +\mathbb{E}_k [L
  (\xi^k) \| x^{k + 1} - x^k \|] + \frac{\tau}{2} \| \hat{x}^k - x^k \|^2
  \label{eqn:a-1} \\
  ={} & \mathbb{E}_k [f (x^{k + 1})] - f (x^k) +\mathbb{E}_k [L (\xi^k) \| x^{k
  + 1} - x^k \|] + \frac{\tau}{2} \| \hat{x}^k - x^k \|^2 \label{eqn:a-2} \\
  \leq{} & L_f \mathbb{E}_k [\| x^{k + 1} - x^k \|] +\mathbb{E}_k [L (\xi^k) \|
  x^{k + 1} - x^k \|] + \frac{\tau}{2} \| \hat{x}^k - x^k \|^2 \label{eqn:a-3}
\end{align}

where \eqref{eqn:a-1} uses $L_f (\xi)$-Lipschitzness of $f_{x^k} (x, \xi)$; \eqref{eqn:a-2} uses quadratic bound from \ref{A5} and \eqref{eqn:a-3} applies $L_f$-Lipschitzness of $f (x)$ \cite{davis2019stochastic}. Re-arranging the terms, we have
\begin{align}
  & \frac{\gamma_k - \kappa}{2} \mathbb{E}_k [\| x^{k + 1} -
  \hat{x}^k \|^2] \nonumber\\
  \leq{} & \frac{\gamma_k - \rho + \tau}{2} \| \hat{x}^k - x^k \|^2 -
  \frac{\gamma_k - \rho}{2} \mathbb{E}_k [\| x^{k + 1} - x^k \|^2]
  +\mathbb{E}_k [(L (\xi^k) + L_f) \| x^{k + 1} - x^k \|] \nonumber\\
  \leq{} & \frac{\gamma_k - \rho + \tau}{2} \| \hat{x}^k - x^k \|^2 + \frac{2
  L_f^2}{\gamma_k - \rho} \label{eqn:a-4} \\
  ={} & \frac{\gamma_k - \kappa}{2} \| \hat{x}^k - x^k \|^2 - \frac{
  \rho - \tau - \kappa}{2} \| \hat{x}^k - x^k \|^2 + \frac{2 L_f^2}{\gamma_k
  - \rho} \nonumber
\end{align}
where \eqref{eqn:a-4} uses the relation $- \frac{a}{2} x^2 + b x \leq \frac{b^2}{2 a}$ and
$\mathbb{E}_{\xi} [(L (\xi) + L_f)^2] \leq 4 L_f^2$. Dividing both sides by
$\frac{\gamma_k  - \kappa}{2}$,
\[ \mathbb{E}_k [\| x^{k + 1} - \hat{x}^k \|^2] \leq \| \hat{x}^k - x^k \|^2 -
   \frac{\rho - \tau - \kappa}{\gamma_k - \kappa} \| \hat{x}^k -
   x^k \|^2 + \frac{4 L_f^2}{(\gamma_k - \rho) (\gamma_k - \kappa)} \]
and the potential function is reduced by
\begin{align}
  \mathbb{E}_k [\psi_{1 / \rho} (x^{k + 1})] ={} & \min_x  \{ f (x) +
  \omega (x) + \frac{\rho}{2} \| x - x^{k + 1} \|^2 \} \nonumber\\
  \leq{} & f (\hat{x}^k) + \omega (\hat{x}^k) + \frac{\rho}{2} \| \hat{x}^k -
  x^{k + 1} \|^2 \nonumber\\
  \leq{} & f (\hat{x}^k) + \omega (\hat{x}^k) + \frac{\rho}{2} \| \hat{x}^k -
  x^k \|^2 - \frac{\rho (\rho - \tau - \kappa)}{2 (\gamma_k -
  \kappa)} \| \hat{x}^k - x^k \|^2 + \frac{2 \rho L_f^2}{(\gamma_k - \rho)
  (\gamma_k - \kappa)} \nonumber\\
  ={} & \psi_{1 / \rho} (x^k) - \frac{\rho (\rho - \tau - \kappa)}{2
  (\gamma_k - \kappa)} \| \hat{x}^k - x^k \|^2 + \frac{2 \rho
  L_f^2}{(\gamma_k - \rho) (\gamma_k - \kappa)}, \nonumber
\end{align}
which completes the proof.

\subsection{Proof of Theorem \ref{thm:1}}

Given fixed stepsize $\gamma_k \equiv \gamma = \rho +  \kappa +
\alpha \sqrt{K}$, where we have, after telescoping, that
\begin{align}
  \frac{\rho (\rho - \tau - \kappa)}{2 (\gamma - \kappa)} \sum_{k =
  1}^K \Ebb[ \| \hat{x}^k - x^k \|^2] \leq{} & \psi_{1 / \rho} (x^1) -\mathbb{E}
  [\psi_{1 / \rho} (x^{K + 1})] + \frac{2 \rho L_f^2 K}{(\gamma - \rho)
  (\gamma - \kappa)} . \nonumber
\end{align}
Re-arranging the terms and summing over $k = 1, \ldots, K$, we have
\begin{align}
\min_{1\leq k \leq K}  \mathbb{E} [\| \nabla \psi_{1 / \rho} (x^{k}) \|^2] \leq{} & \frac{2
  \rho}{\rho - \tau - \kappa}  \bigg[ \frac{(\gamma - \kappa) D}{K}
  + \frac{2 \rho L_f^2}{\gamma - \rho} \bigg] \nonumber\\
  \leq{} & \frac{2 \rho}{\rho - \tau - \kappa} \bigg[ \frac{\rho D}{K} +
  \frac{\alpha D}{\sqrt{K}} + \frac{2 \rho L_f^2}{\alpha \sqrt{K}} \bigg],
  \nonumber
\end{align}
where $D = \psi_{1/\rho}(x^1) -  \inf_x \psi(x) \geq \psi_{1/\rho}(x^1) -  \Ebb[\psi_{1/\rho}(x^K)]$ and this completes the proof.

\section{Proof of results in Section \ref{sec:known-nonlip}}
\label{app:known-nolip}
For brevity of notation, in the proof we define 
\begin{equation}
	\mathsf{G}_k \assign \mathcal{G} (\| x^k \|) \label{eqn:b-1}
\end{equation}

and use them interchangeably in this section. We also note that $\mathsf{G}_k$ is a random variable whose randomness comes from samples from previous iterations $\xi^1,\ldots,\xi^{k-1}$.

\subsection{Proof of Lemma \ref{lem:2}}
Firstly, we still use optimality condition to get, for a given $\xi^k$, that
\[ f_{x^k} (x^{k + 1}, \xi^k) + \omega (x^{k + 1}) + \frac{\gamma_k}{2} \|
   x^{k + 1} - x^k \|^2 \leq f_{x^k} (\hat{x}^k, \xi^k) + \omega (\hat{x}^k) +
   \frac{\gamma_k}{2} \| \hat{x}^k - x^k \|^2 - \frac{\gamma_k - \kappa}{2} \| x^{k +
   1} - \hat{x}^k \|^2 \]
\[ f (\hat{x}^k) + \omega (\hat{x}^k) + \frac{\rho}{2} \| \hat{x}^k - x^k \|^2
   \leq f (x^k) + \omega (x^k) \]
Summing over the above two relations, we deduce that
\begin{align}
  & \frac{\gamma_k}{2} \| x^{k + 1} - x^k \|^2 - \frac{\gamma_k - \rho}{2} \|
  \hat{x}^k - x^k \|^2 + \frac{\gamma_k - \kappa}{2} \| x^{k + 1} - \hat{x}^k \|^2
  \nonumber\\
  \leq{} & f (x^k) - f_{x^k} (x^{k + 1}, \xi^k) + f_{x^k} (\hat{x}^k, \xi^k) - f
  (\hat{x}^k) + L_{\omega} \| x^{k + 1} - x^k \| \label{eqn:b-2-0}\\
  \leq{} & f (x^k) - f (x^k, \xi^k) + f_{x^k} (\hat{x}^k, \xi^k) - f (\hat{x}^k)
  + (\mathsf{G}_k L_f (\xi^k) + L_{\omega}) \| x^{k + 1} - x^k \|, \label{eqn:b-2-1}
\end{align}
where \eqref{eqn:b-2-0} applies $L_\omega$-Lipschitz continuity of $\omega(x)$; \eqref{eqn:b-2-1} applies \ref{C1}.
Dividing both sides by $\frac{\gamma_k - \kappa}{2}$ and re-arranging the terms, we have
\begin{align}
  \| x^{k + 1} - \hat{x}^k \|^2 
  \leq{} & \frac{\gamma_k - \rho}{\gamma_k - \kappa} \| \hat{x}^k - x^k \|^2 -
  \frac{\gamma_k}{\gamma_k - \kappa} \| x^{k + 1} - x^k \|^2 \nonumber\\
   & + \frac{2}{\gamma_k - \kappa} [f (x^k) - f (x^k, \xi^k) + f_{x^k}
  (\hat{x}^k, \xi^k) - f (\hat{x}^k) + (\mathsf{G}_k L_f (\xi^k) + L_{\omega})
  \| x^{k + 1} - x^k \|] \nonumber
\end{align}

Next conditioned on $\xi^1,\ldots,\xi^{k-1}$, taking expectation with respect to $\xi^{k}$, and recalling that $\mathcal{G} (\| x^k \|)$, and therefore $\gamma_k$ is fixed given
$\xi^1, \ldots, \xi^{k-1}$, we have
\begin{align}
  &  \mathbb{E}_k [\| x^{k + 1} - \hat{x}^k \|^2]
  \nonumber\\
  \leq{} & \frac{\gamma_k - \rho + \tau}{\gamma_k - \kappa} \| \hat{x}^k - x^k \|^2
  + \frac{2}{\gamma_k - \kappa} \mathbb{E}_k [ (\mathsf{G}_k L_f (\xi^k) + L_{\omega}) \| x^{k + 1} - x^k \| -
  \frac{\gamma_k}{2} \| x^{k + 1} - x^k \|^2 ] \label{eqn:b-2-2} \\
  \leq{} & \frac{\gamma_k - \rho + \tau}{\gamma_k - \kappa} \| \hat{x}^k - x^k \|^2 + \frac{(\mathsf{G}_k
  L_f + L_{\omega})^2}{\gamma_k(\gamma_k - \kappa)} \label{eqn:b-2} \\
  ={} & \| \hat{x}^k - x^k \|^2 - \frac{\rho - \tau -\kappa}{\gamma_k - \kappa} \|
  \hat{x}^k - x^k \|^2 + \frac{(\mathsf{G}_k L_f + L_{\omega})^2}{\gamma_k(\gamma_k - \kappa)},
  \nonumber
\end{align}
where \eqref{eqn:b-2-2} applies $f(x^k) - \Ebb_k[f(x^k, \xi^k)] = 0, \Ebb_k[f_{x^k}(\hat{x}^k, \xi^k)] - f(\hat{x}^k) \leq \frac{\tau}{2} \|x^k - \hat{x}^k\|^2$;
  \eqref{eqn:b-2} applies the relation $- \frac{a}{2} x^2 + b x \leq \frac{b^2}{2 a}$ and $\mathbb{E} [L_f (\xi)]^2 \leq \mathbb{E} [L_f
(\xi)^2] \leq L_f^2$. 
Now we can deduce reduction of the potential function by
\begin{align}
  \mathbb{E}_k [\psi_{1 / \rho} (x^{k + 1})] ={} & \min_x  \{ f (x) +
  \omega (x) + \frac{\rho}{2} \| x - x^{k + 1} \|^2 \} \nonumber\\
  \leq{} & f (\hat{x}^k) + \omega (\hat{x}^k) + \frac{\rho}{2} \| \hat{x}^k -
  x^{k + 1} \|^2 \nonumber\\
  \leq{} & f (\hat{x}^k) + \omega (\hat{x}^k) + \frac{\rho}{2} \| \hat{x}^k -
  x^k \|^2 - \frac{\rho (\rho - \tau - \kappa)}{2 (\gamma_k - \kappa)} \| \hat{x}^k - x^k \|^2 +
  \frac{\rho (\mathsf{G}_k L_f + L_{\omega})^2}{2 \gamma_k(\gamma_k - \kappa)} \nonumber\\
  ={} & \psi_{1 / \rho} (x^k) - \frac{\rho (\rho - \tau - \kappa)}{2 (\gamma_k - \kappa)} \|
  \hat{x}^k - x^k \|^2 + \frac{\rho (\mathsf{G}_k L_f + L_{\omega})^2}{2 \gamma_k (\gamma_k - \kappa)}
  \nonumber
\end{align}

and this completes the proof.

\subsection{Proof of Theorem \ref{thm:2}}

First we introduce the following lemma.

\begin{lem}[Robbins-Siegmund \cite{robbins1971convergence}] \label{lem:robbins}
  Let $A_k, B_k, C_k$ and $V_k$ be nonnegative
  random variables adapted to the filtration $\mathcal{F}_k$ and satisfying
  $\mathbb{E} [V_{k + 1} |\mathcal{F}_k] \leq (1 + A_k) V_k + B_k - C_k$.
  Then on the event $\{ \sum_{k = 1}^{\infty} A_k < \infty, \sum_{k =
  1}^{\infty} B_k < \infty \}$, there is a random variable $V_{\infty}$
  such that $V_k \overset{a.s.}{\longrightarrow} V_{\infty}$ and $\sum_{k =
  0}^{\infty} C_k < \infty$ almost surely.
\end{lem}

Now we get down to the proof. Recall that in \textbf{Lemma \ref{lem:2}} we have shown that
\[ \mathbb{E}_k [\psi_{1 / \rho} (x^{k + 1})] \leq \psi_{1 / \rho} (x^k) -
   \frac{\rho (\rho - \tau - \kappa)}{2 (\gamma_k - \kappa)} \| \hat{x}^k - x^k \|^2 + \frac{\rho
   (\mathsf{G}_k L_f + L_{\omega})^2}{2 \gamma_k(\gamma_k - \kappa)} \]
and we can bound 
\begin{equation}
\frac{\mathsf{G}_k L_f + L_{\omega}}{\gamma_k - \kappa} = \frac{\mathsf{G}_k L_f +
L_{\omega}}{\rho + \tau  + k^{\zeta} (\mathsf{G}_k + 1)} \leq \frac{\mathsf{G}_k L_f +
L_{\omega}}{ k^{\zeta} (\mathsf{G}_k + 1)} \leq \frac{L_f + L_{\omega}}{
k^{\zeta}}	\label{eqn:bound-error}
\end{equation}
to get
\[ \mathbb{E}_k [\psi_{1 / \rho} (x^{k + 1}) + \Lambda] \leq [\psi_{1 / \rho}
   (x^k) + \Lambda] - \frac{\rho (\rho - \kappa - \tau)}{2 (\gamma_k - \kappa)} \| \hat{x}^k - x^k
   \|^2 + \frac{\rho }{2 k^{2 \zeta}}(L_f + L_{\omega})^2 \]
Then we invoke \textbf{Lemma \ref{lem:robbins}}, plugging in the relation 
\begin{equation}
A_k = 0, \quad B_k = \frac{\rho (L_f + L_{\omega})^2}{2  k^{2 \zeta}}, \quad C_k = \frac{\rho (\rho - \kappa - \tau)}{2 (\gamma_k - \kappa)} \| \hat{x}^k - x^k \|^2, \quad  V_k = \psi_{1 / \rho} (x^k) + \Lambda \geq 0
\end{equation}
Then
with $\zeta \in ( \frac{1}{2}, 1)$ \ $\sum_{k = 1}^{\infty} B_k =
\frac{\rho (L_f + L_{\omega})^2}{2 k^{2 \zeta}} < \infty$ we know that  $\{ \psi_{1 / \rho} (x^k) + \Lambda \} \rightarrow \psi_{1 / \rho} (x^{\infty}) + \Lambda < \infty$ and that $\sum_{k =
1}^{\infty} \frac{\rho (\rho - \kappa - \tau)}{2 (\gamma_k - \kappa)} \| \hat{x}^k - x^k \|^2 <
\infty$. By \ref{A7}, $\| x^k \|$ is bounded with probability
1 and $\mathsf{G}_k$ is bounded almost surely. Finally $\sum_{k = 1}^{\infty}
\frac{1}{\gamma_k - \kappa} = \infty \Rightarrow \inf_{j \leq k} \| \hat{x}^j - x^j \| \rightarrow 0$
almost surely and this completes the proof since $\| \hat{x}^k - x^k \|=\rho^{-1} \|\nabla \psi_{1/\rho}(x^k) \|$.

\subsection{Proof of Lemma \ref{lem:2.22} and \ref{lem:2.2}}

Following \eqref{eqn:bound-error}, we first we bound the error of potential reduction by $\frac{\rho (\mathsf{G}_k L_f + L_{\omega})^2}{2 \gamma_k(\gamma_k - \kappa)} \leq \frac{\rho (L_f + L_{\omega})^2}{\alpha^2
K}$. 

Then a telescopic sum gives, for all $2 \leq k \leq K$ that
\[ \mathbb{E} [\psi_{1 / \rho} (x^k) + \Lambda] \leq \psi_{1 / \rho} (x^1) +
   \Lambda + \frac{\rho(L_f + L_{\omega})^2}{\alpha^2} =: \Delta. \]
Here $\Delta$ is a constant that only depends on the initialization of the algorithm.
Next we consider one step of the algorithm
\[ f_{x^k} (x^{k + 1}, \xi^k) + \omega (x^{k + 1}) + \frac{\gamma_k}{2} \|
   x^{k + 1} - x^k \|^2 \leq f_{x^k} (x^k, \xi^k) + \omega (x^k) \]
and a re-arrangement gives
\begin{align}
 \frac{\gamma_k}{2} \| x^{k + 1} - x^k \|^2 \nonumber \leq{} & f_{x^k} (x^k, \xi^k) - f_{x^k} (x^{k + 1}, \xi^k) + \omega (x^k) -
  \omega (x^{k + 1}) \nonumber\\
  \leq{} & (L_f (\xi^k) \mathsf{G}_k + L_{\omega}) \| x^{k + 1} - x^k \|, \nonumber
\end{align}
where we use \ref{C1} and Lipschitz continuity of $\omega(x)$. Dividing both sides by $\| x^{k + 1} - x^k \|$, we have
\[ \| x^{k + 1} - x^k \| \leq \frac{2(L_f (\xi^k) \mathsf{G}_k + L_{\omega})}{\gamma_k} \leq
   \frac{2(L_f (\xi^k) \mathsf{G}_k + L_{\omega})}{\alpha (\mathsf{G}_k + 1) \sqrt{K}} . \]
Conditioned on $\xi^1, \ldots, \xi^{k-1}$ and taking expectation with respect to $\xi^k$, we get
\[ \mathbb{E}_k [\| x^{k + 1} - x^k \|] \leq \frac{2\mathbb{E}_{\xi^k} [L_f (\xi^k)] \mathsf{G}_k
   + 2L_{\omega}}{\alpha (\mathsf{G}_k + 1) \sqrt{K}} \leq \frac{2L_f \mathsf{G}_k +
   2L_{\omega}}{\alpha (\mathsf{G}_k + 1) \sqrt{K}} \leq \frac{2(L_f + L_{\omega})}{\alpha
   \sqrt{K}} . \]
   This completes the proof of \textbf{Lemma \ref{lem:2.22}}.  By Markov's inequality, we know that, for any $2\leq k\leq K$, the following bound holds
\[ \mathbb{P}_{\xi^k \sim \Xi} \Big\{ \| x^{k + 1} - x^k \| \leq \frac{4
   (L_f + L_{\omega})}{\alpha \sqrt{K}} \big| \xi^1, \ldots, \xi^{k - 1}  \Big\}
   \geq \frac{1}{2} \]
   
and without loss of generality we let $\mathsf{Z} = \frac{4
   (L_f + L_{\omega})}{\alpha \sqrt{K}}$, and clearly $\mathsf{Z} = \mathcal{O}(1/\sqrt{K}) = \mathcal{O}(1)$.
\par\noindent\rule{\textwidth}{0.6pt}\vspace{5pt}
This relation says ``it's likely that $x^k$ and $x^{k + 1}$ are
 close'', and we leverage this intuition to derive a tail-bound on $\| x^k
\|$. So far, we have the following properties in hand:
\begin{enumerate}
  \item $\mathbb{E} [\psi_{1 / \rho} (x^k)]$ is bounded by a constant $\Delta - \Lambda$ for all $k$
  
  \item If $\| x^k \| \geq \mathsf{B}_{v}$, then $\psi_{1 / \rho} (x^k) \geq v$
  
  \item If $\| x^k \| \geq \mathsf{B}_{v}$, it's likely that $\| x^{k + 1} \| \geq \mathsf{B}_{v}
  -\mathcal{O} ( 1 / \sqrt{K} )$.
\end{enumerate}
Our reasoning is as follows: given large $a > 1$, conditioned on $\| x^k \| \geq
\mathsf{B}_{a \Delta}$, then it is likely that $\| x^{k + 1} \| \approx \mathsf{B}_{a \Delta}$ since $\|x^{k+1}-x^k\|$ is likely to be small. And it implies $\Ebb[\psi_{1 / \rho} (x^{k + 1})] \geq a \Delta > \Delta$. However, we
know that $\mathbb{E} [\psi_{1 / \rho} (x^{k+1})] \leq \Delta$, and this will therefore reversely bound the probability that $\| x^k \| \geq \mathsf{B}_{a \Delta} +\mathcal{O}
( 1 / \sqrt{K} )$. We formalize the proof as follows.
\par\noindent\rule{\textwidth}{0.6pt}
First recall that by \ref{A7}, $\| x^k \| \geq \mathsf{B}_{v}$ implies $\psi_{1 /
\rho} (x^k) \geq v$. Taking $v = a \Delta, a > 1$, we have $\| x^k \| \geq \mathsf{B}_{a \Delta} \Rightarrow \psi_{1 / \rho} (x^k) \geq a v$. Now we consider the
event $\| x^k \| \geq \mathsf{B}_{a \Delta} + \mathsf{Z}$ and apply law of total expectation to get
\begin{align}
  \Delta \geq{} & \mathbb{E} [\psi_{1 / \rho} (x^{k + 1}) + \Lambda] \nonumber\\
  ={} & \mathbb{E} [ \psi_{1 / \rho} (x^{k + 1}) + \Lambda | \| x^k \| \geq
  \mathsf{B}_{a \Delta} + \mathsf{Z} ] \cdot \mathbb{P}
  \{ \| x^k \| \geq \mathsf{B}_{a \Delta} + \mathsf{Z} \} \nonumber\\
  & +\mathbb{E} [ \psi_{1 / \rho} (x^{k + 1}) + \Lambda | \| x^k \| \leq
  \mathsf{B}_{a \Delta} + \mathsf{Z} ] \cdot \mathbb{P}
  \{ \| x^k \| \leq \mathsf{B}_{a \Delta} + \mathsf{Z} \} \nonumber\\
  \geq{} & \mathbb{E} [ \psi_{1 / \rho} (x^{k + 1}) + \Lambda | \| x^k \|
  \geq{} \mathsf{B}_{a \Delta} + \mathsf{Z} ] \cdot
  \mathbb{P} \{ \| x^k \| \geq \mathsf{B}_{a \Delta} + \mathsf{Z} \}, \label{eqn:b-3}
\end{align}

where \eqref{eqn:b-3} uses $\psi_{1 / \rho} (x) + \Lambda \geq 0$ for all $x$.
Next we consider the expectation 
\begin{equation}
\mathbb{E} \big[ \psi_{1 / \rho} (x^{k +
1}) + \Lambda | \| x^k \| \geq \mathsf{B}_{a \Delta} + \mathsf{Z}\big],	
\end{equation}
and successively deduce that
\begin{align}
  & \mathbb{E} [ \psi_{1 / \rho} (x^{k + 1}) + \Lambda | \| x^k \| \geq
  \mathsf{B}_{a \Delta} + \mathsf{Z} ]
  \nonumber\\
  ={} & \mathbb{E} [ \psi_{1 / \rho} (x^{k + 1}) + \Lambda | \| x^k \|
  \geq{} \mathsf{B}_{a \Delta} + \mathsf{Z}, \| x^{k +
  1} - x^k \| \leq \mathsf{Z} ] \cdot
  \mathbb{P} \{ \| x^{k + 1} - x^k \| \leq \mathsf{Z}\}\nonumber\\
  & +\mathbb{E} [ \psi_{1 / \rho} (x^{k + 1}) + \Lambda | \| x^k \| \geq
  \mathsf{B}_{a \Delta} + \mathsf{Z}, \| x^{k + 1} -
  x^k \| \geq \mathsf{Z} ] \cdot \mathbb{P}
  \{ \| x^{k + 1} - x^k \| \geq \mathsf{Z} \} \nonumber\\
  \geq{} & \frac{1}{2}\mathbb{E} [ \psi_{1 / \rho} (x^{k + 1}) + \Lambda | \| x^k \|
  \geq{} \mathsf{B}_{a \Delta} + \mathsf{Z}, \| x^{k +
  1} - x^k \| \leq \mathsf{Z} ] 
  \label{eqn:b-4-0}\\
  \geq{} & \frac{1}{2}(a \Delta + \Lambda) , \label{eqn:b-4}
\end{align}

where \eqref{eqn:b-4-0} is by Markov's inequality $\mathbb{P} \{ \| x^{k +
1} - x^k \| \leq \mathsf{Z} \} \geq
0.5$ and \eqref{eqn:b-4} uses the triangle inequality
\[ \| x^{k + 1} \| = \| x^{k + 1} - x^k + x^k \| \geq \| x^k \| - \| x^{k + 1}
   - x^k \| \geq \| x^k \| - \mathsf{Z}
   \geq \mathsf{B}_{a \Delta} \]
and we recall that $\|x^{k + 1}\| \geq \mathsf{B}_{a \Delta}$ implies $\psi_{1 / \rho} (x^{k + 1}) \geq a\Delta$. Chaining the above inequalities, we arrive at
\[ \Delta \geq \big( \tfrac{a \Delta + \Lambda}{2} \big) \cdot \mathbb{P} \{
   \| x^k \| \geq \mathsf{B}_{a \Delta} + \mathsf{Z}
   \} \]
Dividing both sides by $\frac{a \Delta + \Lambda}{2}$ gives the desired tail
bound
\[ \mathbb{P} \Big\{ \| x^k \| \geq \mathsf{B}_{a \Delta} + \mathsf{Z} \Big\} \leq \frac{2 \Delta}{a \Delta +
   \Lambda} . \]
Since $\Lambda \geq 0$, the bound is nontrivial when $a > 2$, and this
completes the proof.

\subsection{Proof of Theorem \ref{thm:2.2}}

Given $\delta \in (0, 1/4)$, we know, from \textbf{Lemma \ref{lem:2.2}} that for $2\leq k \leq K$,
\[ \mathbb{P} \left\{ \| x^k \| \geq \mathsf{B}_{\delta^{- 1} \Delta} + \mathsf{Z} \right\} \leq \frac{2 \Delta}{\delta^{- 1}
   \Delta + \Lambda} \leq 2 \delta . \]
Then denote $I_k =\mathbb{I} \big\{ \| x^k \| \leq \mathsf{B}_{\delta^{- 1} \Delta} +
\mathsf{Z} \big\}$, and we have $\sum_{k =
1}^K \mathbb{E} [1 - I_k] \leq 2 \delta K$. By Markov's inequality, given $p
\in (2\delta, 1)$,
\[ \mathbb{P} \left\{ \sum_{k = 1}^K (1 - I_k) \geq \frac{2\delta K}{p}
   \right\} \leq \frac{2 \delta K}{2 p^{- 1} \delta K} = p \]
and a re-arrangement gives $\mathbb{P} \big\{ \sum_{k = 1}^K I_k \geq K (1 - 2 p^{- 1} \delta)
\big\} \geq 1 - p$. 
Now we telescope \tmtextbf{Lemma \ref{lem:2}} and get
\begin{align}
  \sum_{k = 1}^K \mathbb{E} \left[ \frac{\rho (\rho - \kappa- \tau)}{2 \gamma_k} \|
  \hat{x}^k - x^k \|^2 \right] \leq{} & \psi_{1 / \rho} (x^1) -\mathbb{E}
  [\psi_{1 / \rho} (x^{K + 1})] + \frac{\rho }{2 \alpha^2}(L_f + L_{\omega})^2
  \label{eqn:b-5}\\
  \leq{} & \psi_{1 / \rho} (x^1) - \inf_x \psi (x) + \frac{\rho}{2 \alpha^2} (L_f +
  L_{\omega})^2, \label{eqn:b-5-1}\\
  ={} & D + \frac{\rho }{2 \alpha^2} (L_f + L_{\omega})^2 \nonumber
\end{align}

where \eqref{eqn:b-5} uses the previously established bound $\Ebb_k\Big[\frac{\rho (\mathsf{G}_k L_f +
L_{\omega})^2}{2 \gamma_k(\gamma_k - \kappa)}\Big] \leq \frac{\rho (L_f + L_{\omega})^2}{\alpha^2
K}$. \eqref{eqn:b-5-1} uses the fact that $\psi_{1 / \rho} (x^{K + 1}) \geq \inf_x \psi (x)$.
Next we notice, conditioned on the event $\sum_{k = 1}^K I_k \geq K (1 - 2
p^{- 1} \delta)$, defining $\mathsf{G}_\delta \assign \max_z \mathcal{G} (z), z \leq
\mathsf{B}_{\delta^{- 1} \Delta} + \mathsf{Z}$, that
\begin{align}
  & \sum_{k = 1}^K \mathbb{E} \left[ \frac{\rho (\rho - \kappa - \tau)}{2 \gamma_k} \|
  \hat{x}^k - x^k \|^2 \right] \nonumber\\
  ={} & {\sum_{k \in \{ j : I_j = 0 \}}}  \mathbb{E} \left[ \frac{\rho (\rho - \kappa -
  \tau)}{2 \gamma_k} \| \hat{x}^k - x^k \|^2 \right] + {\sum_{k \in \{ j : I_j
  ={} 1 \}}}  \mathbb{E} \left[ \frac{\rho (\rho - \kappa - \tau)}{2 \gamma_k} \|
  \hat{x}^k - x^k \|^2 \right] \nonumber\\
  \geq{} & {\sum_{k \in \{ j : I_j = 1 \}}}  \mathbb{E} \left[ \frac{\rho (\rho
  - \kappa - \tau)}{2 \gamma_k} \| \hat{x}^k - x^k \|^2 \right] \nonumber\\
  \geq{} & {\sum_{k \in \{ k : I_j = 1 \}}} \frac{\rho (\rho
  - \kappa - \tau)}{2 ( \rho + \kappa + \tau + \alpha (\mathsf{G}_\delta + 1) \sqrt{K} )} \mathbb{E}[ \|
  \hat{x}^k - x^k \|^2 ] \label{eqn:b-6} \\
  \geq{} & \frac{\rho (\rho - \kappa - \tau) (1 - 2 p^{- 1} \delta) K}{2 ( \rho +
  \kappa + \tau + \alpha (\mathsf{G}_\delta + 1) \sqrt{K} )} \min_{k \in \{ j : I_j = 1 \}}
  \mathbb{E} [\| \hat{x}^k - x^k \|^2], \label{eqn:b-7}
\end{align}
where \eqref{eqn:b-6} applies the fact that if $I_j=1$, $\|x^j\|\leq \mathsf{G}_\delta$; \eqref{eqn:b-7} utilizes the fact that $|\{k:I_j = 1\}| \geq K(1-2p^{-1}\delta)$. Putting the inequalities together, we have
\begin{align}
  \min_{k \in \{ j : I_j = 1 \}} \mathbb{E} [\| \hat{x}^k - x^k \|^2] \leq{} &
  \frac{2 ( \rho + \kappa + \tau + \alpha (\mathsf{G}_\delta + 1) \sqrt{K} )}{\rho (\rho
 - \kappa - \tau) (1 - 2 p^{- 1} \delta) K} \left[ D + \frac{\rho (L_f +
  L_{\omega})^2}{2 \alpha^2} \right] \nonumber\\
  ={} & \frac{2 ( D + \frac{\rho (L_f + L_{\omega})^2}{2 \alpha^2}
  )}{\rho (\rho - \kappa - \tau) (1 - 2 p^{- 1} \delta)} \left[ \frac{\rho +
  \kappa + \tau}{K} + \frac{\alpha (\mathsf{G}_\delta + 1)}{\sqrt{K}} \right], \nonumber
\end{align}

Recalling that $\| \hat{x}^k - x^k \|=\rho^{-1} \|\nabla \psi_{1/\rho}(x^k) \|$, at least with probability $1-p$,
\begin{align}
  \min_{1 \leq k \leq K} \mathbb{E} [\| \nabla \psi_{1 / \rho} (x^k) \|^2]
  \leq{} & \min_{k \in \{ k : I_k = 1 \}} \mathbb{E} [\| \nabla \psi_{1 / \rho}
  (x^k) \|^2] \nonumber\\
  \leq{} & \frac{p}{p - 2 \delta} \cdot \frac{2 \rho}{\rho - \tau - \kappa}
  \left[ D + \frac{\rho (L_f + L_{\omega})^2}{2 \alpha^2} \right] \left(
  \frac{\rho + \tau + \kappa}{K} + \frac{\alpha (\mathsf{G}_{\delta} +
  1)}{\sqrt{K}} \right) \nonumber
\end{align}

and this completes the proof.

\section{Proof of results in Section \ref{sec:unknown-nolip}}
\label{app:unknown-nolip}
For brevity of expression, we define $L_f^k := \textup{\textsf{Lip}} (x^k, \xi^k), L_f' := \textup{\textsf{Lip}} (x^k, \xi')$ ($k$ is hidden when clear from the context) and use them
interchangeably.

\subsection{Auxiliary lemmas}
\begin{lem}
  \label{lem:3.3}Given independent nonnegative random variables $X$ and $Y$. If $\mathbb{E}_{X, Y} [| X - Y |^2] \leq \sigma^2,$ then
  \[ \mathbb{E}_{X, Y} \left[ \tfrac{X^2}{\max \{ Y^2, \alpha^2 \}} \right]
     \leq \big( \tfrac{\sigma + \alpha}{\alpha} \big)^2, \quad
  \mathbb{E}_{X, Y} \left[ \tfrac{X}{\max \{ Y^2, \alpha^2 \}} \right] \leq
     \tfrac{\sigma}{\alpha^2} + \tfrac{1}{\alpha}, \quad \mathbb{E}_{X, Y} \left[ \tfrac{X}{\max \{ Y, \alpha \}} \right] \leq
     \tfrac{\sigma}{\alpha} + 1 \]
     where $\alpha > 0$.
\end{lem}
\begin{proof}
For the first relation, we successively deduce that
\begin{align}
  \mathbb{E}_{X, Y} \left[ \tfrac{X^2}{\max \{ Y^2, \alpha^2 \}} \right] ={} &
  \mathbb{E}_{X, Y} \left[ \tfrac{(X - Y + Y)^2}{\max \{ Y^2, \alpha^2 \}}
  \right] \nonumber\\
    ={} & \mathbb{E}_{X, Y} \left[  \tfrac{(X - Y)^2 + Y^2 + 2 Y (X -
  Y)}{\max \{ Y^2, \alpha^2 \}} \right] \nonumber\\
  ={} & \mathbb{E}_{X, Y} \left[ \tfrac{(X - Y)^2}{\max \{ Y^2, \alpha^2 \}}
  \right] +\mathbb{E}_Y \left[ \tfrac{Y^2}{\max \{ Y^2, \alpha^2 \}} \right]
  +\mathbb{E}_{X, Y} \left[ \tfrac{2 Y (X - Y)}{\max \{ Y^2, \alpha^2 \}}
  \right] \nonumber\\
  \leq{} & \tfrac{\sigma^2}{\alpha^2} + 1 +\mathbb{E}_{X, Y} \left[ \tfrac{2 Y | X
  - Y |}{\max \{ Y, \alpha \} \cdot \max \{ Y, \alpha \}} \right]
  \label{eqn:c-1} \\
  \leq{} & \tfrac{\sigma^2}{\alpha^2} + \tfrac{2 \sigma}{\alpha} + 1 = \big(
  \tfrac{\sigma + \alpha}{\alpha} \big)^2, \nonumber
\end{align}
where \eqref{eqn:c-1} uses $\tfrac{a}{\max\{b, c\}} \leq \tfrac{a}{c}$ and the last inequality applies 

\begin{equation}
\tfrac{2 Y | X
  - Y |}{\max \{ Y, \alpha \} \cdot \max \{ Y, \alpha \}} = \tfrac{2 Y }{\max \{ Y, \alpha \}}\cdot\tfrac{| X
  - Y |}{\max \{ Y, \alpha \}} \leq 2\cdot \tfrac{|X - Y|}{\alpha}.
\end{equation}
Similarly we can deduce that
\begin{align}
  \mathbb{E}_{X, Y} \left[ \tfrac{X}{\max \{ Y^2, \alpha^2 \}} \right] \leq{}
  \mathbb{E}_{X, Y} \left[ \tfrac{| X - Y |}{\max \{ Y^2, \alpha^2 \}} \right]
  +\mathbb{E}_Y \left[ \tfrac{Y}{\max \{ Y^2, \alpha^2 \}} \right] 
  \leq{} & \tfrac{\sigma}{\alpha^2} + \tfrac{1}{\alpha}, \nonumber
\end{align}

\[ \mathbb{E}_{X, Y} \left[ \tfrac{X}{\max \{ Y, \alpha \}} \right] \leq
   \mathbb{E}_{X, Y} \left[ \tfrac{| X - Y |}{\max \{ Y, \alpha \}} \right]
   +\mathbb{E}_{X, Y} \left[ \tfrac{Y}{\max \{ Y, \alpha \}} \right] \leq
   \tfrac{\sigma}{\alpha} + 1, \]
which completes the proof.	
\end{proof}

\subsection{Proof of Lemma \ref{lem:3.4}}
By the optimality condition, we have
\[ f_{x^k} (x^{k + 1}, \xi^k) + \omega (x^{k + 1}) + \frac{\gamma_k}{2} \|
   x^{k + 1} - x^k \|^2 \leq f_{x^k} (\hat{x}^k, \xi^k) + \omega (\hat{x}^k) +
   \frac{\gamma_k}{2} \| \hat{x}^k - x^k \|^2 - \frac{\gamma_k - \kappa}{2} \| x^{k +
   1} - \hat{x}^k \|^2 \]
\[ f (\hat{x}^k) + \omega (\hat{x}^k) + \frac{\rho}{2} \| \hat{x}^k - x^k \|^2
   \leq f (x^k) + \omega (x^k) \]
and summation over the two relations gives
\begin{align}
  & \frac{\gamma_k}{2} \| x^{k + 1} - x^k \|^2 - \frac{\gamma_k - \rho}{2} \|
  \hat{x}^k - x^k \|^2 + \frac{\gamma_k - \kappa}{2} \| x^{k + 1} - \hat{x}^k \|^2
  \nonumber\\
  \leq{} & f (x^k) - f_{x^k} (x^{k + 1}, \xi^k) + f_{x^k} (\hat{x}^k, \xi^k) - f
  (\hat{x}^k) + L_{\omega} \| x^{k + 1} - x^k \| \nonumber\\
  \leq{} & f (x^k) - f (x^k, \xi^k) + f_{x^k} (\hat{x}^k, \xi^k) - f (\hat{x}^k)
  + (L_f^k + L_{\omega}) \| x^{k + 1} - x^k \|,\label{eqn:c-2}
\end{align}
where \eqref{eqn:c-2} applies \ref{D1}.
Fixing $\xi'$, we divide both sides by $\frac{\gamma_k - \kappa}{2}$ and get
\begin{align}
  \| x^{k + 1} - \hat{x}^k \|^2 \leq{} & \frac{\gamma_k - \rho}{\gamma_k -
  \kappa} \| \hat{x}^k - x^k \|^2 - \frac{\gamma_k}{\gamma_k - \kappa} \| x^{k
  + 1} - x^k \|^2 \nonumber\\
  & + \frac{2}{\gamma_k - \kappa} [f (x^k) - f (x^k, \xi^k) + f_{x^k}
  (\hat{x}^k, \xi^k) - f (\hat{x}^k) + (L_f^k + L_{\omega}) \| x^{k + 1} - x^k
  \|] \nonumber
\end{align}
Conditioned on $\xi^1, \ldots, \xi^{k-1}$ and taking expectation with respect to $\xi^k$, we have
\begin{align}
  & \mathbb{E}_k [\| x^{k + 1} - \hat{x}^k \|^2]
  \nonumber\\
  \leq{} & \frac{\gamma_k - \rho + \tau}{\gamma_k - \kappa} \| \hat{x}^k - x^k \|^2 + \frac{2}{\gamma_k - \kappa}\mathbb{E}_k
  [  (L_f^k + L_{\omega}) \| x^{k + 1} - x^k
  \| -
  \frac{\gamma_k}{2} \| x^{k + 1} - x^k \|^2] \label{eqn:c-3-1} \\
  \leq{} & \frac{\gamma_k - \rho + \tau}{\gamma_k - \kappa} \| \hat{x}^k - x^k \|^2
  +\mathbb{E}_k \Big[ \frac{1}{\gamma_k (\gamma_k - \kappa)} (L_f^k + L_{\omega})^2\Big] \label{eqn:c-3}
\end{align}
where \eqref{eqn:c-3-1} again uses $f(x^k) - \Ebb_k[f(x^k, \xi^k)] = 0, \Ebb_k[f_{x^k}(\hat{x}^k, \xi^k)] - f(\hat{x}^k) \leq \frac{\tau}{2} \|x^k - \hat{x}^k\|^2$; \eqref{eqn:c-3} uses the relation $- \frac{a}{2} x^2 + b x \leq
\frac{b^2}{2 a}$. Now we have
\[ \mathbb{E}_k [\| x^{k + 1} - \hat{x}^k \|^2] \leq \| \hat{x}^k - x^k \|^2 -
   \frac{\rho - \tau - \kappa}{\gamma_k - \kappa} \| \hat{x}^k - x^k \|^2 +\mathbb{E}_k \Big[
   \frac{1}{\gamma_k(\gamma_k - \kappa)}(L_f^k + L_{\omega})^2 \Big] . \]
In view of our potential function, we have
\begin{align}
  \mathbb{E}_k [\psi_{1 / \rho} (x^{k + 1})] ={} & \min_x  \{ f (x) +
  \omega (x) + \frac{\rho}{2} \| x - x^{k + 1} \|^2 \} \nonumber\\
  \leq{} & f (\hat{x}^k) + \omega (\hat{x}^k) + \frac{\rho}{2} \| \hat{x}^k -
  x^{k + 1} \|^2 \nonumber\\
  \leq{} & f (\hat{x}^k) + \omega (\hat{x}^k) + \frac{\rho}{2} \| \hat{x}^k -
  x^k \|^2 - \frac{\rho (\rho - \tau - \kappa)}{2 (\gamma_k -\kappa)} \| \hat{x}^k - x^k \|^2
  +\mathbb{E}_k \Big[ \frac{\rho}{2 \gamma_k (\gamma_k - \kappa)}(L_f^k + L_{\omega})^2
  \Big] \nonumber\\
  ={} & \psi_{1 / \rho} (x^k) - \frac{\rho (\rho - \tau - \kappa)}{2 (\gamma_k -\kappa)} \|
  \hat{x}^k - x^k \|^2 +\mathbb{E}_k \Big[ \frac{\rho}{2 \gamma_k(\gamma_k - \kappa)}(L_f^k +
  L_{\omega})^2 \Big] \nonumber
\end{align}

and this completes the proof.

\subsection{Proof of Theorem \ref{thm:3}}
Our reasoning is the same as in \tmtextbf{Theorem \ref{thm:2}}, and we start by bounding $\mathbb{E}_{\xi'} \mathbb{E}_k \Big[ \frac{(L_f^k + L_{\omega})^2}{\gamma_k(\gamma_k - \kappa)} \Big]$. 

For brevity we omit $\Ebb_k[\cdot]$  and notice that, for $\gamma_k \geq 2 \kappa$, that,
\begin{align}
\frac{(L_f^k + L_{\omega})^2}{\gamma_k(\gamma_k - \kappa)}\leq \frac{2(L_f^k + L_{\omega})^2}{\gamma_k^2} 
  = \frac{(L_f^k)^2}{\gamma_k^2}
  + \frac{2 L_f^k L_{\omega}}{\gamma_k^2}  +
  \frac{L_{\omega}^2}{\gamma_k^2} \nonumber
\end{align}
so that we can bound the three terms respectively using
\[ \mathbb{E}_{\xi'} \Big[ \tfrac{(L_f^k)^2}{\gamma_k^2} \Big] \leq
   \mathbb{E}_{\xi'} \Big[ \tfrac{(L_f^k)^2}{\max \{ (L_f')^2, \alpha^2 \} k^{2
   \zeta}} \Big] \leq \Big( \frac{\alpha + \sigma}{\alpha} \Big)^2
   \frac{1}{k^{2 \zeta}} \]
where we invoked {\tmstrong{Lemma \ref{lem:3.3}}} and take $X = L_f^k, Y =
L_f'$ and we recall that by reference Lipschitz property \ref{D2}: $\mathbb{E} [| L_f^k
- L_f' |^2] \leq \sigma^2$. Similarly, we can deduce that
\[ \mathbb{E}_{\xi'} \Big[ \tfrac{2 L_f^k L_{\omega}}{\gamma_k^2} \Big]
   \leq \mathbb{E}_{\xi'} \Big[ \tfrac{2 L_f^k L_{\omega}}{\max \{
   \alpha L_f', \alpha^2 \} k^{2 \zeta}} \Big]
   \leq \Big( \frac{\alpha + \sigma}{\alpha^2} \Big) \frac{2
   L_{\omega}}{k^{2 \zeta}} \]
and $\frac{L_{\omega}^2}{\gamma_k^2} \leq \frac{L_{\omega}^2}{\alpha^2 k^{2
   \zeta}}$ since $\gamma_k \geq \alpha k^\zeta$. Putting the things together, we have
\[ \mathbb{E}_{\xi'} \Big[ \tfrac{(L_f^k + L_{\omega})^2}{\gamma_k^2} \Big]
   \leq \frac{(\alpha + \sigma)^2 + 2 L_{\omega} (\alpha + \sigma) +
   L_{\omega}^2}{\alpha^2 k^{2 \zeta}} = \frac{(\alpha + \sigma +
   L_{\omega})^2}{\alpha^2 k^{2 \zeta}} . \]
and
\[ \mathbb{E}_k [\psi_{1 / \rho} (x^{k + 1})] \leq \psi_{1 / \rho} (x^k) -
   \Ebb_{\xi'} \biggl[\frac{\rho (\rho - \tau - \kappa)}{2 (\gamma_k - \kappa)}  \biggl] \| \hat{x}^k - x^k \|^2 + \frac{\rho
   (\alpha + \sigma + L_{\omega})^2}{\alpha^2 k^{2 \zeta}}, \]
or
\[ \mathbb{E}_k [\psi_{1 / \rho} (x^{k + 1}) + \Lambda] \leq [\psi_{1 / \rho}
   (x^k) + \Lambda] - \Ebb_{\xi'} \biggl[\frac{\rho (\rho - \tau - \kappa)}{2 (\gamma_k - \kappa)}  \biggl] \| \hat{x}^k - x^k
   \|^2 + \frac{\rho (\alpha + \sigma + L_{\omega})^2}{\alpha^2 k^{2
   \zeta}} . \]
Invoking \textbf{Lemma \ref{lem:robbins}}, plugging in the relation $A_k = 0, V_k = \psi_{1 / \rho} (x^k) + \Lambda \geq
0$, $\mathsf{B}_k = \frac{\rho (\alpha + \sigma + L_{\omega})^2}{\alpha^2 k^{2
\zeta}}$ and $C_k = \Ebb_{\xi'} \big[ \frac{\rho (\rho - \tau - \kappa)}{\gamma_k - \kappa}\big] \| \hat{x}^k - x^k
\|^2$. Note that since $\Ebb_{\xi'}\big[ \frac{\rho (\rho - \tau - \kappa)}{\gamma_k - \kappa}\big]$ is determined by $x^k$, we can view $C_k = g(\|x^k\|) \|\hat{x}^k - x^k\|^2$  for some function $g$ and  the rest of reasoning is the same as in \tmtextbf{Lemma \ref{thm:2}}.

\subsection{Proof of Lemma \ref{lem:3.5}}

Similar to \tmtextbf{Lemma \ref{lem:2.2}} we first bound the error of
potential reduction $\mathbb{E}_k \Big[ \tfrac{\rho (L_f^k + L_{\omega})^2}{2
\gamma_k(\gamma_k - \kappa)} \Big]$, and according to the proof of \textbf{Lemma \ref{thm:3}},
\[ \mathbb{E}_{\xi'} \mathbb{E}_k \Bigg[ \frac{\rho (L_f^k + L_{\omega})^2}{2 \gamma_k(\gamma_k - \kappa)}
   \Bigg] \leq \frac{\rho (\alpha + \sigma + L_{\omega})^2}{\alpha^2 K} .
\]
Telescoping the relation $\mathbb{E} [\psi_{1 / \rho} (x^{k + 1})] \leq
\psi_{1 / \rho} (x^k) + \frac{\rho (\alpha + \sigma + L_{\omega})^2}{2
\alpha^2 K}$ gives
\[ \mathbb{E} [\psi_{1 / \rho} (x^k) + \Lambda] \leq \psi_{1 / \rho} (x^1) +
   \Lambda + \frac{\rho (\alpha + \sigma + L_{\omega})^2}{\alpha^2}
   =: \Delta . \]
Next we show that $\mathbb{E}_k [\| x^{k + 1} - x^k \|]$ is bounded. By the
optimality condition we have
\[ f_{x^k} (x^{k + 1}, \xi^k) + \omega (x^{k + 1}) + \frac{\gamma_k}{2} \|
   x^{k + 1} - x^k \|^2 \leq{} f_{x^k} (x^k, \xi^k) + \omega (x^k) -
   \frac{\gamma_k - \kappa}{2} \| x^{k + 1} - x^k \|^2 \]
and
\begin{align}
  \frac{\gamma_k}{2} \| x^{k + 1} - x^k \|^2 \leq{} & f_{x^k} (x^k, \xi^k) - f_{x^k} (x^{k
  + 1}, \xi^k) + \omega (x^k) - \omega (x^{k + 1}) \nonumber\\
  \leq{} & (L_f^k + L_{\omega}) \| x^{k + 1} - x^k \| . \nonumber
\end{align}

Dividing both sides by $\| x^{k + 1} - x^k \|$, we get
$\| x^{k + 1} - x^k \| \leq \frac{2(L_f^k + L_{\omega})}{\gamma_k} \leq
   \frac{2L_f^k}{\gamma_k} + \frac{2L_{\omega}}{\alpha \sqrt{K}}.$
Taking expectation on both sides, we have $\mathbb{E}_{\xi'} \mathbb{E}_k [L_f^k / \gamma_k] \leq \frac{\alpha + \sigma}{\alpha \sqrt{K}}$,
where we invoke \textbf{Lemma \ref{lem:3.3}} with $X = L_f^k, Y = L_f'$ again. Now 
\[ \mathbb{E}_k [\| x^{k + 1} - x^k \|] \leq \frac{2(\alpha + \sigma +
   L_{\omega})}{\alpha \sqrt{K}} . \]
The rest of the reasoning is consistent with \tmtextbf{Lemma \ref{lem:2.2}} up
to difference in constants. And we still present them for completeness. Applying Markov's inequality, we know that
\[ \mathbb{P}_{\xi^k, \xi' \sim \Xi} \{ \| x^{k + 1} - x^k \| \leq \tfrac{4
   (\alpha + \sigma + L_{\omega})}{\alpha \sqrt{K}} \big| \xi_1, \ldots, \xi_{k -
   1} \} \geq \frac{1}{2} . \]

By \ref{A7}, $\| x^k \| \geq \mathsf{B}_v$ implies $\psi_{1 / \rho} (x^k) \geq v$.
Taking $v = a \Delta$, we have $\| x^k \| \geq \mathsf{B}_{a \Delta} \Rightarrow
\psi_{1 / \rho} (x^k) \geq a v$. Without loss of generality, let $\mathsf{Z} = \frac{4
(\alpha + \sigma + L_{\omega})}{\alpha \sqrt{K}} > 0$, and we condition on the
event $\| x^k \| \geq \mathsf{B}_{a \Delta} + \mathsf{Z}$ to deduce that
\begin{align}
  \Delta \geq{} & \mathbb{E} [\psi_{1 / \rho} (x^{k + 1}) + \Lambda] \nonumber\\
  ={} & \mathbb{E} [\psi_{1 / \rho} (x^{k + 1}) + \Lambda | \| x^k \| \geq \mathsf{B}_{a
  \Delta} + \mathsf{Z}] \cdot \mathbb{P} \{ \| x^k \| \geq \mathsf{B}_{a \Delta} +\mathsf{Z} \} \nonumber\\
  & +\mathbb{E} [\psi_{1 / \rho} (x^{k + 1}) + \Lambda | \| x^k \| \leq{} \mathsf{B}_{a
  \Delta} + \mathsf{Z}] \cdot \mathbb{P} \{ \| x^k \| \leq \mathsf{B}_{a \Delta} +\mathsf{Z} \}
  \nonumber\\
  \geq{} & \mathbb{E} [\psi_{1 / \rho} (x^{k + 1}) + \Lambda | \| x^k \| \geq
  \mathsf{B}_{a \Delta} + \mathsf{Z}] \cdot  \mathbb{P} \{ \| x^k \| \geq \mathsf{B}_{a \Delta} +\mathsf{Z} \}, \label{eqn:c-5}
\end{align}

where \eqref{eqn:c-5} uses $\psi_{1 / \rho} (x) + \Lambda \geq 0$ for all $x$.
Next we consider the expectation $\mathbb{E} [\psi_{1 / \rho} (x^{k + 1}) +
\Lambda | \| x^k \| \geq \mathsf{B}_{a \Delta} + \mathsf{Z}]$ and we successively deduce that
\begin{align}
  & \mathbb{E} [\psi_{1 / \rho} (x^{k + 1}) + \Lambda | \| x^k \| \geq \mathsf{B}_{a
  \Delta} + \mathsf{Z}] \nonumber\\
  ={} & \mathbb{E} [\psi_{1 / \rho} (x^{k + 1}) + \Lambda | \| x^k \| \geq
  \mathsf{B}_{a \Delta} + \mathsf{Z}, \| x^{k + 1} - x^k \| \leq \mathsf{Z}] \cdot \mathbb{P} \{ \| x^{k + 1} -
  x^k \| \leq \mathsf{Z} \} \nonumber\\
  & +\mathbb{E} [\psi_{1 / \rho} (x^{k + 1}) + \Lambda | \| x^k \| \geq \mathsf{B}_{a
  \Delta} + \mathsf{Z}, \| x^{k + 1} - x^k \| \geq \mathsf{Z}] \cdot \mathbb{P} \{ \| x^{k + 1} - x^k
  \| \geq  \mathsf{Z} \} \nonumber\\
  \geq{} & \Big( \frac{a \Delta + \Lambda}{2} \Big) \cdot \mathbb{P} \{ \| x^{k +
  1} - x^k \| \leq \mathsf{Z} \}, \label{eqn:c-6}
\end{align}

where \eqref{eqn:c-6} is by $\mathbb{P} \{ \| x^{k + 1} -
x^k \| \leq \mathsf{Z} \} \geq 0.5$ and that, conditioned on $\| x^{k + 1} - x^k \|
\leq  \mathsf{Z}$,
\[ \| x^{k + 1} \| = \| x^{k + 1} - x^k + x^k \| \geq \| x^k \| - \| x^{k + 1}
   - x^k \| \geq \| x^k \| - \mathsf{Z} \geq \mathsf{B}_{a \Delta} . \]
Chaining the above inequalities, we arrive at
\[ \Delta \geq \Big( \frac{a \Delta + \Lambda}{2} \Big) \cdot \mathbb{P} \{ \|
   x^k \| \geq \mathsf{B}_{a \Delta} + \mathsf{Z} \} \]

Dividing both sides by $\frac{a \Delta + \Lambda}{2}$ gives the desired tail
bound.

\subsection{Proof of Theorem \ref{thm:3.2}}

The proof again follows the clue of \tmtextbf{Theorem \ref{thm:2.2}}. Recall that $\mathsf{Z}
= \frac{4 (\alpha + \sigma + L_{\omega})}{\alpha \sqrt{K}}$ and given $\delta
\in (0, 1/4)$,
\[ \mathbb{P} \{ \| x^k \| \geq \mathsf{B}_{\delta^{- 1} \Delta} + \mathsf{Z} \} \leq \frac{2
   \Delta}{\delta^{- 1} \Delta + \Lambda} \leq 2 \delta . \]
Denoting $I_k =\mathbb{I} \{ \| x^k \| \leq \mathsf{B}_{\delta^{- 1} \Delta} + \mathsf{Z} \}$,
we have $\sum_{k = 1}^K \mathbb{E} [1 - I_k] \leq 2 \delta K$ and by Markov's
inequality, given $p \in (2\delta, 1)$, we have
\[ \mathbb{P} \Big\{ \sum_{k = 1}^K (1 - I_k) \geq \frac{2\delta K}{p}
   \Big\} \leq \frac{2 \delta K}{2 p^{- 1} \delta K} = p \]
and $\mathbb{P} \{ \sum_{k = 1}^K I_k \geq K (1 - 2 p^{- 1} \delta)
\} \geq 1 - p$.  Now we telescope over \tmtextbf{Lemma \ref{lem:3.4}} and get

\begin{align}
  \sum_{k = 1}^K \mathbb{E} \left[ \frac{\rho (\rho - \tau - \kappa)}{2
  \gamma_k} \| \hat{x}^k - x^k \|^2 \right] \leq{} & \psi_{1 / \rho} (x^1)
  -\mathbb{E} [\psi_{1 / \rho} (x^{K + 1})] + \frac{\rho (\alpha + \sigma +
  L_{\omega})^2}{2 \alpha^2} \nonumber\\
  \leq{} & \psi_{1 / \rho} (x^1) - \inf_x \psi (x) + \frac{\rho}{2 \alpha^2} (\alpha + \sigma
  + L_{\omega})^2 \nonumber\\
  ={} & D + \frac{\rho}{2 \alpha^2} (\alpha + \sigma + L_{\omega})^2. \nonumber
\end{align} 

Next define
\[ \mathsf{G}_\delta \assign \max_x ~\sup_{\xi \sim \Xi}~ \textup{\textsf{Lip}} (x, \xi) \text{\quad subject to\quad} \| x \| \leq
   \mathsf{B}_{\delta^{- 1} \Delta} + \mathsf{Z}, \]
and we have, conditioned on the event $\sum_{k = 1}^K I_k \geq K (1 - 2 p^{-
1} \delta)$,
\begin{align}
  & \sum_{k = 1}^K \mathbb{E} \left[ \frac{\rho (\rho - \tau - \kappa)}{2
  \gamma_k} \| \hat{x}^k - x^k \|^2 \right] \nonumber\\
  ={} & {\sum_{k \in \{ j : I_j = 0 \}}}  \mathbb{E} \left[ \frac{\rho (\rho -
  \tau - \kappa)}{2 \gamma_k} \| \hat{x}^k - x^k \|^2 \right] + {\sum_{k \in
  \{ j : I_j = 1 \}}}  \mathbb{E} \left[ \frac{\rho (\rho - \tau - \kappa)}{2
  \gamma_k} \| \hat{x}^k - x^k \|^2 \right] \nonumber\\
  \geq{} & {\sum_{k \in \{ j : I_j = 1 \}}}  \mathbb{E} \left[ \frac{\rho (\rho
  - \tau - \kappa)}{2 \gamma_k} \| \hat{x}^k - x^k \|^2 \right] \nonumber\\
  \geq{} & {\sum_{k \in \{ j : I_j = 1 \}}}  \mathbb{E} \left[ \frac{\rho (\rho
  - \tau - \kappa)}{2 ( \rho + \kappa + \tau + (\alpha + \mathsf{G}_{\delta}) \sqrt{K}
  )} \| \hat{x}^k - x^k \|^2 \right] \nonumber\\
  ={} & \frac{\rho (\rho - \tau - \kappa)}{2 ( \rho + \kappa + \tau + (\alpha + \mathsf{G}_{\delta}) \sqrt{K} )} {\sum_{k \in \{ j : I_j = 1 \}}} 
  \mathbb{E} [\| \hat{x}^k - x^k \|^2] \nonumber\\
  \geq{} & \frac{\rho (\rho - \tau - \kappa) (1 - 2 p^{- 1} \delta) K}{2 (
  \rho + \kappa + \tau + (\alpha + \mathsf{G}_{\delta}) \sqrt{K} )} \min_{k \in \{ j
  : I_j = 1 \}} \mathbb{E} [\| \hat{x}^k - x^k \|^2] \nonumber
\end{align}
Re-arranging the terms, we have, at least with probability $1 - p$, that
\begin{align}
  \min_{1 \leq k \leq K} \mathbb{E} [\| \nabla \psi_{1 / \rho} (x^k) \|^2]
  \leq{} & \min_{k \in \{ k : I_k = 1 \}} \mathbb{E} [\| \nabla \psi_{1 / \rho}
  (x^k) \|^2] \nonumber\\
  \leq{} & \frac{p}{p - 2 \delta} \cdot \frac{2 \rho}{\rho - \tau - \kappa}
  \left[ D + \frac{\rho}{\alpha^2} (\alpha + \sigma + L_{\omega})^2 \right]
  \left( \frac{\rho + \lambda}{K} + \frac{\alpha + \mathsf{G}_{\delta}}{\sqrt{K}} \right) \nonumber
\end{align}
and this completes the proof after re-arrangement.

\section{Stochastic convex optimization beyond Lipschitz continuity} \label{app:cvx}

In this section, we consider applying the above mentioned ideas to convex
optimization. When dealing with convex optimization problems,
instead of relying on Moreau envelope smoothing, we have a better potential
function $\| x - x^{\star} \|$ directly relevant to distance to optimal set
$\mathcal{X}^{\star}$. This turns out greatly simplifies our assumptions.
\begin{enumerate}[leftmargin=35pt,itemsep=2pt,label=\textbf{E\arabic*:},ref=\rm{\textbf{E\arabic*}}]
  \item $f(x, \xi)$ is convex for all $\xi\sim\Xi$. $\lambda = \kappa = 0$ and  $\tau = 0$.\label{E1}
\end{enumerate}
It is rather straight-forward to extend our results to convex optimization. And we remark that our analysis is different from \cite{mai2021stability}, where the authors focus on subgradient method and assume quadratic growth condition. 

\subsection{Convex optimization under standard Lipschitzness} \label{sec:cvx-known-lip}

\begin{lem} \label{cvx-lem:1}
Suppose that \ref{A1} to \ref{A5}, \ref{E1} as well as \ref{B1} holds, then given $\gamma_k > 0$
\begin{equation}
	\mathbb{E}_k [\| x^{k + 1} - x^{\star} \|^2] \leq \| x^k - x^{\star} \|^2 -
     \frac{2}{\gamma_k} [\psi (x^k) - \psi (x^{\star})] + \frac{(L_f +
     L_{\omega})^2}{\gamma_k^2}, \label{eqn:d-1}
\end{equation}
  where $x^{\star} \in \mathcal{X}^{\star}$ is any optimal solution.
\end{lem}

\begin{thm}\label{cvx-thm:1}
  Under the same assumptions as \textbf{Lemma \ref{cvx-lem:1}}, if we take $\gamma_k \equiv \gamma = \alpha \sqrt{K}$, then
  \[ \min_{1 \leq k \leq K} \mathbb{E} [\psi (x^{k}) - \psi (x^{\star})] \leq \frac{1}{2\sqrt{K}} \Big [\| x^1 - x^{\star} \|^2 \alpha + \frac{(L_f +  L_{\omega})^2}{\alpha} \Big], \]
  where $x^\star \in \Xcal^\star$ is an optimal solution.
\end{thm}

\begin{rem}
  We observe the same trade-off as in weakly convex optimization, where we
  have, given telescopic sum of \eqref{eqn:d-1}, that
  \[ \frac{1}{K} \sum_{k = 1}^K \mathcal{O} (\gamma_k^{- 1}) \Ebb[\psi (x^k) -
     \psi (x^{\star})] \leq \mathcal{O} \Big( \frac{1}{K} \Big) +
     \frac{1}{K} \sum_{k = 1}^K \mathcal{O}(L_f^2\gamma_k^{-2}) . \]
     Compared with weakly convex case
     \[ \frac{1}{K}\sum_{k=1}^K \mathcal{O}(\gamma_k^{-1}) \Ebb [\|\nabla\psi_{1/\rho}(x^k) \|^2] \leq \mathcal{O}\Big(\frac{1}{K}\Big) + \frac{1}{K}\sum_{k=1}^K \mathcal{O}(L_f^2\gamma_k^{-2}), \]
  this resemblance implies our previous analysis for weakly convex
  optimization is immediately applicable.
\end{rem}

\subsection{Convex optimization under generalized Lipschitzness} \label{sec:cvx-known-nolip}

\begin{lem} \label{cvx:lem:2}
Suppose \ref{A1} to \ref{A5}, \ref{E1} as well as \ref{C1} holds, then given $\gamma_k > 0$,
  \[ \mathbb{E}_k [\| x^{k + 1} - x^{\star} \|^2] \leq \| x^k - x^{\star} \|^2 -
     \frac{2}{\gamma_k} [\psi (x^k) - \psi (x^{\star})] + \frac{(\mathcal{G}
     (\| x^k \|) L_f + L_{\omega})^2}{\gamma_k^2}, \]
where $x^\star \in \Xcal^\star$ is an optimal solution.
\end{lem}

\begin{thm} \label{cvx:thm:2}
  With the same conditions as \textbf{Lemma \ref{cvx:lem:2}}, if
  $\gamma_k =  (\mathcal{G} (\| x^k \|) + 1) k^{\zeta}, \zeta \in \left(
  \frac{1}{2}, 1 \right)$, then as $k \rightarrow \infty$, $\{ \| x^k \| \}$
  is bounded with probability 1 and $\{ \inf_{j \leq k} f (x^j) - f (x^{\star}) \}$ converges to 0 almost surely.
\end{thm}

\begin{lem} \label{cvx:lem:3}
Under the same conditions as \textbf{Lemma \ref{cvx:lem:2}}, if we take $\gamma_k = \alpha (\mathcal{G} (\| x^k \|) + 1) \sqrt{K}$, then the tail bound
  \[ \mathbb{P} \left\{ \| x^k - x^{\star} \| \geq \frac{2 (L_f +
     L_{\omega})}{\alpha \sqrt{K}} + a \right\} \leq \frac{2 \Delta}{a}, \]
holds for all $2\leq k\leq K$, where $\Delta = \| x^1 - x^{\star} \| + \frac{L_f + L_{\omega}}{\alpha}$.
\end{lem}

\begin{thm} \label{cvx:thm:3}
Under the same conditions as \textbf{Lemma \ref{cvx:lem:2}}, given $\delta\in(0,1/4), p \in (2\delta, 1)$, $(1-2p^{-1}\delta) K$ iterations will lie in the ball centered around $x^\star$ with radius $\mathsf{R}(\delta) = \delta^{-1}  \Delta + \frac{2(L_f + L_{\omega})}{\alpha \sqrt{K}}$ and 
\begin{equation}
	\min_{1 \leq k \leq K} \mathbb{E} [\psi (x^k) - \psi (x^{\star})]
	\leq \frac{p}{p - 2 \delta} \cdot \frac{\mathsf{G}_{\delta} + 1}{2
   \sqrt{K}} \left[ \| x^1 - x^{\star} \|^2 \alpha + \frac{(L_f +
   L_{\omega})^2}{\alpha}, \right]
\end{equation}
where $\mathsf{G}_{\delta} \assign \max_z \mathcal{G} (z), \| z - x^{\star} \| \leq
   \delta^{- 1} \Delta + \frac{2 (L_f +
     L_{\omega})}{\alpha \sqrt{K}}$.
\end{thm}

\begin{rem}
We note that $x^\star$ is actually arbitrary. Therefore we can take it to be a minimum norm optimal solution to get a tighter bound.
\end{rem}

\subsection{Convex optimization beyond Lipschitzness} \label{sec:cvx-unknown-nolip}

\begin{lem} \label{cvx:lem:4}
Suppose that \ref{A1} to \ref{A5}, \ref{E1} as well as \ref{D1}, \ref{D2} hold, then given $\gamma > 0$, 

\begin{equation}
	\mathbb{E}_k [\| x^{k + 1} - x^{\star} \|^2] \leq \| x^k - x^{\star} \|^2 -
   \frac{2}{\gamma_k} [\psi (x^k) - \psi (x^{\star})] +\mathbb{E}_k \bigg[
   \frac{(\mathsf{Lip}(x^k, \xi^k)) + L_{\omega})^2}{\gamma_k^2} \bigg],
\end{equation}
where $\gamma_k$ is chosen to be independent of $\xi^k$ and is considered deterministic here.
\end{lem}

\begin{thm} \label{cvx:thm:4}
  With the same conditions as \textbf{Lemma \ref{cvx:lem:4}}, if
  $\gamma_k =  \max \{ \mathsf{Lip}(x^k, \xi^k)), \alpha  \} k^{\zeta}, \zeta \in (
  \frac{1}{2}, 1 )$, then as $k \rightarrow \infty$, $\{ \| x^k \| \}$
  is bounded with probability 1 and $\{ \inf_{j \leq k} f (x^j) - f (x^{\star}) \}$ converges to 0 almost surely. 
\end{thm}

\begin{lem} \label{cvx:lem:5}
Under the same conditions as \textbf{Lemma \ref{cvx:lem:4}}, if we take $\gamma_k = \max \{ \mathsf{Lip}(f(x^k, \xi^k), \alpha  \} \sqrt{K}$, then the tail bound
  \[ \mathbb{P} \left\{ \| x^k - x^{\star} \| \geq \frac{2 (\alpha + \sigma +
     L_{\omega})}{\alpha \sqrt{K}} + a \right\} \leq \frac{2 \Delta}{a}, \]
holds for all $2\leq k\leq K$, where $\Delta = \| x^1 - x^{\star} \| + \frac{\alpha + \sigma + L_{\omega}}{\alpha}$.
\end{lem}

\begin{rem}
Now that in convex optimization our potential function $\|x - x^\star\|^2$ itself already defines a bounded set. The proof of \textbf{\ref{cvx:lem:5}} can also be done based on a conditional probability argument.
\end{rem}

\begin{thm} \label{cvx:thm:5}
Under the same conditions as \textbf{Lemma \ref{cvx:lem:4}}, given $\delta\in(0,1/4), p \in (2\delta, 1)$, $(1-2p^{-1}\delta) K$ iterations will lie in the ball centered around $x^\star$ with radius $\mathsf{R}(\delta) = \delta^{-1}  \Delta + \frac{2(\alpha + \sigma + L_{\omega})}{\alpha \sqrt{K}}$ and 
\begin{equation}
	\min_{1 \leq k \leq K} \mathbb{E} [\psi (x^k) - \psi (x^{\star})]
	\leq \frac{p}{p - 2 \delta} \cdot \frac{\mathsf{G}_{\delta} + \alpha }{2
   \sqrt{K}} \left[ \| x^1 - x^{\star} \|^2 \alpha + \frac{(\alpha + \sigma +
   L_{\omega})^2}{\alpha}, \right]
\end{equation}
where $\mathsf{G}_{\delta} \assign \max_x \sup_{\xi \sim \Xi} \mathsf{Lip}(x, \xi) , \| x - x^{\star} \| \leq
   \delta^{- 1} \Delta + \frac{2 (\alpha + \sigma +
     L_{\omega})}{\alpha \sqrt{K}}$.
\end{thm}

\subsection{Proof of results in Subsection \ref{sec:cvx-known-lip}}

\subsubsection{Proof of Lemma \ref{cvx-lem:1}}

Let $x^{\star} \in \mathcal{X}^{\star}$ be an optimal solution to the problem.
Then by three-point lemma, we have
\begin{align}
  & f_{x^k} (x^{k + 1}, \xi^k) + \omega (x^{k + 1}) + \frac{\gamma_k}{2} \|
  x^{k + 1} - x^k \|^2 \nonumber\\
  \leq{} & f_{x^k} (x^{\star}, \xi^k) + \omega (x^{\star}) + \frac{\gamma_k}{2} \|
  x^k - x^{\star} \|^2 - \frac{\gamma_k}{2} \| x^{k + 1} - x^{\star} \|^2 .
  \nonumber
\end{align}
Re-arranging the terms, we deduce that
\begin{align}
  & \frac{\gamma_k}{2} \| x^{k + 1} - x^{\star} \|^2 \nonumber\\
  \leq{} & \frac{\gamma_k}{2} \| x^k - x^{\star} \|^2 - \frac{\gamma_k}{2} \|
  x^{k + 1} - x^k \|^2 + f_{x^k} (x^{\star}, \xi^k) + \omega (x^{\star}) -
  f_{x^k} (x^{k + 1}, \xi^k) - \omega (x^{k + 1}) \nonumber\\
  \leq{} & \frac{\gamma_k}{2} \| x^k - x^{\star} \|^2 - \frac{\gamma_k}{2} \|
  x^{k + 1} - x^k \|^2 + (L_f (\xi^k) + L_{\omega}) \| x^{k + 1} - x^k \|
  \label{eqn:d-2} \\
  & + f_{x^k} (x^{\star}, \xi^k) + \omega (x^{\star}) - f (x^k) - \omega (x^k),
  \nonumber
\end{align}
where \eqref{eqn:d-2} applies \ref{B1} to get $f_{x^k} (x^{k + 1}, \xi^k)-f_{x^k} (x^{k}, \xi^k) \leq L_f(\xi^k) \|x^{k+1} - x^k\|$. Dividing both sides by $\frac{\gamma_k}{2}$,
\begin{align}
  \| x^{k + 1} - x^{\star} \|^2 \leq{} & \| x^k - x^{\star} \|^2 - \| x^{k + 1} -
  x^k \|^2 + \frac{2 (L_f (\xi^k) + L_{\omega})}{\gamma_k} \| x^{k + 1} - x^k
  \| \nonumber\\
  & + \frac{2}{\gamma_k} [f (x^{\star}, \xi^k) + \omega (x^{\star}) - f (x^k) -
  \omega (x^k)] + \frac{2}{\gamma_k} [f_{x^k}(x^\star, \xi) - f(x^\star, \xi)].   \nonumber
\end{align}
Conditioned on $x^1, \ldots, x^k$ and taking expectation with respect to
$\xi^k$, we successively deduce that
\begin{align}
  & \mathbb{E}_k [\| x^{k + 1} - x^{\star} \|^2] \nonumber\\
  \leq{} & \| x^k - x^{\star} \|^2 -\mathbb{E}_k [\| x^{k + 1} - x^k \|^2]
  +\mathbb{E}_k [2 \gamma_k^{- 1} (L_f (\xi^k) + L_{\omega}) \| x^{k + 1} -
  x^k \|] \nonumber \\
  & + \frac{2}{\gamma_k} [f (x^{\star}) + \omega (x^{\star}) - f (x^k) - \omega
  (x^k)] \label{eqn:d-3}\\
  \leq{} & \| x^k - x^{\star} \|^2 + \frac{(L_f + L_{\omega})^2}{\gamma_k^2} +
  \frac{2}{\gamma_k} [f (x^{\star}) + \omega (x^{\star}) - f (x^k) - \omega
  (x^k)] \label{eqn:d-4} \\
  ={} & \| x^k - x^{\star} \|^2 + \frac{(L_f + L_{\omega})^2}{\gamma_k^2} +
  \frac{2}{\gamma_k} [\psi (x^{\star}) - \psi (x^k)], \nonumber
\end{align}
where \eqref{eqn:d-3} applies \ref{E1} to get $\Ebb_{\xi^k}[f_{x^k}(x^\star, \xi) - f(x^\star, \xi)] \leq 0$; \eqref{eqn:d-4} uses $- \frac{a}{2} x^2 +b x \leq \frac{b^2}{2a}$ and that $\Ebb[L_f(\xi)^2] \leq L_f^2$. Re-arranging the terms, we arrive at
\[ \mathbb{E}_k [\| x^{k + 1} - x^{\star} \|^2] \leq \| x^k - x^{\star} \|^2 -
   \frac{2}{\gamma_k} [\psi (x^k) - \psi (x^{\star})] + \frac{(L_f +
   L_{\omega})^2}{\gamma_k^2} \]
and this completes the proof.

\subsubsection{Proof of Theorem \ref{cvx-thm:1}}
Taking $\gamma_k \equiv \gamma = \alpha \sqrt{K}$ and telescoping from $1, \ldots, K$, we have
\begin{align}
  \min_{1 \leq k \leq K} \mathbb{E} [\psi (x^{k}) - \psi (x^{\star})] \leq{} & \frac{1}{K}
  \sum_{k = 1}^K \mathbb{E} [\psi (x^k) - \psi (x^{\star})] \nonumber\\
  \leq{} & \frac{1}{2\sqrt{K}} \Big [\| x^1 - x^{\star} \|^2 \alpha + \frac{(L_f +  L_{\omega})^2}{\alpha} \Big] \nonumber
\end{align}
and this completes the proof.

\subsection{Proof of results in Subsection \ref{sec:cvx-known-nolip}}

\subsubsection{Proof of Lemma \ref{cvx:lem:2}}
Define $\mathsf{G}_k \assign \mathcal{G} (\| x^k \|)$. Let $x^{\star} \in
\mathcal{X}^{\star}$ be an optimal solution to the problem. Similarly, we have
\begin{align}
  & f_{x^k} (x^{k + 1}, \xi^k) + \omega (x^{k + 1}) + \frac{\gamma_k}{2} \|
  x^{k + 1} - x^k \|^2 \nonumber\\
  \leq{} & f_{x^k} (x^{\star}, \xi^k) + \omega (x^{\star}) + \frac{\gamma_k}{2} \|
  x^k - x^{\star} \|^2 - \frac{\gamma_k}{2} \| x^{k + 1} - x^{\star} \|^2 .
  \nonumber
\end{align}
Re-arranging the terms, for $\xi^k \sim \Xi$, we deduce that
\begin{align}
  & \frac{\gamma_k}{2} \| x^{k + 1} - x^{\star} \|^2 \nonumber\\
  \leq{} & \frac{\gamma_k}{2} \| x^k - x^{\star} \|^2 - \frac{\gamma_k}{2} \|
  x^{k + 1} - x^k \|^2 + f_{x^k} (x^{\star}, \xi^k) + \omega (x^{\star}) -
  f_{x^k} (x^{k + 1}, \xi^k) - \omega (x^{k + 1}) \nonumber\\
  \leq{} & \frac{\gamma_k}{2} \| x^k - x^{\star} \|^2 - \frac{\gamma_k}{2} \|
  x^{k + 1} - x^k \|^2 + (\mathsf{G}_k L_f (\xi^k) + L_{\omega}) \| x^{k + 1} - x^k \|
  \label{eqn:d-5} \\
  & + f (x^{\star}, \xi^k) + \omega (x^{\star}) - f (x^k) - \omega (x^k),
  \nonumber
\end{align}
where \eqref{eqn:d-5} applies convexity. Dividing both sides by $\frac{\gamma_k}{2}$, we have
\begin{align}
  \| x^{k + 1} - x^{\star} \|^2 \leq{} & \| x^k - x^{\star} \|^2 - \| x^{k + 1} -
  x^k \|^2 + \frac{2 (\mathsf{G}_k L_f (\xi^k) + L_{\omega})}{\gamma_k} \| x^{k + 1} -
  x^k \| \nonumber\\
  & + \frac{2}{\gamma_k} [f (x^{\star}, \xi^k) + \omega (x^{\star}) - f (x^k) -
  \omega (x^k)] \nonumber
\end{align}
Next, conditioned on $x^1, \ldots, x^k$ and taking expectation with respect to
$\xi^k$, we successively deduce that
\begin{align}
  & \mathbb{E}_k [\| x^{k + 1} - x^{\star} \|^2] \nonumber\\
  \leq{} & \| x^k - x^{\star} \|^2 -\mathbb{E}_k [\| x^{k + 1} - x^k \|^2]
  +\mathbb{E}_k [2 \gamma_k^{- 1} (\mathsf{G}_k L_f (\xi^k) + L_{\omega}) \| x^{k + 1}
  - x^k \|] \nonumber\\
  & + \frac{2}{\gamma_k} [f (x^{\star}) + \omega (x^{\star}) - f (x^k) - \omega
  (x^k)] \nonumber\\
  \leq{} & \| x^k - x^{\star} \|^2 + \frac{(\mathsf{G}_k L_f + L_{\omega})^2}{\gamma_k^2}
  + \frac{2}{\gamma_k} [f (x^{\star}) + \omega (x^{\star}) - f (x^k) - \omega
  (x^k)] \label{eqn:d-6} \\
  ={} & \| x^k - x^{\star} \|^2 + \frac{(\mathsf{G}_k L_f + L_{\omega})^2}{\gamma_k^2} +
  \frac{2}{\gamma_k} [\psi (x^{\star}) - \psi (x^k)], \nonumber
\end{align}
where \eqref{eqn:d-6} again uses $- \frac{a}{2} x^2 +b x \leq \frac{b^2}{2a}$ and the assumption $\mathbb{E}_{\xi} [L_f (\xi)^2] \leq L_f^2$. Re-arranging the terms, we get
\begin{equation}
	\mathbb{E}_k [\| x^{k + 1} - x^{\star} \|^2] \leq \| x^k - x^{\star} \|^2 -
   \frac{2}{\gamma_k} [\psi (x^k) - \psi (x^{\star})] + \frac{(\mathsf{G}_k L_f +
   L_{\omega})^2}{\gamma_k^2} \label{eqn:telesum}
\end{equation}
and this completes the proof.

\subsubsection{Proof of Theorem \ref{cvx:thm:2}}

Now that our recursive potential reduction is changed into
\[ \mathbb{E}_k [\| x^{k + 1} - x^{\star} \|^2] \leq \| x^k - x^{\star} \|^2 -
   \frac{2}{\gamma_k} [\psi (x^k) - \psi (x^{\star})] + \frac{(\mathsf{G}_k L_f +
   L_{\omega})^2}{\gamma_k^2} \]
and we bound
\[ \frac{\mathsf{G}_k L_f + L_{\omega}}{\gamma_k} = \frac{\mathsf{G}_k L_f + L_{\omega}}{(\mathsf{G}_k +
   1) k^{\zeta}} \leq \frac{L_f + L_{\omega}}{k^{\zeta}}, \]
giving
\[ \mathbb{E}_k [\| x^{k + 1} - x^{\star} \|^2] \leq \| x^k - x^{\star} \|^2 -
   \frac{2}{\gamma_k} [\psi (x^k) - \psi (x^{\star})] + \frac{(L_f +
   L_{\omega})^2}{k^{2 \zeta}} . \]
Invoking \textbf{Lemma \ref{lem:robbins}} with $A_k = 0, V_k = \| x^k - x^{\star} \|^2 \geq 0, B_k =
\frac{(L_f + L_{\omega})^2}{k^{2 \zeta}}$ and $C_k = \frac{2}{\gamma_k} [\psi
(x^k) - \psi (x^{\star})]$, we complete the proof with the same argument as in \textbf{Theorem \ref{thm:2}}.

\subsubsection{Proof of Lemma \ref{cvx:lem:3}} 
Our proof is a duplicate of \textbf{Lemma \ref{lem:2.2}} using a different potential function.
We start by bounding the error of potential reduction by $\frac{(\mathsf{G}_k L_f +
L_{\omega})^2}{\gamma_k^2} \leq \frac{(L_f + L_{\omega})^2}{\alpha^2 K} .$ Then telescoping of \eqref{eqn:telesum} gives us, for all $2 \leq k \leq K$, that
\[ \mathbb{E} [\| x^k - x^{\star} \|]^2 \leq \mathbb{E} [\| x^{k} -
   x^{\star} \|^2] \leq \| x^1 - x^{\star} \|^2 + \frac{(L_f +
   L_{\omega})^2}{\alpha^2} \leq \Big( \|x^1 - x^\star\| + \frac{L_f +
   L_{\omega}}{\alpha} \Big)^2 =: \Delta^2 \]
and $\mathbb{E} [\| x^k - x^{\star} \|] \leq \Delta$. Next consider, by optimality condition, that
\[ f_{x^k} (x^{k + 1}, \xi^k) + \omega (x^{k + 1}) + \frac{\gamma_k}{2} \|
   x^{k + 1} - x^k \|^2 \leq f_{x^k} (x^k, \xi^k) + \omega (x^k) -
   \frac{\gamma_k}{2} \| x^{k + 1} - x^k \|^2. \]
A re-arrangement gives
\[ \gamma_k \| x^{k + 1} - x^k \|^2 \leq (L_f (\xi^k) \mathsf{G}_k +
   L_{\omega}) \| x^{k + 1} - x^k \| . \]
Dividing both sides by $\| x^{k + 1} - x^k \|$,
\[ \| x^{k + 1} - x^k \| \leq \frac{L_f (\xi^k) \mathsf{G}_k +
   L_{\omega}}{\gamma_k} = \frac{L_f (\xi^k) \mathsf{G}_k + L_{\omega}}{\alpha
   (\mathsf{G}_k + 1) \sqrt{K}} . \]
Taking expectation, we get $\mathbb{E}_k [\| x^{k + 1} - x^k \|] \leq \frac{L_f +
L_{\omega}}{\alpha \sqrt{K}} .$ Then by Markov's inequality,
\[ \mathbb{P}_{\xi^k \sim \Xi} \left\{ \| x^{k + 1} - x^k \| \leq \frac{2
   (L_f + L_{\omega})}{\alpha \sqrt{K}} | \xi^1, \ldots, \xi^{k - 1} \right\}
   \geq \frac{1}{2} \]
and without loss of generality, we let $\mathsf{Z} = \frac{2 (L_f +
L_{\omega})}{\alpha \sqrt{K}} > 0$. Then we successively deduce that
\begin{align}
  \Delta \geq{} & \mathbb{E} [\| x^{k + 1} - x^{\star} \|] \nonumber\\
  ={} & \mathbb{E} [ \| x^{k + 1} - x^{\star} \| | \| x^k - x^{\star} \| \leq
  \mathsf{Z} + z ] \cdot \mathbb{P} \{ \| x^k - x^{\star} \| \leq
  \mathsf{Z} + z \} \nonumber\\
  & +\mathbb{E} [ \| x^{k + 1} - x^{\star} \| | \| x^k - x^{\star} \| \geq
  \mathsf{Z} + z ] \cdot \mathbb{P} \{ \| x^k - x^{\star} \| \geq
  \mathsf{Z} + z \} \nonumber\\
  \geq{} & \mathbb{E} [\| x^{k + 1} - x^{\star} \| | \| x^k - x^{\star} \| \geq
  \mathsf{Z} + z]. \nonumber
\end{align}
Next, consider the expectation $\mathbb{E} [\| x^{k + 1} - x^{\star} \| | \| x^k
- x^{\star} \| \geq \mathsf{Z} + z]$ and we successively
deduce that
\begin{align}
  & \mathbb{E} [\| x^{k + 1} - x^{\star} \| | \| x^k - x^{\star} \| \geq
  \mathsf{Z} + z] \nonumber\\
  ={} & \mathbb{E} [\| x^{k + 1} - x^{\star} \| | \| x^k - x^{\star} \| \geq
  \mathsf{Z} + z, \| x^{k + 1} - x^k \| \leq \mathsf{Z}] \cdot \mathbb{P} \{
  \| x^{k + 1} - x^k \| \leq \mathsf{Z} \} \nonumber\\
  & +\mathbb{E} [\| x^{k + 1} - x^{\star} \| | \| x^k - x^{\star} \| \geq
  \mathsf{Z} + z, \| x^{k + 1} - x^k \| \geq \mathsf{Z}] \cdot \mathbb{P} \{
  \| x^{k + 1} - x^k \| \geq \mathsf{Z} \} \nonumber\\
  \geq{} & \frac{z}{2} \cdot \mathbb{P} \{ \| x^{k + 1} - x^k \| \geq
  \mathsf{Z} \} \label{eqn:d-7}
\end{align}

where \eqref{eqn:d-7} is by $\mathbb{P} \{ \| x^{k + 1} - x^k \| \leq \mathsf{Z}
\} \geq 0.5$ and that, conditioned on $\| x^{k + 1} - x^k \| \leq \mathsf{Z}$,
\[ \| x^{k + 1} - x^{\star} \| = \| x^{k + 1} - x^k + x^k - x^{\star} \| \geq
   \| x^k - x^{\star} \| - \| x^{k + 1} - x^k \| \geq z. \]
Chaining the above inequalities, we arrive at
\[ \Delta \geq \frac{z}{2} \cdot \mathbb{P} \{ \| x^k - x^{\star} \| \geq
   \mathsf{Z} + z \} . \]
Dividing both sides by $\frac{z}{2}$ and taking $z=a$ gives the desired tail bound.

\subsubsection{Proof of Theorem \ref{cvx:thm:3}}
Given $\delta \in (0, 1 / 4)$, we have
\[ \mathbb{P} \{ \| x^k - x^{\star} \| \geq \mathsf{Z} + \delta^{- 1} \Delta \}
   \leq \frac{2 \Delta}{\delta^{- 1} \Delta} \leq 2 \delta . \]
Denote $I_k =\mathbb{I} \{ \| x^k - x^{\star} \| \leq \delta^{- 1} \Delta +
\mathsf{Z} \}$. We have $\sum_{k = 1}^K \mathbb{E} [1 - I_k] \leq 2 \delta K$
and by Markov's inequality, given $p \in (2 \delta, 1)$,
\[ \mathbb{P} \left\{ \sum_{k = 1}^K I_k \geq K (1 - 2 p^{- 1} \delta)
   \right\} \geq 1 - p. \]
Now we telescope over \tmtextbf{Lemma \ref{cvx:lem:3}} and deduce that
\begin{align}
  \sum_{k = 1}^K \mathbb{E} \left[ \frac{2}{\gamma_k} (\psi (x^k) - \psi
  (x^{\star})) \right] \leq{} & \| x^1 - x^{\star} \|^2 -\mathbb{E} [\| x^{K + 1}
  - x^{\star} \|^2] + \frac{(L_f + L_{\omega})^2}{\alpha^2} \nonumber\\
  \leq{} & \| x^1 - x^{\star} \|^2 + \frac{(L_f + L_{\omega})^2}{\alpha^2} .
  \nonumber
\end{align}
Next we condition on the event $\sum_{k = 1}^K I_k \geq K (1 - 2 p^{- 1}
\delta)$, define
$\mathsf{G}_{\delta} \assign \max_z \mathcal{G} (z), \| z - x^{\star} \| \leq
   \delta^{- 1} \Delta + \mathsf{Z}$,
and successively deduce that
\begin{align}
  & \sum_{k = 1}^K \mathbb{E} \left[ \frac{2}{\gamma_k} (\psi (x^k) - \psi
  (x^{\star})) \right] \nonumber\\
  ={} & \sum_{k \in \{ j : I_j = 0 \}} \mathbb{E} \left[ \frac{2}{\gamma_k}
  (\psi (x^k) - \psi (x^{\star})) \right] + \sum_{k \in \{ j : I_j = 1 \}}
  \mathbb{E} \left[ \frac{2}{\gamma_k} (\psi (x^k) - \psi (x^{\star})) \right]
  \nonumber\\
  \geq{} & \sum_{k \in \{ j : I_j = 1 \}} \mathbb{E} \left[ \frac{2}{\gamma_k}
  (\psi (x^k) - \psi (x^{\star})) \right] \nonumber\\
  \geq{} & \sum_{k \in \{ j : I_j = 1 \}} \mathbb{E} \left[ \frac{2}{\alpha
  (\mathsf{G}_{\delta} + 1) \sqrt{K}} (\psi (x^k) - \psi (x^{\star})) \right]
  \nonumber\\
  \geq{} & \frac{2 (1 - 2 p^{- 1} \delta) \sqrt{K}}{\alpha (\mathsf{G}_{\delta}
  + 1)} \min_{k \in \{ j : I_j = 1 \}} \mathbb{E} [(\psi (x^k) - \psi
  (x^{\star}))], \label{eqn:d-8}
\end{align}
where \eqref{eqn:d-8} follows from the event $\sum_{k = 1}^K I_k \geq K (1 - 2 p^{- 1}
\delta)$. Putting the results together, we have
\begin{align}
	 \min_{1 \leq k \leq K} \mathbb{E} [\psi (x^k) - \psi (x^{\star})]
	\leq{} & \min_{k \in \{ j : I_j = 0 \}} \mathbb{E} [\psi (x^k) - \psi (x^{\star})]\\
   \leq{} & \frac{p}{p - 2 \delta} \cdot \frac{\mathsf{G}_{\delta} + 1}{2
   \sqrt{K}} \left[ \| x^1 - x^{\star} \|^2 \alpha + \frac{(L_f +
   L_{\omega})^2}{\alpha} \right]
\end{align}
and this completes the proof.

\subsection{Proof of results in Subsection \ref{sec:cvx-unknown-nolip}}

In this section, we again define $L_f^k := \textup{\textsf{Lip}} (x^k, \xi^k), L_f' := \textup{\textsf{Lip}} (x^k, \xi')$  to simplify notation.

\subsubsection{Proof of Lemma \ref{cvx:lem:3}}

By the optimality condition we have
\begin{align}
  f_{x^k} (x^{k + 1}, \xi^k) + \omega (x^{k + 1}) + \frac{\gamma_k}{2} \|
  x^{k + 1} - x^k \|^2 
  \leq  f_{x^k} (x^{\star}, \xi^k) + \omega (x^{\star}) + \frac{\gamma_k}{2} \|
  x^k - x^{\star} \|^2 - \frac{\gamma_k}{2} \| x^{k + 1} - x^{\star} \|^2 .
  \nonumber
\end{align}
Re-arranging the terms, we get
\begin{align}
  & \frac{\gamma_k}{2} \| x^{k + 1} - x^{\star} \|^2 \nonumber\\
  \leq{} & \frac{\gamma_k}{2} \| x^k - x^{\star} \|^2 - \frac{\gamma_k}{2} \|
  x^{k + 1} - x^k \|^2 + f_{x^k} (x^{\star}, \xi^k) + \omega (x^{\star}) -
  f_{x^k} (x^{k + 1}, \xi^k) - \omega (x^{k + 1}) \nonumber\\
  \leq{} & \frac{\gamma_k}{2} \| x^k - x^{\star} \|^2 - \frac{\gamma_k}{2} \|
  x^{k + 1} - x^k \|^2 + (L_f^k + L_{\omega}) \| x^{k + 1} - x^k \|
  \nonumber\\
  & + f (x^{\star}, \xi^k) + \omega (x^{\star}) - f (x^k) - \omega (x^k) .
  \nonumber
\end{align}
Dividing both sides by $\frac{\gamma_k}{2}$,
\begin{align}
  \| x^{k + 1} - x^{\star} \|^2 \leq{} & \| x^k - x^{\star} \|^2 - \| x^{k + 1} -
  x^k \|^2 + \frac{2 (L_f^k + L_{\omega})}{\gamma_k} \| x^{k + 1} - x^k \|
  \nonumber\\
  & + \frac{2}{\gamma_k} [f (x^{\star}, \xi^k) + \omega (x^{\star}) - f (x^k) -
  \omega (x^k)] \nonumber
\end{align}
Conditioned on $x^1, \ldots, x^k$, \ taking expectation with respect to
$\xi^k$, we have
\begin{align}
  & \mathbb{E}_k [\| x^{k + 1} - x^{\star} \|^2] \nonumber\\
  \leq{} & \| x^k - x^{\star} \|^2 -\mathbb{E}_k [\| x^{k + 1} - x^k \|^2]
  +\mathbb{E}_k [2 \gamma_k^{- 1} (L_f^k + L_{\omega}) \| x^{k + 1} - x^k \|]
  \nonumber\\
  & + \frac{2}{\gamma_k} [f (x^{\star}) + \omega (x^{\star}) - f (x^k) - \omega
  (x^k)] \label{eqn:d-9}\\
  \leq{} & \| x^k - x^{\star} \|^2 + \frac{(L_f^k + L_{\omega})^2}{\gamma_k^2} +
  \frac{2}{\gamma_k} [f (x^{\star}) + \omega (x^{\star}) - f (x^k) - \omega
  (x^k)] \nonumber\\
  ={} & \| x^k - x^{\star} \|^2 + \frac{(L_f^k + L_{\omega})^2}{\gamma_k^2} +
  \frac{2}{\gamma_k} [\psi (x^{\star}) - \psi (x^k)], \nonumber
\end{align}

where \eqref{eqn:d-9} use the fact that $\gamma_k$ does not inherit randomness from
$\xi^k$. Re-arranging the terms, we get
\[ \mathbb{E}_k [\| x^{k + 1} - x^{\star} \|^2] \leq \| x^k - x^{\star} \|^2 -
   \frac{2}{\gamma_k} [\psi (x^k) - \psi (x^{\star})] +\mathbb{E}_k \Big[
   \tfrac{(L_f^k + L_{\omega})^2}{\gamma_k^2} \Big] \]
and this completes the proof.

\subsubsection{Proof of Theorem \ref{cvx:thm:4}}

We start by bounding $\mathbb{E}_k \Big[ \frac{(L_f^k +
L_{\omega})^2}{\gamma_k^2} \Big]$ and notice that
\[ \frac{(L_f^k + L_{\omega})^2}{\gamma_k^2} = \frac{L_f^k}{\gamma_k^2} +
   \frac{2 L_f^k L_{\omega}}{\gamma_k^2} + \frac{L_{\omega}^2}{\gamma_k^2} \]
so that we can bound
\[ \mathbb{E}_{\xi'} \Big[ \tfrac{L_f^k}{\gamma_k^2} \Big]
   =\mathbb{E}_{\xi'} \Big[ \tfrac{(L_f^k)^2}{\max \{ (L_f')^2, \alpha
   \} k^{2 \zeta}} \Big] \leq \Big( \frac{\alpha + \sigma}{\alpha}
   \Big)^2 \frac{1}{k^{2 \zeta}}, \]
where we invoked \tmtextbf{Lemma \ref{lem:3.3}} and take $X = L_f^k, Y = L_f'$, and we
recall that by \ref{D2}, $\mathbb{E} [| L_f^k - L_f' |^2] \leq \sigma^2$.
Similarly, we can deduce that
\[ \mathbb{E}_{\xi'} \Big[ \tfrac{2 L_f^k L_{\omega}}{\gamma_k^2}
   \Big] \leq \mathbb{E}_{\xi'} \Big[ \tfrac{2 L_f^k L_{\omega}}{\max
   \{ L_f', \alpha \}^2} \Big] \leq \mathbb{E}_{\xi'} \Big[ \tfrac{2
   L_f^k}{\alpha L_f'} \Big] \leq \Big( \frac{\alpha + \sigma}{\alpha^2}
   \Big) \frac{2 L_{\omega}}{k^{2 \zeta}} \]

\[ \mathbb{E}_{\xi'} \Big[ \tfrac{L_{\omega}^2}{\gamma_k^2} \Big]
   \leq \tfrac{L_{\omega}^2}{\alpha^2 k^{2 \zeta}} . \]
Putting the bounds together,
\[ \mathbb{E}_{\xi'} \Big[ \tfrac{(L_f^k + L_{\omega})^2}{\gamma_k^2}
   \Big] \leq \tfrac{(\alpha + \sigma + L_{\omega})^2}{\alpha^2 k^{2 \zeta}}
   . \]
Then we have
\[ \mathbb{E}_k [\| x^{k + 1} - x^{\star} \|^2] \leq \| x^k - x^{\star} \|^2 -
   \frac{2}{\gamma_k} [\psi (x^k) - \psi (x^{\star})] + \frac{(\alpha + \sigma
   + L_{\omega})^2}{\alpha^2 k^{2 \zeta}} \]
and we complete the proof by invoking \tmtextbf{Lemma \ref{lem:robbins}}.

\subsubsection{Proof of Lemma \ref{cvx:lem:5}}

We start by bounding $\mathbb{E}_k \Big[ \tfrac{(L_f^k +
L_{\omega})^2}{\gamma_k^2} \Big] \leq \frac{(\alpha + \sigma +
L_{\omega})^2}{\alpha^2 K}$. Telescoping the relation
\[ \mathbb{E}_k [\| x^{k + 1} - x^{\star} \|^2] \leq \| x^k - x^{\star} \|^2 -
   \frac{2}{\gamma_k} [\psi (x^k) - \psi (x^{\star})] + \frac{(\alpha + \sigma
   + L_{\omega})^2}{\alpha^2 K} \]
gives, for all $2 \leq k \leq K$, that
\[ \mathbb{E} [\| x^k - x^{\star} \|^2] \leq \| x^1 - x^{\star} \|^2 +
   \frac{(\alpha + \sigma + L_{\omega})^2}{\alpha^2} \leq \Big( \| x^1 -
   x^{\star} \| + \frac{\alpha + \sigma + L_{\omega}}{\alpha} \Big)^2
   =: \Delta^2 . \]

and $\mathbb{E} [\| x^k - x^{\star} \|] \leq \Delta$. Next consider
\[ f_{x^k} (x^{k + 1}, \xi^k) + \omega (x^{k + 1}) + \frac{\gamma_k}{2} \|
   x^{k + 1} - x^k \|^2 \leq f_{x^k} (x^k, \xi^k) + \omega (x^k) -
   \frac{\gamma_k}{2} \| x^{k + 1} - x^k \|^2 \]
and re-arrangement gives $\gamma_k \| x^{k + 1} - x^k \|^2 \leq (L_f^k +
L_{\omega}) \| x^{k + 1} - x^k \| .$ Dividing both sides by $\| x^{k + 1} -
x^k \|$, we get $\| x^{k + 1} - x^k \| \leq \frac{L_f^k +
L_{\omega}}{\gamma_k} .$ Taking expectation, $\mathbb{E}_k [\| x^{k + 1} - x^k
\|] \leq \frac{\alpha + \sigma + L_{\omega}}{\alpha \sqrt{K}} .$ Then by
Markov's inequality, $\mathbb{P}_{\xi^k, \xi' \sim \Xi} \left\{ \| x^{k + 1} - x^k
\| \leq \frac{2 (\alpha + \sigma + L_{\omega})}{\alpha \sqrt{K}} | \xi^1,
\ldots, \xi^{k - 1} \right\} \geq \frac{1}{2}$. and without loss of
generality, let $\mathsf{Z} = \frac{2 (\alpha + \sigma + L_{\omega})}{\alpha
\sqrt{K}} > 0$. By the same reasoning, we arrive at $\Delta \geq \frac{z}{2}
\cdot \mathbb{P} \{ \| x^k - x^{\star} \| \geq \mathsf{Z} + z \}$ and
dividing both sides by $\frac{z}{2}$ gives the desired tail bound.

\subsubsection{Proof of Theorem \ref{cvx:thm:5}}

Following the same argument as \tmtextbf{Theorem \ref{cvx:thm:4}}, we have, for $\delta \in
(0, 1 / 4)$, that $\mathbb{P} \{ \| x^k - x^{\star} \| \geq \mathsf{Z} +
\delta^{- 1} \Delta \} \leq 2 \delta$. Define $I_k =\mathbb{I} \{ \| x^k -
x^{\star} \| \leq \delta^{- 1} \Delta + \mathsf{Z} \}$. We have, by Markov's
inequality, that $\mathbb{P} \{ \sum_{k = 1}^K I_k \geq K (1 - 2 p^{- 1}
\delta) \} \geq 1 - p$. Then telescoping over \tmtextbf{Lemma \ref{cvx:lem:5}} gives
\begin{align}
  \sum_{k = 1}^K \mathbb{E} \left[ \frac{2}{\gamma_k} (\psi (x^k) - \psi
  (x^{\star})) \right] \leq{} & \| x^1 - x^{\star} \|^2 -\mathbb{E} [\| x^{K + 1}
  - x^{\star} \|^2] + \frac{(\alpha + \sigma + L_{\omega})^2}{\alpha^2} \nonumber\\
  \leq{} & \| x^1 - x^{\star} \|^2 + \Big( \frac{\alpha + \sigma +
  L_{\omega}}{\alpha} \Big)^2 . \nonumber
\end{align}

Conditioned on $\sum_{k = 1}^K I_k \geq K (1 - 2 p^{- 1} \delta)$, we get
\[ \sum_{k = 1}^K \mathbb{E} \left[ \frac{2}{\gamma_k} (\psi (x^k) - \psi
   (x^{\star})) \right] \geq \frac{2 (1 - 2 p^{- 1} \delta) \sqrt{K}}{\alpha +
   \mathsf{G}_{\delta}} \min_{k \in \{ j : I_j = 1 \}} \mathbb{E} [\psi (x^k)
   - \psi (x^{\star})] \]
where $\mathsf{G}_{\delta} \assign \max_x \sup_{\xi \sim \Xi}  \mathsf{Lip} (x, \xi), \| x - x^{\star} \| \leq \delta^{- 1} \Delta +
\mathsf{Z}$. Combining two inequalities completes the proof.

\end{document}